\documentclass{article}
\usepackage{amsthm}
\usepackage{fancyhdr}
\usepackage{amsmath}
\usepackage{indentfirst}
\usepackage{titlesec}
\usepackage{amsfonts}
\usepackage{amssymb}

\newcommand{\cA}{\mathcal A}
\newcommand{\cB}{\mathcal B}

\newcommand{\cF}{\mathcal F}
\newcommand{\cG}{\mathcal G}
\newcommand{\cH}{\mathcal H}
\newcommand{\cL}{\mathcal L}
\newcommand{\cO}{\mathcal O}
\newcommand{\cR}{\mathcal R}

\newcommand{\N}{\mathbb N}
\newcommand{\Q}{\mathbb Q}
\newcommand{\R}{\mathbb R}
\newcommand{\Z}{\mathbb Z}
\newcommand{\cov}{\textrm{\rm Cov}}
\newcommand{\corr}{\textrm{\rm Corr}}
\newcommand{\var}{\textrm{\rm Var}}

\newcommand{\zhfb}{\hfil\break}
\newcommand{\zeps}{\varepsilon}
\newcommand{\bbC}{\mathbb C}
\newcommand{\bbN}{\mathbb N}
\newcommand{\bbP}{\mathbb P}
\newcommand{\bbR}{\mathbb R}
\newcommand{\bbZ}{\mathbb Z}

 \newtheorem{theorem}{Theorem}[section]

 \newtheorem{lemma}[theorem]{Lemma}

 \theoremstyle{definition}
 \newtheorem{definition}[theorem]{Definition}

\newtheorem{remark}[theorem]{Remark}

 \newtheorem{notations}[theorem]{Notations}

 \newtheorem{background}[theorem]{Background}

 \numberwithin{equation}{section}

 \numberwithin{theorem}{section}

\begin{document}

\title{\textbf{On a boundary of the central limit theorem for 
strictly stationary, reversible Markov chains}}

\newcommand{\orcidauthorA}{0000-0000-000-000X} 

\author{Richard C.\ Bradley \\
Department of Mathematics, 
Indiana University, Bloomington, IN 47405, U.S.A.\\ bradleyr@iu.edu}


\maketitle

\begin{abstract}
     Consider the class of (functions of) strictly stationary Markov chains in which 
(i) the second moments are finite and (ii)  absolute regularity ($\beta$-mixing) is
satisfied with exponential mixing rate.
For (functions of) Markov chains in that class that are also reversible, 
the central limit theorem holds,
as a well known byproduct of results of Roberts, Rosenthal, and 
Tweedie in two papers in 1997 and 2001 involving reversible Markov chains.
In contrast, for (functions of) Markov chains in that class 
that are not reversible, the central limit theorem may fail to hold, as is known from counterexamples, including ones with arbitrarily fast
mixing rate (for absolute regularity).
Here it will be shown that for Markov chains in that class that are reversible, the
``borderline'' class of mixing rates (for absolute regularity) for the central limit theorem 
is in fact exponential.
That will be shown here with a class of counterexamples — strictly stationary,
 countable-state Markov chains that are reversible, have finite second moments, 
 and satisfy absolute regularity with mixing rates that can be arbitrarily close to 
 (but not quite) exponential, but fail to satisfy the central limit theorem.
 \end{abstract}

\textbf{Keywords:}\ \ strictly stationary, reversible Markov chain;
\hfil\break 
\null \hskip 1.2 in 
strong mixing; absolute regularity; central limit question
\medskip

\textbf{AMS 2020 Mathematics Subject Classifications:}\ \ 
60J10, 60G10




\section{Introduction}
\label{sc1}

   In this paper, for simplicity, every Markov chain
(i) is indexed by the set $\bbZ$ of all integers,
(ii) is strictly stationary, and
(outside of Section \ref{sc8})
(iii) has as its state space either the set $\bbR$ of all real numbers 
(accompanied by its Borel $\sigma$-field
$\cR$) or a finite or countably infinite subset of $\bbR$.
\smallskip

   Over twenty years ago, in the papers of 
Roberts and Rosenthal \cite {ref-journal-RR} and 
Roberts and Tweedie \cite {ref-journal-RT} together,
some special theory for reversible Markov chains was developed which included
 (among many other things) the following ``key result'' (in much greater 
 generality than the statement given here):    
For a given (say real) strictly stationary Markov chain 
$(X_k, k \in \bbZ)$ that is reversible and is also ``irreducible''
(in a certain sense equivalent to Harris recurrence), the
following two conditions are equivalent: \zhfb
(a) geometric ergodicity; \zhfb
(b) a certain ``$\cL^2$ spectral gap'' condition.
\smallskip

     Of course the term ``reversible'' simply means that
the distribution of the Markov chain as a whole is invariant
under a reversal of the ``direction of time''.   
(More on that in Section \ref {sc2}.)\ \
The definitions of the three conditions (a), (b), and ``irreducible'' will not 
be given here directly; instead,  Section \ref {sc2} will include a well 
known equivalent formulation for each one in terms of 
the ``dependence coefficients'' associated with certain mixing conditions.
In particular, we mention here (as well as in Section \ref {sc2}) that the 
``geometric ergodicity'' condition is equivalent to absolute 
regularity ($\beta$-mixing) with mixing rate at least exponentially fast.
Such use of mixing conditions and their dependence coefficients
in the treatment of the aforementioned ``key result'', has long 
been a common practice; see e.g.\
\cite{ref-journal-KM}, \cite{ref-journal-Longla14},
\cite{ref-journal-LP}, or the primarily expository paper
\cite{ref-journal-Bradley2021}.
\vskip 0.1 in
 
     The aforementioned ``key result'' has the following well known consequence:
For (functions of) strictly stationary Markov chains that satisfy geometric 
ergodicity --- that is, absolute regularity with at least exponentially fast 
mixing rate --- the extra property of reversibility, if satisfied, provides
considerable extra leverage in the development of central limit theory.
(That is reviewed in a bit of detail in Section \ref {sc2}.)\ \   
A natural question is the following:
For (functions of) strictly stationary Markov chains that satisfy absolute 
regularity but with a mixing rate that is slower than exponential, 
does the extra property of reversibility, if satisfied, still provide considerable 
extra leverage in the development of central limit theory?
\smallskip

      For mixing rates of power-type $O(n^{-p})$ for $p > 1$, 
the answer to that question seems to be ``essentially no''.
That was illustrated with a class of Markov chains constructed by
Doukhan et.\ al.\ [\cite {ref-journal-DMR1994},
Theorems 2 and 5 and Corollary 1] that were in a certain 
sense as sharp as possible.
(Those examples are reversible, though that property does not appear to be 
mentioned explicitly in that paper.)
\smallskip

     Here we shall examine that question for mixing rates that are close to
exponential.  
Our main result, Theorem \ref{thm3.3} in Section \ref{sc3}, will describe a class
of counterexamples that will show that under reversibility, 
finite second moments, and absolute regularity with mixing rates arbitrarily 
close to (but short of) exponential, the central limit theorem may fail to hold. 
\smallskip

   The rest of this paper is organized as follows:
Section \ref{sc2} will provide relevant background information on certain 
mixing conditions, with special emphasis on Markov chains.
In Section \ref{sc3}, some key terminology will be developed, and the
main result (Theorem \ref{thm3.3}) will be stated.
Section \ref{sc4} will provide some background information on a partial
limiting distribution for (normalized) partial sums in the construction
for the main result. 
Section \ref{sc5} will develop a class of building blocks — strictly stationary,
reversible, 3-state Markov chains of a relatively simple structure — that will be
used in the construction for the main result.
Section \ref{sc6} will lay out some miscellaneous further preliminary material
that will play a role in the proof of the main result.
Then Section \ref{sc7} will give the proof of the main result.
Finally, Section \ref{sc8}, an Appendix, will treat a topic that provides a little 
additional background motivation for the main result.

\section{Mixing conditions and Markov chains}
\label{sc2}

   In this section, we shall lay out some notations and review some 
background facts regarding mixing conditions and Markov chains.
Part of the exposition in this section (essentially all of the part from here 
through Remark \ref{rem2.4} below) has been taken practically 
verbatim from the introduction of a recent paper by the 
author \cite {ref-journal-Bradley2023}.

\begin{notations}
\label{nt2.1}
First, here are a few really basic notations:\hfil\break
$\bbN$ denotes the set of all positive integers, \hfil\break
$\bbZ$ denotes the set of all integers, \hfil\break
$\bbR$ denotes the set of all real numbers, \hfil\break
$\cR$ denotes the Borel $\sigma$-field on $\bbR$, and \hfil\break
$\bbC$ denotes the set of all complex numbers.
\hfil\break
The usual notations such as $\bbR^\bbZ$ and $\cR^\bbZ$ will be used for 
Cartesian products of sets and for product $\sigma$-fields.
\smallskip

     We use an old-fashioned convention that for sets $A$ and $B$, the notation
$A \subset B$ simply means that any element of $A$ is an element of $B$; it
includes the case $A=B$.
\smallskip 
   
     In what follows, the setting will be a probability space
$(\Omega, \cF, P)$, rich enough to accommodate all random variables declared.
In this paper, all random variables will be real-valued, possibly discrete.
The $\sigma$-field (on $\Omega$) generated by a family $(Y_i, i \in I)$ of  
random variables (with the index set $I$ being nonempty), will be denoted
$\sigma(Y_i, i \in I)$.
\smallskip

The phrase ``central limit theorem'' will sometimes be abbreviated as CLT.
\end{notations}

\begin{notations}
\label{nt2.2}
Suppose $\cA$ and $\cB$ are any 
two $\sigma$-fields $\subset \cF$.
Define the following three measures of dependence:
First,
\begin{equation}
\label{eq2.21}
\alpha(\cA,\cB) := \sup_{A\in \cA, B\in \cB} |P(A\cap B)-P(A) P(B)|.
\end{equation}
Next, 
\begin{equation}\label{eq2.22}
\beta(\cA,\cB):= 
\sup \frac{1}{2} \sum^I_{i=1}\sum^J_{j=1} |P(A_i\cap B_j) - P(A_i)P(B_j)|
\end{equation}
where the supremum is taken over all pairs of finite partitions 
$\{A_1,A_2,\dots,A_I\}$ and $\{B_1,B_2,\dots,B_J\}$ of 
$\Omega$ such that $A_i\in \cA$ for each $i$ and $B_j\in \cB$ for each~$j$.
[The factor of $1/2$ in (\ref{eq2.22}) is not of special significance, 
but has become customary in order to make 
certain inequalities a little nicer.]\ \ 
Finally, define the ``maximal correlation coefficient''
\begin{equation}
\label{eq2.23}
\rho(\cA,\cB) := \sup | \corr(Y,Z)|
\end{equation}
where the supremum is taken over all pairs of square-integrable 
random variables $Y$ and $Z$ such that $Y$ is $\cA$-measurable 
and $Z$ is $\cB$-measurable.
\medskip 

The following inequalities are elementary and well known:
\begin{align}
\label{eq2.24}
0 &\le 2\alpha(\cA,\cB) \le \beta(\cA,\cB) \le 1; \quad {\rm and} \\
\label{eq2.25}
0 &\le 4\alpha(\cA,\cB) \le \rho(\cA,\cB) \le 1. 
\end{align}
(See e.g.\  [\cite{ref-journal-Bradley2007}, Vol.\ 1, Proposition 3.11].)\ \ 
The quantities $\alpha(\cA,\cB)$, $\beta(\cA,\cB)$, and $\rho(\cA,\cB)$ are 
all equal to $0$ if the $\sigma$-fields $\cA$ and $\cB$ are independent, 
and are all positive otherwise.
\end{notations}

\begin{notations}
\label{nt2.3}
Suppose $X:= (X_k$, $k\in \Z)$ is a (not necessarily Markovian) 
strictly stationary sequence of (real-valued) random variables.
Refer to Notations \ref{nt2.1}.
For each integer $j$, define the notations $\cF^j_{-\infty} :=\sigma(X_k, k\le j)$ and 
$\cF^\infty_j:= \sigma(X_k$, $k\ge j)$.
\medskip

     For each positive integer $n$, define the following three ``dependence coefficients'':
\begin{align}
\label{eq2.31}
\alpha(n) &= \alpha_X(n) := \alpha(\cF^0_{-\infty}, \cF^\infty_n); \\
\label{eq2.32}
\beta(n) &= \beta_X(n) := \beta(\cF^0_{-\infty}, \cF^\infty_n);\,\,\mbox{and}\\
\label{eq2.33}
\rho(n) &= \rho_X(n) := \rho(\cF^0_{-\infty}, \cF^\infty_n).
\end{align}
For each positive integer $n$, one has by strict stationarity that 
$\alpha(n) = \alpha(\cF^j_{-\infty},\cF^\infty_{j+n})$ for every integer 
$j$, and the analogous comment holds for $\beta(n)$ and 
for $\rho(n)$ as well.

     Also, one (trivially) has that $\alpha(1) \ge \alpha(2)\ge \alpha(3)\ge \dots\,$; and the analogous comment holds for the numbers $\beta(n)$ and for the numbers $\rho(n)$.
\medskip

      The strictly stationary sequence $X$ is said to satisfy

\noindent ``strong mixing'' (or ``$\alpha$-mixing'') if $\alpha_X(n) \to 0$ as $n\to\infty$;

\noindent ``absolute regularity'' (or ``$\beta$-mixing'') if $\beta_X(n) \to 0$ as $n\to\infty$;

\noindent ``$\rho$-mixing'' if $\rho_X(n) \to 0$ as $n\to \infty$.
\medskip

The strong mixing ($\alpha$-mixing) condition is due to 
Rosenblatt \cite{ref-journal-Rosenblatt1956}.
The absolute regularity ($\beta$-mixing) condition was first studied by 
Volkonskii and Rozanov \cite{ref-journal-VolkRoz}, and was attributed there 
to Kolmogorov.
The $\rho$-mixing condition is due to 
Kolmogorov and Rozanov \cite{ref-journal-KolRoz}.
(The maximal correlation coefficient $\rho(\cA,\cB)$ itself, for 
 $\sigma$-fields $\cA$ and $\cB$, was first studied earlier by 
 Hirschfeld \cite{ref-journal-Hirschfeld} in a statistical context that had no 
 particular connection with stochastic processes.)
\medskip 

    By (\ref{eq2.24}) and (\ref{eq2.25}), one has that for each positive integer $n$,
\begin{align}
\label{eq2.34}
0&\le 2\alpha(n) \le \beta(n) \le 1,\,\,\mbox{and}\\
\label{eq2.35}
0 &\le 4\alpha(n) \le \rho(n) \le 1.
\end{align}
By (\ref{eq2.34}), absolute regularity ($\beta$-mixing) implies strong mixing 
($\alpha$-mixing); and 
by (\ref{eq2.35}), $\rho$-mixing implies strong mixing ($\alpha$-mixing).    
\end{notations}

\begin{remark}
\label{rem2.4}
Now suppose that 
$X:=(X_k, k\in \Z)$ is a strictly stationary {\em Markov chain} (with the random 
variables $X_k$, $k \in \Z$ being real-valued, possibly discrete).
(No assumption of ``reversibility'' yet.)
\medskip

      As a well known consequence of the Markov property, for each positive 
integer $n$, eqs.\ (\ref{eq2.31})-(\ref{eq2.33}) hold in the following augmented 
forms for the given (strictly stationary) Markov chain $X$:
\begin{align}
\label{eq2.41}
\alpha(n) &= \alpha_X(n) =  
\alpha(\cF^0_{-\infty}, \cF^\infty_n) = \alpha\bigl(\sigma(X_0),\sigma(X_n)\bigl);\\
\label{eq2.42}
\beta(n) &= \beta_X(n) =  
\beta(\cF^0_{-\infty},\cF^\infty_n) = \beta\bigl(\sigma(X_0), \sigma(X_n)\bigl);\\
\label{eq2.43}
\rho(n) &= \rho_X(n) = 
\rho(\cF^0_{-\infty}, \cF^\infty_n) = \rho\bigl(\sigma(X_0), \sigma(X_n)\bigl).
\end{align}
(See e.g.\  [\cite{ref-journal-Bradley2007}, Vol.\ 1, Theorem 7.3].)
\medskip

      By strict stationarity and (\ref{eq2.41})-(\ref{eq2.43}), one has that for 
any integer $j$ and any positive integer $n$, the (strictly stationary) 
Markov chain $X$ satisfies \hfil\break 
(i)\ $\alpha(n) = \alpha(\sigma(X_j),\sigma(X_{j+n}))$, \hfil\break
(ii)\ $\beta(n) = \beta(\sigma(X_j), \sigma(X_{j+n}))$, and \hfil\break 
(iii) $\rho(n) = \rho(\sigma(X_j), \sigma(X_{j+n}))$. \hfil\break
 
      Here are some special facts involving the 
dependence coefficients $\rho(n)$, $n \in \N$. 
As an elementary consequence of (\ref{eq2.23}), for any two 
$\sigma$-fields $\cA$ and $\cB$, 
$\rho(\cA, \cB) = \sup \|E(Y|\cB)\|_2/\|Y\|_2$ 
where the supremum is taken over all square-integrable, $\cA$-measurable
random variables $Y$ with mean 0.  
(When necessary, interpret 0/0 := 0.)\ \ 
As a well known application of that fact
(together with the equality in (iii) in the preceding paragraph),
for the given strictly stationary Markov chain
$X:=(X_k, k\in \Z)$, one has that 
for any pair of positive integers $m$ and $n$, 
\begin{equation}
\label{eq2.44}
\rho(m+n) \le \rho(m) \cdot \rho(n).
\end{equation}
As a well known consequence, for the given strictly stationary Markov chain $X$,
if $\rho(m) < 1$ for some $m \in \N$, then the Markov chain $X$ is $\rho$-mixing,
with $\rho(n) \to 0$ at least exponentially fast as $n \to \infty$.
In Remark \ref{rem2.5}, Remark \ref{rem2.9}, and (indirectly) 
Remark \ref{rem2.10} below, the case $\rho(1) < 1$ will come into play. 
\end{remark}

\begin{remark}
\label{rem2.5}
For any given (say real) strictly stationary Markov chain
$(X_k, k \in \bbZ)$, the following statements are well known:
\smallskip

   (A) From the theory developed by 
Nummelin and Tweedie \cite {ref-journal-NTweed} and
Nummelin and Tuominen \cite {ref-journal-NTuo},
one has (among many other results) that the following two 
conditions are equivalent: \zhfb
(w) geometric ergodicity; \zhfb
(x) absolute regularity ($\beta$-mixing) with 
$\beta(n) \to 0$ at least exponentially fast as $n \to \infty$. 
\smallskip

   (B) From some basic functional analysis, one has that the
following two conditions are equivalent: \zhfb
(y) the ``$\cL^2$ spectral gap'' condition alluded to in Section {\sc1}; 
\zhfb
(z) $\rho(1) < 1$.
\smallskip

   (C) The following three conditions are equivalent: \zhfb
(i) the ``irreducibility'' condition alluded to in Section \ref {sc1};
\zhfb
(ii) Harris recurrence; \zhfb
(iii) ergodicity and $\lim_{n \to \infty} \beta(n) < 1$.
\smallskip

\noindent If those three equivalent conditions (i), (ii), (iii) are satisfied, 
then for some positive integer $p$, one has that
$\lim_{n \to \infty} \beta(n) = 1 - (1/p)$ and the Markov chain
has period $p$.
(If $p =1$, then the Markov chain is aperiodic and satisfies 
absolute regularity.)\ \ 
\smallskip

   More details on comments (A), (B), and  (C) can be found 
e.g.\ in the (mainly) expository paper 
\cite {ref-journal-Bradley2021}, 
and (for (A) and (C)) in 
[\cite {ref-journal-Bradley2007}, Vol. 2, Chapter 21]. 
Much further relevant information can be found in references such as
the books by Orey \cite {ref-journal-Orey} and by  
Meyn and Tweedie \cite {ref-journal-MT}.     
\end{remark}

\begin{remark}
\label{rem2.6}
There is a considerable history of the use of (functions of) strictly stationary 
Markov chains to illuminate the ``borderline'' of central limit theorems 
under weak dependence.
\smallskip

   (A) For example, consider the following central limit theorem of 
Ibragimov \cite{ref-journal-Ibragimov1962},
stated here as in Theorem 18.5.4 in the book by Ibragimov and Linnik
\cite{ref-journal-IbrLin}:
\smallskip

\noindent {\it Suppose $(X_k, k \in \Z)$ is a (not necessarily Markovian) strictly 
stationary sequence of {\rm bounded} random variables with mean 0, such that
$\sum_{n=1}^\infty \alpha(n) < \infty$.  
Then $\sigma^2 := E(X_0^2) + 2 \sum_{n=1}^\infty E(X_0X_n)$ exists in 
$[0, \infty)$, the sum being absolutely convergent; and if $\sigma^2 > 0$ then 
$(1/(n^{1/2} \sigma))\sum_{k=1}^n X_k$ converges in distribution to $N(0,1)$
as $n \to \infty$.}
\smallskip

   (B) Davydov [\cite{ref-journal-Davydov1973}, Example 2] constructed a class of 
examples of bounded, strictly stationary, countable-state Markov chains, for which
the partial sums, suitably normalized, converge in distribution to a symmetric 
stable law with index (exponent) strictly between 1 and 2.
By the Theorem of Types (see e.g.\ Theorem 14.2 in \cite{ref-journal-Bill}), 
in those examples, regardless of what normalization is employed, 
normalized partial sums cannot converge in distribution to a nondegenerate
normal law.  
As Davydov \cite{ref-journal-Davydov1973} pointed out, for any given $\zeps > 0$, 
those examples include ones for which [see (\ref{eq2.34})] 
$\alpha(n) \leq \beta(n) = O(n^{-1 + \zeps})$ as $n \to \infty$.
[Davydov's discussion there focused on the dependence coefficient
$\alpha(n)$; however, with calculations on pp.\ 327-328 in that paper
\cite{ref-journal-Davydov1973}, he implicitly verified the same mixing rate 
for $\beta(n)$ as for $\alpha(n)$ (modulo a constant factor).]\ \ 
With those examples, Davydov showed that the central limit theorem 
of Ibragimov stated in (A) above is very sharp.
\smallskip

   (C) Ibragimov \cite{ref-journal-Ibragimov1962} also proved a CLT for 
(not necessarily Markovian) strictly stationary {\it un\/}bounded sequences of 
random variables having certain finite moments of higher than 
second order, with (as compensation) certain mixing rates for $\alpha(n)$ faster 
than the summable rate in the CLT (for bounded sequences) stated in (A) above.
In essentially the same way as in (B) above, 
Davydov [\cite{ref-journal-Davydov1973}, Example 1] showed with a class of 
(unbounded) strictly stationary, countable-state Markov chains that those 
CLTs are very sharp.
[The counterexamples alluded to here and in (B) above 
were structured along the lines of similar Markov chains that were studied by 
Chung \cite{ref-journal-Chung1967} without emphasis on mixing conditions.]
\smallskip

   (D) For (not necessarily Markovian) strictly stationary sequences of (say) 
unbounded random variables, 
Doukhan et.\ al.\ \cite{ref-journal-DMR1994}   
proved a very sharp central limit theorem (in fact a functional CLT), involving a 
trade-off between (i) quantile functions of the marginal distribution and
(ii) the mixing rate for $\alpha(n)$.
In Theorems 2 and 5 and Corollary 1 of that paper \cite{ref-journal-DMR1994}, 
they also constructed a class of (functions of) strictly stationary Markov chains 
that showed that with respect to both $\alpha(n)$ and $\beta(n)$, at least for 
power type mixing rates $O(n^{-p})$ for $p > 1$, their CLT is as sharp as possible  
(in a manner somewhat different from the sharpness in the examples of 
Davydov discussed in (B) and (C) above).
\end{remark}

\begin{remark}
\label{rem2.7}
This paper here will focus on (functions of) strictly stationary Markov chains
with an assumption of finite second moments (and, for convenience, mean 0),
but with no assumption of any other condition involving (generalized) moments 
or quantiles.
\smallskip

There exist strictly stationary, countable-state Markov chains
$X := (X_k, k \in \bbZ)$ and functions $h: \bbR \to \bbR$ satisfying
$E[\, [h(X_0)]^2] < \infty$ and $\var[h(X_0)] > 0$,
such that $\beta_X(n) \to 0$ at least exponentially fast
as $n \to \infty$, but no normalized versions of the partial sums 
$\sum_{k=1}^n h(X_k)$ converge in distribution to $N(0,1)$ as $n \to \infty$.
\smallskip

That was shown with constructions by the author 
\cite {ref-journal-Bradley1983}, 
by Herrndorf \cite {ref-journal-Herrndorf}, and 
by H\"aggstrom \cite {ref-journal-Haggstrom}.
In the examples in \cite {ref-journal-Bradley1983} and
\cite {ref-journal-Herrndorf}, the mixing rate for absolute regularity
can be made arbitrarily fast (short of $m$-dependence).
In the example in \cite {ref-journal-Bradley1983},
(i) the partial sums, suitably normalized, converge in distribution along
subsequences to all infinitely divisible laws, and 
(ii) (with a suitable choice of certain parameters) the variances of those 
partial sums diverges to $\infty$ almost (but not quite) as fast as $n^2$. 
For more on those examples (especially regarding hidden properties of the
example in \cite {ref-journal-Herrndorf}), see 
[\cite {ref-journal-Bradley2007}, Vol.\ 3, Chapter 31].
Section 8 of this paper here, an Appendix, gives a detailed, gentle exposition of one 
particular method of a conversion of the example in \cite{ref-journal-Bradley1983} 
to a Markov chain (not just a ``function'' of a Markov chain), without
weakening the other stated properties.     
\end{remark}

\begin{notations}
\label{nt2.8}
Now we turn to strictly stationary, {\it reversible\/} Markov chains.
\smallskip

   (A) A given (say real) strictly stationary Markov chain 
$X := (X_k, k \in \Z)$
is said to be ``reversible'' if the distribution 
[on $(\bbR^\Z, \cR^\Z)$] of the ``time-reversed''
Markov chain $(X_{-k}, k \in \Z)$ is the same as that
of the sequence $X$ itself.
\smallskip

   (B) By a standard argument, a given (say real) strictly stationary Markov 
chain $X := (X_k, k \in \Z)$ is reversible if and only if
the distribution [on $(\bbR^2, \cR^2)$] of the random vector
$(X_1, X_0)$ is the same as that of the random vector 
$(X_0, X_1)$.  
\smallskip

   (C) The counterexamples constructed by 
Doukhan et.\ al.\ \cite{ref-journal-DMR1994} 
alluded to in Remark \ref{rem2.6}(D) above are reversible 
(although that property does not seem to be mentioned in that paper).
The other counterexamples alluded to in Remarks \ref{rem2.6}
and \ref{rem2.7} are not reversible.
\end{notations}

\begin{remark}
\label{rem2.9}
(A) Section 1 alluded to a ``key theorem'' that came from the theory developed 
in \cite {ref-journal-RR} and \cite {ref-journal-RT} for reversible Markov chains.
As is well known, in light of Remark \ref {rem2.5}(A)(B)(C), that ``key theorem''
— or at least a restricted version of it  — 
can be formulated in terms of the dependence coefficients associated 
with the absolute regularity and $\rho$-mixing conditions, 
as follows (here with some redundancy):
\smallskip

\noindent {\em If $X := (X_k, k \in \bbZ)$ is a strictly stationary
Markov chain which is reversible and also satisfies 
ergodicity and $\lim_{n \to \infty} \beta(n) < 1$, 
then the following two conditions are
equivalent: \zhfb
(i) absolute regularity with $\beta(n) \to 0$ at least exponentially fast 
as $n \to \infty$; \zhfb
(ii) $\rho(1) < 1$.}
\medskip

    (B)    With calculations that seem to be ``close cousins'' of ones in 
\cite {ref-journal-RR} and \cite {ref-journal-RT}, it was shown in 
\cite{ref-journal-Bradley2021} that
for strictly stationary Markov chains that are reversible
(but do not necessarily satisfy absolute regularity or
related properties such as Harris recurrence),
the condition $\rho(1) < 1$ is equivalent to strong mixing with
$\alpha(n) \to 0$ at least exponentially fast as $n \to \infty$.   
\end{remark}

\begin{remark}
\label{rem2.10}
It is well known that the ``key theorem'' mentioned in Remark \ref{rem2.9}, 
and extensions of it, provide some special leverage in the development of 
central limit theory for strictly stationary, reversible Markov chains.
Compare Remark \ref{rem2.7}, involving (functions of) strictly stationary
non-reversible Markov chains, with (A) and (B) below,
involving (functions of) strictly stationary reversible Markov chains.
Comments (A) and (B) below were discussed (in greater detail, among much else) 
by Michael Lin at a conference at the University of Cincinnati in May 2024, 
in a talk on joint work of his with Christophe Cuny. 
\smallskip

     (A) Suppose one combines the ``key theorem'' from 
\cite {ref-journal-RR} and \cite {ref-journal-RT}
that was given in Remark \ref{rem2.9}(A) above, 
with a classic central limit theorem for (functions of) strictly stationary 
$\rho$-mixing Markov chains (such as in  
\cite {ref-journal-Rosenblatt1970},
\cite {ref-journal-Rosenblatt1971}, or 
\cite {ref-journal-Lifshits}),
or with an appropriate CLT for strictly stationary, not necessarily Markovian
random sequences, (such as in 
[\cite{ref-journal-Ibragimov1975}, Theorem 2.2]).
It is well known that one thereby obtains the following result: 
\smallskip

\noindent {\it Suppose $X := (X_k, k \in \Z)$ is a strictly stationary Markov chain
that is reversible, and $h: \R \to \R$ is a Borel function such that  
$E[[h(X_0)]^2] < \infty$ and $E[h(X_0)] = 0$.
Suppose this Markov chain $X$ satisfies absolute regularity
 with $\beta_X(n) \to 0$ at least exponentially fast as $n \to \infty$.
Then 
$\sigma^2 := E[[h(X_0)]^2] \allowbreak 
+ 2 \sum_{n=1}^\infty E[h(X_0)\, h(X_n)]$ exists in 
$[0, \infty)$, the sum being absolutely convergent; and if $\sigma^2 > 0$ then \break 
$(1/(n^{1/2} \sigma))\sum_{k=1}^n h(X_k)$ converges in distribution to $N(0,1)$
as $n \to \infty$.} 
\smallskip

    (B) By combining the comment in Remark \ref{rem2.9}(B) with references such 
as the ones mentioned here in (A) above, one obtains this CLT under
the mixing assumption of strong mixing (instead of absolute regularity)    
with $\alpha(n) \to 0$ at least exponentially fast as $n \to \infty$. 
\end{remark}

   The purpose of this paper here is to show that for (functions of) strictly stationary,
reversible Markov chains in the CLT reviewed in Remark \ref{rem2.10}(A) above,
in essence the (exponential) mixing rate cannot be slowed down at all
(without strengthening or adding other assumptions). 
For any proposed mixing rate (for absolute regularity or for strong mixing) that is 
slower than exponential, there exists a (nontrivial) counterexample 
whose mixing rate is at least as fast.
That result, the main result of this paper,
will be formulated in Theorem \ref{thm3.3} in the next section.

\section{Formulation of the main result}
\label{sc3}

   The main result of this paper will be formulated in Theorem \ref{thm3.3} below,
with the use of material presented in Remark \ref{rem3.1}
and Notations \ref{nt3.2}.

\begin{remark}
\label{rem3.1}
The new examples, as described in Theorem \ref{thm3.3} below,
will be strictly stationary random sequences 
$(X_k, k \in \bbZ)$ (reversible Markov chains) 
that satisfy the conditions
\begin{equation}
\label{eq3.11}
E(X_0^2) < \infty \indent {\rm and} \indent EX_0 = 0,
\end{equation} 
as well as the condition
\begin{equation}
\label{eq3.12}
\lim_{n \to \infty} n^{-1}
\var\biggl(\, \sum_{k=1}^nX_k\biggl)\ 
=\ \infty
\end{equation}
and even the (stronger) condition
\begin{equation}
\label{eq3.13}
{\rm for\ every}\ b > 0, \quad 
\lim_{n \to \infty} \biggl[\, 
\sup_{r \in \bbR}\,
P\biggl(r-b\ <\ n^{-1/2}\sum_{k=1}^nX_k\ <\ r+ b\biggl) \biggl]\ =\ 0.
\end{equation}
Of course (\ref {eq3.13}) [in conjunction with 
(\ref {eq3.11})] implies (\ref {eq3.12}) by
Chebyshev's inequality.
Eq.\ (\ref {eq3.13}) might be referred to as the ``complete dissipation'' of the 
random variables $n^{-1/2}\sum_{k=1}^nX_k$ as $n \to \infty$.  
For a given fixed positive integer $n$, the term in the brackets 
in (\ref {eq3.13}) is, as a function of $b > 0$, a well known 
``concentration function'' of the random variable $n^{-1/2}\sum_{k=1}^nX_k$.
\smallskip

   Of course, by a trivial argument, (\ref{eq3.13}) is equivalent to the 
(at first sight seemingly weaker) condition that there {\it exists\/} a 
positive number $b$ such that the equality in (\ref{eq3.13}) holds.
In particular, to verify (\ref{eq3.13}), it suffices to carry out 
the argument for (say) just the case $b = 1$, allowing one 
to slightly simplify the arithmetic.     
\end{remark}

\begin{notations}
\label{nt3.2}
In the examples constructed in this paper, a certain non-degenerate, 
non-normal distribution, labeled below as ``$\mu_{\rm P1sL}$'', will 
play a role as a ``partial limiting distribution'' — a limit law for
convergence in distribution along subsequences. 
That distribution is formulated in part (B) below, employing brief 
preliminary material in part (A).
\medskip 

   (A)\ The Laplace (``double exponential'') distribution is the
absolutely continuous distribution (probability measure on
$(\R, \cR)$) with probability density function 
$f_{\rm Laplace}$ given by 
$f_{\rm Laplace}(x) = (1/2) e^{-|x|}$ for $x \in \bbR$.
It is well known (see e.g.\ the chart 
on page 348 of  \cite{ref-journal-Bill}) that its characteristic 
function $\varphi_{\rm Laplace}$ is given
by $\varphi_{\rm Laplace}(t) = 1/(1 + t^2)$ for $t \in \bbR$.
\medskip

   (B)\ Let $\mu_{\rm P1sL}$ denote the distribution [probability 
measure on $(\bbR, \cR)$] of a random variable $Y$ of the form
\begin{equation}
\label{eq3.21}
Y(\omega)\ :=\ \sum_{k=1}^{N(\omega)} \eta_k(\omega)
\indent {\rm for}\ \omega \in \Omega 
\end{equation}       
where (i) $(\eta_1, \eta_2, \eta_3, \dots)$ is a sequence of independent, 
identically distributed random variables
with the Laplace distribution (see (A) above), and 
(ii) the random variable $N$ is independent of the sequence
$(\eta_1, \eta_2, \eta_3, \dots)$ and is nonnegative 
integer-valued and has the Poisson distribution with mean 1.
\smallskip

   In the subscript of the notation $\mu_{\rm P1sL}$, 
the symbols P, 1, s, and L stand for 
``{\bf P}oisson with mean {\bf 1}'' 
mixture of ``{\bf s}ums of {\bf L}aplace'' random variables.
In (\ref{eq3.21}), it is tacitly understood that for 
$\omega \in \Omega$ such that $N(\omega) = 0$, 
the value $Y(\omega)$ is defined to be 0, with the 
``empty sum'' in the right hand side of that equation
understood to take the value 0.
\medskip

   (C) In the class of examples constructed in this paper
for Theorem \ref{thm3.3}, the
partial sums, suitably normalized, will converge 
in distribution along subsequences to the law 
$\mu_{\rm P1sL}$.
To help verify that in arguments later on, the 
following quick calculation will be useful for later reference:
\smallskip  
   
   The characteristic function $\varphi_{\rm P1sL}$ of the distribution 
$\mu_{\rm P1sL}$ can be seen from paragraphs (A) and (B) above, 
as follows:  
For $t \in \bbR$,
\begin{align} 
\label{eq3.22}
\varphi_{\rm P1sL}(t)\ &=\ E[e^{itY}]\ =
\sum_{n=0}^\infty 
\biggl( P(N=n) \cdot E\Bigl[e^{itY}\Bigl|N=n\Bigl]\, \biggl)\
 =\ \sum_{n=0}^\infty
\biggl( \frac{1}{n!} e^{-1} \cdot 
E\biggl[\exp\biggl( it \sum_{k=1}^n \eta_k \biggl) \biggl | 
N=n\biggl]\, \biggl)   \nonumber \\
&=\ 
\sum_{n=0}^\infty  \biggl( \frac{1}{n!} e^{-1} \cdot 
E\biggl[\exp\biggl( it \sum_{k=1}^n \eta_k \biggl) \biggl]
\, \biggl)\   
=\ \sum_{n=0}^\infty \biggl[\, \frac{1}{n!} e^{-1} \cdot  
\Bigl( E\bigl[\exp(it\eta_1)\bigl]\, \Bigl)^n \biggl] \nonumber \\  
&=\ e^{-1} \sum_{n=0}^\infty 
\biggl[\, \frac{1}{n!} \cdot \biggl(\,  
\frac {1} {1 + t^2} \biggl)^n\, \biggl]\  
=\ \exp \biggl( \frac {1} {1 + t^2}\ -\ 1 \biggl)\ . \\
\nonumber   
\end{align}
\end{notations}

   Now here is our main result:

\begin{theorem}
\label{thm3.3}
Suppose $(\zeta_1, \zeta_2, \zeta_3 \dots)$ is a sequence
of positive numbers such that \zhfb
(a) $\lim_{n \to \infty} \zeta_n = 0$, and \zhfb
(b) for every $c > 0$, 
$\lim_{n \to \infty} [\, \zeta_n/e^{-cn}] = \infty$. 

   Then there exists a (real) strictly stationary, countable-state
Markov chain $X := (X_k, k \in \bbZ)$ which is reversible, such that 
the following four conditions hold: \zhfb
(i) eq.\ (\ref{eq3.11}) holds: $E(X_0^2) < \infty$ and $EX_0 = 0$; \zhfb
(ii) $\alpha_X(n)\ \leq\ \beta_X(n)\  \leq\ \zeta_n$ for\ every $n \in \bbN$; \zhfb
(iii) eqs.\ (\ref{eq3.12}) and (\ref{eq3.13}) hold; and \zhfb
(iv) there exists a strictly increasing sequence 
$\bigl(J(1), J(2), J(3), \dots\bigl)$ of positive integers, 
and a sequence of positive constants
$(b_1, b_2, b_3, \dots)$ satisfying $b_n \to \infty$
as $n \to \infty$, such that
\begin{equation}
\label{eq3.31}
b_n^{-1} \sum_{k=1}^{J(n)}X_k\ \to\ \mu_{\rm P1sL}\ \
{\rm in\ distribution\ as}\ n \to \infty. 
\end{equation}    
\end{theorem}

\begin{remark}
\label{rem3.4}
Of course assumption (b) is that the convergence of the positive numbers
$\zeta_n$ to 0 is at a rate that is slower than any exponential decay.
\smallskip

By (\ref{eq3.31}) and the Theorem of Types 
(see e.g.\ Theorem 14.2 in \cite{ref-journal-Bill} again), 
for the Markov chain $X$ in this theorem, 
regardless of what normalization is employed, 
the normalized partial sums cannot converge in distribution to a 
nondegenerate normal law.
\smallskip 

The pertinence of property (iii) in the conclusion of this theorem may warrant 
a comment.
In the study of CLTs for sequences $(X_k, k \in \Z)$ of dependent random variables
(say strictly stationary, with mean 0 and typically with finite second moments),
there is sometimes an interest in whether $n^{-1/2} (X_1 + X_2 + \dots + X_n)$
— that is, with ``conventional'' normalizing constants $n^{-1/2}$ — 
converges to 0 in probability as $n \to \infty$ 
[which is equivalent to convergence in distribution 
to the ``degenerate normal law'' $N(0,0)$ (the point mass at 0)].
In the theorem here, such behavior is ruled out by property (iii) --- specifically, by
eq.\ (\ref{eq3.13}). 
\end{remark}

\section{Background on the partial limit distributions}
\label{sc4}

\begin{notations}\label{nt4.1}
In the convergence in distribution alluded to in the first paragraph of 
Notations \ref{nt3.2}, a key role will be played by the following 
class of probability functions for discrete random variables.
\smallskip

     (A) For each pair of real numbers $a,p \in (0,1)$ (the open unit interval), 
define the function ${\bf g}_{a,p}: \Z \to [0,1]$ as follows:
\begin{equation}
\label{eq4.1A1}
{\bf g}_{a,p}(0)\ :=\ 1- a, 
\indent {\rm and\ for\ each}\ n \in \N, \indent
{\bf g}_{a,p}(n)\ :=\ {\bf g}_{a,p}(-n)\ :=\ (a/2) p(1-p)^{n-1}\ . 
\end{equation}
For any given pair of numbers $a,p \in (0,1)$, the
numbers ${\bf g}_{a,p}(k)$, $k \in \Z$ are nonnegative and
(by a simple calculation) add up to 1.
That is, ${\bf g}_{a,p}$ is a probability function on the set $\Z$ of all integers.
For typographical convenience, the probability function ${\bf g}_{a,p}$ will 
sometimes be denoted as ${\bf g}[a,p]$; that is,
${\bf g}[a,p](k) = {\bf g}_{a,p}(k)$ for $k \in \Z$. 
\medskip

   (B)\ Next, as an elementary tool for use below, 
for each complex number $c$,
define the complex number $\upsilon_c$ by
\begin{equation}
\label{eq4.1B1}
\upsilon_c\ :=\ e^c - (1 + c)\ =\ \sum_{k=2}^\infty c^k/k!\ .
\end{equation}
Of course $\upsilon_0 = 0$.  
For each complex number $c$,
$|\upsilon_c| \leq \sum_{k=2}^\infty |c|^k/k!
\leq |c|^2 \sum_{k=2}^\infty |c|^{k-2}/(k-2)!
= |c|^2 \cdot e^{|c|}$.
Hence
\begin{equation}
\label{eq4.1B2}
\upsilon_c/c \to 0\ \ {\rm as}\ \ 
c \to 0\ ({\rm in}\ \bbC - \{0\}).
\end{equation}

(C) For each pair of real numbers $a,p \in (0,1)$, the characteristic 
function $\varphi_{a,b}$ for the probability function ${\bf g}_{a,p}$ can be 
computed in the following steps:
First, for each $t \in \bbR$, 
employing (\ref{eq4.1A1}) and (\ref{eq4.1B1}),
\begin{align}
\label{eq4.1C1}
\sum_{n=1}^\infty\ e^{itn} {\bf g}_{a,p}(n)\ 
&=\ \sum_{n=1}^\infty \Bigl[ e^{it} \cdot e^{it(n-1)} \cdot 
(a/2) p (1-p)^{n-1} \Bigl] \
=\ (a/2)pe^{it} \sum_{n=1}^\infty [e^{it}(1 - p)]^{n-1}
\nonumber\\
&=\ (a/2)pe^{it} \cdot \frac {1} {1 - e^{it}(1-p) }\
=\ (a/2)p \cdot \frac {1} {e^{-it} - (1-p) }
\nonumber\\
&=\ (a/2)p \cdot \frac {1} {1 - it + \upsilon_{-it} - (1-p) }\
=\ (a/2)p \cdot \frac {1} {p - it + \upsilon_{-it}}\ .    
\end{align} 
By an analogous calculation,  
\begin{equation}
\label{eq4.1C2}
\sum_{n=1}^\infty\ e^{it(-n)} {\bf g}_{a,p}(-n)\
=\ (a/2)p \cdot \frac {1} {p + it + \upsilon_{it}}\ .   
\end{equation}
Now by (\ref{eq4.1C1}) and (\ref{eq4.1C2}) and 
the first equality in (\ref{eq4.1A1})
(together with the trivial fact that $e^{it \cdot 0} = 1$), 
\begin{align}
\label{eq4.1C3}
\varphi_{a,p}(t)\ &=\ 
\sum_{k = -\infty}^\infty e^{itk} {\bf g}_{a,p}(k)\ 
=\ (1-a) \cdot 1\ +\ a \cdot \frac {p} {2} 
\biggl[ \frac {1} {p - it + \upsilon_{-it}}\ 
+\ \frac {1} {p + it + \upsilon_{it}} \biggl]
\nonumber \\
&=\ 1\ +\ a \cdot \biggl(\frac {p} {2}  
\biggl[ \frac {1} {p - it + \upsilon_{-it}}\
+\ \frac {1} {p + it + \upsilon_{it}} \biggl]\ -\ 1 \biggl)\ .\\
\nonumber    
\end{align}       
\end{notations}

\begin{lemma}
\label{lem4.2}
Suppose that for every ordered pair
$(a,p) \in (0,1)^2$, $J(a,p)$ is a positive integer; 
and suppose that 
\begin{equation}
\label{eq4.21}
\lim_{a \to 0+,\thinspace p \to 0+} a \cdot J(a,p)\ =\ 1\ .
\end{equation}
Suppose that for every ordered pair $(a,p) \in (0,1)^2$,
$\zeta^{(a,p)}_1, \zeta^{(a,p)}_2, \dots,
\zeta^{(a,p)}_{J(a,p)}$ are independent, identically distributed
integer-valued random variables, each with the probability function 
${\bf g}_{a,p}$ in eq.\ (\ref{eq4.1A1}) (with the same parameters $a$ and $p$).
Then 
\begin{equation}
\label{eq4.22}
p \cdot \sum_{k=1}^{J(a,p)} \zeta^{(a,p)}_k\
\to\ \mu_{P1sL}\ \ {\rm in\ distribution\ as}\ \  
a \to 0+,\  p \to 0 +.
\end{equation}
\end{lemma}

  {\bf Proof.}\ \ For convenience, for each ordered pair
$(a,p) \in (0,1)^2$, we shall define the random variable
$Y^{(a,p)}$ by  
\begin{equation}
\label{eq4.2P1}
Y^{(a,p)}\ :=\ p \cdot \sum_{k=1}^{J(a,p)} \zeta^{(a,p)}_k\
\end{equation}
[the left side of (\ref{eq4.22})], and we shall employ the  
(slight abuse of) notation $\varphi_{Y(a,p)}$ to denote its characteristic function.
Then for each $(a,p) \in (0,1)^2$,
the characteristic function of the random
variable $\zeta^{(a,p)}_1$ is the function
$\varphi_{a,p}$ in (\ref{eq4.1C3}), and hence the characteristic function of the 
random variable $p \cdot \zeta^{(a,p)}_1$ is
given by $t \mapsto E[\exp(it \cdot p \zeta^{(a,p)}_1)]
= \varphi_{a,p}(tp)$.
Applying (\ref{eq4.1C3}) (with $t$ there replaced by $tp$ here), one obtains that    
\begin{align}
\label{eq4.2P2}
\varphi_{Y(a,b)}(t)\ &=\ 
E\biggl[ \exp\Bigl(itY^{(a,p)}\Bigl) \biggl]\
=\ E\biggl[ \exp\biggl(it \sum_{k=1}^{J(a,p)} 
p\thinspace \zeta^{(a,p)}_k\biggl) \biggl]\
=\ \biggl(
E\Bigl[\exp\bigl(itp\thinspace \zeta^{(a,p)}_1\bigl)\Bigl]\biggl)^{J(a,p)}\nonumber \\
&=\ \biggl[
1\ +\ a \cdot \biggl(\frac {p} {2} 
\biggl[ \frac {1} {p - itp + \upsilon_{-itp}}\
+\ \frac {1} {p + itp + \upsilon_{itp}} \biggl]\ -\ 1 \biggl)   
\biggl]^{J(a,p)} \nonumber\\
&=\ \Biggl[
1\ +\ \frac {1} {J(a,p)} \cdot \Biggl(a \cdot J(a,p) \cdot 
\biggl(\frac {1} {2} 
\biggl[ \frac {1} {(1 - it) + \upsilon_{-itp}/p}\
+\ \frac {1} {(1 + it) + \upsilon_{itp}/p} \biggl]\ -\ 1 \biggl)   
\Biggl) \Biggl]^{J(a,p)}\ . \quad \\
\nonumber     
\end{align}

   Now for the calculations to follow here, suppose $t$ is any 
fixed real number.
   
   Let us first look at the first ``big fraction'' (inside the
``innermost big brackets'') in the last line of (\ref{eq4.2P2}).
Refer again to (\ref{eq4.1B1}) and (\ref{eq4.1B2}).
One has that
$\lim_{p \to 0+} \upsilon_{-itp}/p = 0$.
[For the case $t = 0$, recall that $\upsilon_0 = 0$.  
For the case $t \in \bbR - \{0\}$, simply note that as 
$p \to 0+$, one has that $-itp \to 0$ in $\bbC - \{0\}$ 
and hence 
$\upsilon_{-itp}/p = -it \cdot \upsilon_{-itp}/(-itp) \to 0$ by (\ref{eq4.1B2}).]    
Hence, inside the ``innermost big brackets'' in 
the last line of (\ref{eq4.2P2}), the ``first big fraction'' converges 
to $1/(1 - it)$ as $p \to 0+$.

Similarly, the ``second big fraction'' there converges to 
$1/(1 + it)$ as $p \to 0+$.
Hence the sum of those two ``big fractions''
converges to 
$[1/(1-it)] + [1/(1 + it)] = 2/(1 + t^2)$ as $p \to 0+$.

It now follows from (\ref{eq4.21}) that the entire expression inside 
the ``outer big parentheses'' in the last  
line\ of (\ref{eq4.2P2}) converges to $[1/(1 + t^2)] - 1$ as 
$a \to 0+,\thinspace p \to 0+$.

   Also, by (\ref{eq4.21}), $J(a,p) \to \infty$ as
$a \to 0+$ and $p \to 0+$. 
Now recall that if $c, c_1, c_2, c_3, \dots$ are complex
numbers such that $c_n \to c$ as $n \to \infty$, 
then $[1 + (1/n)c_n]^n \to e^c$ as $n \to \infty$.
It now follows that the last line of (\ref{eq4.2P2}) 
converges to $\exp([1/(1 + t^2)] - 1)$
as $a \to 0+,\thinspace p \to 0+$.   
Thus by (\ref{eq4.2P2}) (now its entire three-line display)
and (\ref{eq3.22}), one has that  
$\varphi_{Y(a,p)}(t) \to \varphi_{P1sL}(t)$
as $a \to 0+,\thinspace p \to 0+$.
\medskip

   That was shown to hold for arbitrary $t \in \R$.
It follows from the Continuity Theorem 
(see e.g.\ Theorem 26.3 on page 349 of
\cite{ref-journal-Bill}) that
$Y^{(a,p)} \to \mu_{P1sL}$ in distribution as
$a \to 0+, \thinspace p \to 0+$.   
Thus by (\ref{eq4.2P1}), eq.\ (\ref{eq4.22}) holds.
That completes the proof of Lemma \ref{lem4.2}.

\section{The ``building block'' Markov chains}
\label{sc5}

In Definition \ref{def5.2} below, we shall define a class of 
fairly simple strictly stationary, 3-state, reversible Markov chains that will 
be used in Section \ref{sc7} as
``building blocks'' in the construction of the Markov chain 
for Theorems \ref{thm3.3}.
\medskip

   The ``building block'' Markov chains here will be somewhat 
similar to the 2-state ``building block'' Markov chains in a 
construction in the paper \cite{ref-journal-Bradley2023}.
However, because of the different problems studied here
and in that paper, the parameterizations will be different.
In particular, the (``almost equal'') parameters $\theta$ and
$\theta^*$ below will roughly correspond to the
quantity $1 - \theta$ in \cite{ref-journal-Bradley2023}.
\medskip    

   In Notations \ref{nt5.1}, we shall formulate and 
develop some key properties of several items that will 
then be, explicitly or implicitly, key components 
in Definition \ref{def5.2}: \hfil\break
the (invariant) probability vectors, in eq.\ (\ref{eq5.131});
\hfil\break
the one-step transition probability matrices, in
eq.\ (\ref{eq5.141}); \hfil\break
the associated $n$-step transition probability matrices, in
eqs.\ (\ref{eq5.161}), with bounds in (\ref{eq5.171}); and
\hfil\break
the  {\it joint\/} probability functions of two consecutive
random variables, in (\ref{eq5.121}).     
\smallskip

   In Lemmas \ref{lem5.3} and \ref{lem5.4},
certain key properties (including reversibility) will be established for 
the Markov chains in Definition \ref{def5.2}.

\begin{notations}
\label{nt5.1}
 In the ``building block'' Markov chains in Definition \ref{def5.2}, the three states
will be $-1$, $0$, and $1$.
Accordingly, here and throughout Section \ref{sc5}, the $3 \times 3$ matrices will
be indexed by $\{-1,0,1\} \times \{-1,0,1\}$ (instead of, say, $\{1,2,3\} \times \{1,2,3\}$),
and the $1 \times 3$ row vectors will be indexed by $\{-1,0,1\}$.
\medskip
 
   {\bf Part 1.  The parameters.}\ \ In the notations that follow, the two main 
parameters will be ``small'' positive numbers 
$\zeps$ and $\theta$.  
As an ongoing reminder that they are intended to be
``small'', we shall impose the conditions
\begin{equation}
\label{eq5.111}
0 < \zeps \leq 1/9 \quad {\rm and} \quad
0 < \theta \leq 1/9\ .                    
\end{equation}
As a function of those two parameters 
satisfying (\ref{eq5.111}), we shall define another,
closely related parameter $\theta^*$ as follows:
\begin{equation}
\label{eq5.112}
\theta^*\ =\ \theta^*(\zeps, \theta)\ :=\
\frac {\theta} {1 - \zeps}\ .
\end{equation}
Of course under (\ref{eq5.111}) and (\ref{eq5.112}),
one has (with redundancy here), for convenient later 
reference, that 
\begin{equation}
\label{eq5.113}
(1 - \zeps)\thinspace \theta^* = \theta, \quad 
{\rm and} \quad
0 < \theta < \theta^* \leq 1/8 \quad {\rm and} \quad
1\ >\ 1 - \theta\ >\ 1 - \theta^*\ \geq\ 7/8\ >\ 0\ . 
\end{equation}

   {\bf Part 2.  Matrices for joint probabilities.}\ \ The following  $3 \times 3$ 
matrices will be used for {\it joint\/} (not conditional) probabilities
(for consecutive random variables in the Markov chains defined in
Definition \ref{def5.2}).
 For any $\zeps$ and $\theta$ satisfying (\ref{eq5.111}), 
 define the $3 \times 3$ matrix
$\Lambda^{(\zeps, \theta)} = 
(\lambda^{(\zeps, \theta)}_{i,j},\ i,j\in \{-1,0,1\})$ 
as follows:
\begin{align}
\label{eq5.121} 
\lambda^{(\zeps, \theta)}_{0,0}  
&=\ 1 - \zeps - \theta\zeps; \nonumber\\
\lambda^{(\zeps, \theta)}_{-1,-1}\ 
=\ \lambda^{(\zeps, \theta)}_{1,1}   
&=\ (1 - \theta)\zeps/2 ; \nonumber\\
\lambda^{(\zeps, \theta)}_{-1,1}\ 
=\ \lambda^{(\zeps, \theta)}_{1,-1}\ 
&=\ 0;\ \ {\rm and} \nonumber \\
\lambda^{(\zeps, \theta)}_{0,-1}\ 
=\ \lambda^{(\zeps, \theta)}_{0,1}\       
= \lambda^{(\zeps, \theta)}_{-1,0}\ 
=\ \lambda^{(\zeps, \theta)}_{1,0}\
&=\ \theta\zeps/2.  
\end{align}      
Note that the nine entries in this matrix are all 
nonnegative, and they add up to 1.
Also, note that this matrix is symmetric:
$\lambda^{(\zeps, \theta)}_{i,j} =
\lambda^{(\zeps, \theta)}_{j,i}$
for $i,j \in \{-1,0,1\}$. 
(There is further obvious symmetry as well.) 
 \medskip
 
     {\bf Part 3.  Vectors for marginal distributions.}\ \ 
For each $\zeps \in (0, 1/9]$ [as in eq.\ (\ref{eq5.111})],
define the $1 \times 3$ row vector 
$\pi^{(\zeps)} := 
[\pi^{(\zeps)}_{-1}, \pi^{(\zeps)}_0, \pi^{(\zeps)}_1]$
as follows:
\begin{equation}
\label{eq5.131}
\pi^{(\zeps)}_0\ =\ 1 - \zeps \indent {\rm and} \indent
\pi^{(\zeps)}_{-1}\ =\ \pi^{(\zeps)}_1\ =\ \zeps/2\ .
\end{equation}
By simple arithmetic, assuming (\ref{eq5.111}) as usual,  
if $(X,Y)$ is a random vector
taking its values in the set $\{-1,0,1\}^2$ with joint
probability function $\Lambda^{(\zeps, \theta)}$
in (\ref{eq5.121}), then for each of the random 
variables $X$ and $Y$ separately, the marginal probability function is 
$\pi^{(\zeps)}$.
\medskip  

     {\bf Part 4.  Matrices for one-step transition probabilities.}\ \ 
The matrices in eq.\ (\ref{eq5.141}) below will be employed later on 
as one-step transition probability matrices 
for the ``building block'' Markov chains in Definition \ref{def5.2}.
\medskip 
     
     For each $\zeps$ and $\theta$ satisfying
(\ref{eq5.111}), define the $3 \times 3$ matrix
$\bbP^{(\zeps, \theta)} := (p^{(\zeps, \theta)}_{ij},\ 
i,j \in \{-1,0,1\})$ as follows 
[employing (\ref{eq5.112}) 
in the first two rows of this display]:
\begin{align}
\label{eq5.141}
p^{(\zeps, \theta)}_{0,0}\ &=\ 1 - \theta^*\zeps;
\nonumber \\
p^{(\zeps, \theta)}_{0,-1}\ =\ p^{(\zeps, \theta)}_{0,1}\ 
&=\ \theta^* \zeps/2;
\nonumber \\  
p^{(\zeps, \theta)}_{-1,-1}\ =\ p^{(\zeps, \theta)}_{1,1}\  
&=\ 1 - \theta; \nonumber \\     
p^{(\zeps, \theta)}_{-1,0}\ =\ p^{(\zeps, \theta)}_{1,0}\ 
&=\ \theta; \nonumber \\ 
p^{(\zeps, \theta)}_{-1,1}\ =\ p^{(\zeps, \theta)}_{1,-1}\ 
&=\ 0. 
\end{align}
By (\ref{eq5.121}), (\ref{eq5.131}), (\ref{eq5.141}), and the initial
equality in (\ref{eq5.113}),  
one has [under (\ref{eq5.111})--(\ref{eq5.112})] that
\begin{equation}
\label{eq5.142}
{\rm for\ every}\ (i,j) \in \{-1,0,1\}^2,\ \ \
\pi^{(\zeps)}_i \cdot p^{(\zeps, \theta)}_{i,j}\ 
=\ \lambda^{(\zeps, \theta)}_{i,j}\ ; 
\end{equation}
and one obtains the matrix product equality
\begin{equation}
\label{eq5.143}
\pi^{(\zeps)} \bbP ^{(\zeps, \theta)} = \pi^{(\zeps)}\ .    
\end{equation}

{\bf Part 5.  Supporting matrices.}\ \ Next we shall work with three simple 
matrices that will be employed together to help simplify some 
calculations in Parts 6 and 7 below that involve powers of the
matrix in (\ref{eq5.141}).
\vskip 0.1 in

    Let $I_3$ denote the $3 \times 3$ identity matrix:
$I_3 := (\delta_{i,j}, i,j \in \{-1,0,1\})$ where
$\delta_{ij} = 1$ (resp.\ 0) if $i=j$ (resp.\ $i \neq j$).
\smallskip

    Next, for each $\zeps \in (0,1/9]$ [as in (\ref{eq5.111})], 
let $A^{(\zeps)} = (a^{(\zeps)}_{ij},\ i,j \in \{-1,0,1\})$ 
be the $3 \times 3$ matrix in which each row is the vector 
$\pi^{(\zeps)}$ in (\ref{eq5.131}).
That is,
\begin{align}
\label{eq5.151}
{\rm for\ each}\ i \in \{-1,0,1\}, \quad
a^{(\zeps)}_{i,0} = \pi^{(\zeps)}_0 = 1 - \zeps\ \ {\rm and}\ \ 
a^{(\zeps)}_{i,-1} = \pi^{(\zeps)}_{-1} = \zeps/2\ \ {\rm and}\ \  
a^{(\zeps)}_{i,1} = \pi^{(\zeps)}_1 = \zeps/2\ . 
\end{align}

   Next, define the $3 \times 3$ matrix 
$C := (c_{i,j},\ i,j\in\{-1,0,1\})$ as follows
\begin{align}
\label{eq5.152}   
&c_{-1,-1}\ =\ c_{1,1}\ =\ 1/2 \indent {\rm and} \indent
c_{-1,1}\ =\ c_{1,-1}\ =\ -1/2; \indent {\rm and} \nonumber\\ 
&c_{i,j}\ =\ 0\ \ \ {\rm if}\ i = 0\ {\rm or}\ j=0.    
\end{align}

   By the three preceding paragraphs together with (\ref{eq5.141}) 
and some simple arithmetic [including again the initial equality in (\ref{eq5.113})],
one has [under (\ref{eq5.111})--(\ref{eq5.112})] that
\begin{align}
\label{eq5.153}
&\bbP^{(\zeps, \theta)}\ 
=\ (1 - \theta^*) I_3\ +\ \theta^* [A^{(\zeps)} + \zeps C]\ ; 
\indent {\rm and} \indent \\
\label{eq5.154}
&{\rm for\ every\ ordered\ pair}\ (i,j) \in \{-1, 0, 1\}^2, \quad
\delta_{i,j} \geq 0, \ \ a^{(\zeps)}_{i,j} > 0,\ \ {\rm and}\ \  
\delta_{i,j} + a_{i,j}^{(\zeps)} + |c_{i,j}|\ \leq\ 2\ .
\end{align}

   Also, by (\ref{eq5.151}) and (\ref{eq5.152}) and  
simple matrix multiplication, $C^2 = C$; and   
for every $\zeps \in (0, 1/9]$ [as in (\ref{eq5.111})],  
$(A^{(\zeps)})^2 = A^{(\zeps)}$;
and also the matrices $A^{(\zeps)}$ and $C$ commute, with each of the 
matrix products $A^{(\zeps)}C$ and $CA^{(\zeps)}$ being equal to the 
zero matrix (the $3 \times 3$ matrix in which all nine entries are 0).
\medskip

   {\bf Part 6.  Matrices for $n$-step transition probabilities.}\ \
Refer again to the matrices $\bbP^{(\zeps, \theta)}$ in (\ref{eq5.141})
[under (\ref{eq5.111})--(\ref{eq5.112})].
For a given positive integer $n$, 
the matrix $(\bbP^{(\zeps, \theta)})^n$, the 
$n^{\rm th}$ power (under matrix multiplication) of $\bbP^{(\zeps, \theta)}$,
will be the $n$-step transition probability matrix for
the Markov chains in Definition \ref{def5.2}.
The entries of the matrix $(\bbP^{(\zeps, \theta)})^n$ will be examined in
Part 7 below. 
\smallskip

    But first our task here is to show that [under (\ref{eq5.111})--(\ref{eq5.112})]
one has that for every positive integer $n$,
\begin{equation}
\label{eq5.161}
(\bbP^{(\zeps, \theta)})^n\ 
=\ (1 - \theta^*)^n I_3\ 
+\ [1 - (1 - \theta^*)^n]A^{(\zeps)}\ 
+\ \bigl[ (1 - \theta)^n - (1 - \theta^*)^n \bigl] C.     
\end{equation}

     {\bf Proof of (\ref{eq5.161}).}\ \ 
Just for this argument, to reduce the clutter in the calculations, let
$y := 1 - \theta$ and $z := 1 - \theta^*$.
Then by (\ref{eq5.111})-(\ref{eq5.112}), $y > z$ and in fact
$y - z = \theta^* - \theta = \theta^* - (1 - \zeps)\theta^* = \zeps \theta^*$
by (\ref{eq5.113}); and hence by (\ref{eq5.153}), 
\begin{equation}
\label{eq5.162}
\bbP^{(\zeps, \theta)}\ 
=\ z I_3\ +\ (1-z) A^{(\zeps)}\  +\ (y-z) C.
\end{equation}
Our goal [rewriting (\ref{eq5.161})] is to show that for each positive integer $n$,
\begin{equation}
\label{eq5.163}
\bigl(\bbP^{(\zeps, \theta)}\bigl)^n\ 
=\ z^n I_3\ +\ (1 - z^n) A^{(\zeps)}\ +\ (y^n - z^n) C.     
\end{equation} 

   For $n=1$, (\ref{eq5.163}) holds by (\ref{eq5.162}).
Now for induction, suppose $n$ is any positive integer such that (\ref{eq5.163})
holds.
Our task is to show that (\ref{eq5.163}) holds with each $n$ replaced by $n+1$.
We shall use the representation
$(\bbP^{(\zeps, \theta)})^{n+1} = (\bbP^{(\zeps, \theta)})^n\, \bbP^{(\zeps, \theta)}$      
(matrix multiplication), and compute this matrix product using the right hand sides of
(\ref{eq5.163}) and (\ref{eq5.162}) (``in that order'').
Refer to the sentence after (\ref{eq5.154}).
The calculation here will produce nine terms, one of them being a scalar multiple 
of the identity matrix $I_3$ , three of them being scalar multiples of the matrix 
$A^{(\zeps)}$, three of them being scalar multiples of the matrix $C$,
and the other two (frivolously) being scalar multiples of the $3 \times 3$ zero matrix.      
Omitting that last frivolous pair of terms, one ends up with
\begin{align}
\nonumber      
(\bbP^{(\zeps, \theta)})^{n+1}\ &=\ (\bbP^{(\zeps, \theta)})^n \bbP^{(\zeps, \theta)}\
=\ z^{n+1}I_3\ +\ 
[z^n(1-z) + (1 - z^n)z + (1 - z^n)(1-z)] A^{(\zeps)} \\
\nonumber
& \hskip 1.8 in    
+\ [z^n(y-z) + (y^n-z^n)z +  (y^n-z^n)(y-z)] C \\
\nonumber
& =\ z^{n+1} I_3\ +\ [z^n(1-z) + (1 - z^n) \cdot 1\, ] A^{(\zeps)}\
+\ [z^n(y-z) + (y^n - z^n) \cdot y\, ] C \\
\nonumber
& =\ z^{n+1}I_3\ +\ (1 - z^{n+1}) A^{(\zeps)}\ +\ (y^{n+1} - z^{n+1}) C,       
\end{align}
which gives (\ref{eq5.163}) with $n$ replaced by $n+1$, as desired.   
Thus by induction, (\ref{eq5.163}) — that is, (\ref{eq5.161}) —
holds for every positive integer $n$.
That completes the proof of (\ref{eq5.161}). 
\medskip

     {\bf Part 7.  Bounds on the entries of the matrices $(\bbP^{(\zeps, \theta)})^n$.}\ \  
For each $\zeps$ and $\theta$ satisfying (\ref{eq5.111})
and each positive integer $n$, 
we use the representation 
$(\bbP^{(\zeps, \theta)})^n = \bigl(p^{(\zeps, \theta)(n)}_{i,j},\ i,j \in \{-1,0,1\}\bigl)$. 
That is, $p^{(\zeps, \theta)(n)}_{i,j}$ is 
simply the $(i,j)$-entry of the $n^{\rm th}$ power (under matrix multiplication) 
of the matrix $\bbP^{(\zeps, \theta)}$.
\smallskip

   Our task here in Part 7 is to show that [under (\ref{eq5.111})--(\ref{eq5.112})]
one has [see (\ref{eq5.131})] that for every positive integer $n$ and every
ordered pair $(i,j) \in \{-1,0,1\}^2$, 
\begin{equation}
\label{eq5.171}
\bigl|p^{(\zeps, \theta)(n)}_{i,j} - \pi^{(\zeps)}_j \bigl|\ 
\leq\ 2(1 - \theta)^n.
\end{equation} 

   {\bf Proof of (\ref{eq5.171}).}\ \
Suppose $\zeps$ and $\theta$ are as in (\ref{eq5.111}), $n$ is a positive integer,
and $(i,j) \in \{-1,0,1\}^2$. 
We shall again use the notations
$y := 1 - \theta$ and $z := 1 - \theta^*$ from the argument in Part 6 just above.
Again recall that $y > z$.
We start by repeating (\ref{eq5.161}) there, in the form given in (\ref{eq5.163}),
with one trivial operation of arithmetic:
$\bigl(\bbP^{(\zeps, \theta)}\bigl)^n - A^{(\zeps)} 
= z^n I_3 -  z^n A^{(\zeps)} + (y^n - z^n) C$.
Thus from the notations in Part 5 above,
\begin{equation}
\label{eq5.172}
p^{(\zeps, \theta)(n)}_{i,j}\ -\ a^{(\zeps)}_{i,j}\ 
=\ z^n\delta_{i,j}\ -\ z^n a^{(\zeps)}_{i,j}\ +\ (y^n-z^n) c_{i,j} . 
\end{equation}  

   Recall from (\ref{eq5.151}) that $a^{(\zeps)}_{i,j} = \pi^{(\zeps)}_j$.
Substituting that into the left side of (\ref{eq5.172}) and then taking absolute
values, and then using the fact that $0 < z< y$ and using
(\ref{eq5.154}), one obtains
\begin{align}
\nonumber
\bigl|p^{(\zeps, \theta)(n)}_{i,j}\ -\ \pi^{(\zeps)}_j \bigl|\ 
&=\ \bigl|z^n\delta_{i,j}\ -\ z^n a^{(\zeps)}_{i,j}\ +\ (y^n-z^n) c_{i,j} \bigl|\
\leq\ \bigl|z^n\delta_{i,j}|\ +\ \bigl|- z^n a^{(\zeps)}_{i,j}\bigl|\ 
+\ \bigl|(y^n-z^n) c_{i,j} \bigl| \\
\nonumber
&=\ z^n\delta_{i,j}\ +\ z^n a^{(\zeps)}_{i,j}\  +\ (y^n - z^n)|c_{i,j}|\
\leq\ y^n\delta_{i,j}\ +\ y^n a^{(\zeps)}_{i,j}\  +\ y^n |c_{i,j}|\
\leq\ 2y^n.      
\end{align}
Thus (\ref{eq5.171}) holds.  That completes the argument.
(End of Notations \ref{nt5.1})
\end{notations}

\begin{definition}
\label{def5.2}
     Suppose $\zeps$ and $\theta$ satisfy (\ref{eq5.111}).
Refer to (\ref{eq5.112}).
A Markov chain $X := (X_k, k \in \Z)$ 
is said to satisfy ``Condition ${\cal H}(\zeps, \theta)$'' if 
$X$ is a strictly stationary Markov chain such that \hfil\break
(i) its state space is $\{-1, 0, 1\}$, \hfil\break
(ii) its (invariant) marginal probability vector is the vector 
$\pi^{(\zeps)}$ in (\ref{eq5.131}), and \hfil\break
(iii) its one-step transition probability matrix is the matrix
$\bbP^{(\zeps, \theta)}$ in (\ref{eq5.141}).     
\smallskip

   That is [under Condition ${\cal H}(\zeps, \theta)$],
for $i,j \in \{-1, 0, 1\}$, 
$P(X_0 = i) = \pi_i^{(\zeps)}$ from (\ref{eq5.131}), 
and $P(X_1=j|X_0 = i) = p^{(\zeps, \theta)}_{i,j}$ from
(\ref{eq5.141}).  
The requirements (ii) and (iii) in this definition are
compatible by (\ref{eq5.143}).
In fact Condition ${\cal H}(\zeps, \theta)$ uniquely 
determines the distribution of the entire Markov chain $X$
on $(\R^\Z, \cR^\Z)$.  
\end{definition}

   As a quick summary, [under Condition 
${\cal H}(\zeps, \theta)$ for $\zeps$ and $\theta$ 
satisfying (\ref{eq5.111})], the Markov chain $X$ is 
strictly stationary, finite-state (in fact 3-state),
and [as an elementary by-product of (\ref{eq5.141})]
irreducible and aperiodic.
As a well known consequence (stated here in less than
full strength), the Markov chain $X$ 
satisfies absolute regularity ($\beta$-mixing),
and it is therefore also ergodic
(see e.g.\ [\cite {ref-journal-Bradley2007}, Vol.\ 1, 
Theorem 7.7 and Chart 5.22]).

\begin{lemma}
\label{lem5.3} 
Suppose $\zeps$ and $\theta$ satisfy (\ref{eq5.111}). 
Suppose $X := (X_k, k \in \Z)$ is a (strictly stationary) Markov chain 
that satisfies Condition ${\cal H}(\zeps, \theta)$ in Definition \ref{def5.2}.
Then the state space of this Markov chain $X$ is
$\{-1,0,1\}$; and the following statements hold: \hfil\break
\noindent (1) The (joint) probability function of the random 
vector $(X_0, X_1)$ is the matrix $\Lambda^{(\zeps, \theta)}$ 
from (\ref{eq5.121}).
\hfil\break
\noindent (2) This Markov chain $X$ is reversible. \hfil\break
\noindent (3a) One has that $P(X_0 = 0) = 1 - \zeps$ and
$P(X_0 = -1) = P(X_0 = 1) = \zeps/2$. \hfil\break
\noindent (3b) For each $s \in \{-1,0,1\}$, 
$P(\{X_0 = 0\} \cap \{X_1 = s\}) > 0$; in fact
$P(\{X_0 = 0\} \cap \{X_1 = 0\}) \geq 1 - 2\zeps$. \zhfb 
\noindent (4) One has that $EX_0 = 0$ and 
$E|X_0| =  E(X_0^2) = \var(X_0) = \zeps > 0$. \hfil\break
\noindent (5) The sequence $Y := (-X_k, k \in \bbZ)$
is a strictly stationary Markov chain that satisfies 
Condition ${\cal H}(\zeps, \theta)$. \hfil\break 
\noindent (6) For each $n \in N$, 
$\beta_X(n) \leq 6 \zeps (1 - \theta)^n 
\leq 6 \zeps \exp(-\theta n)$. \zhfb 
\noindent (7) For each $n \in \bbN$,
$P(X_0 = X_n = 1) = 
(\zeps/2)(1 - \theta)^n + P(\{X_0 = 1\} \cap \{X_n = -1\})$.
\hfil\break
\noindent (8a) For each $n \in \bbN$,
$E(X_0X_n) = \cov(X_0, X_n) = \zeps (1 - \theta)^n > 0$. \hfil\break
\noindent (8b) One has that 
$E(X_0^2) + 2 \sum_{n=1}^\infty E(X_0X_n) 
= \zeps [ (2/\theta) -1]$.
\hfil\break 
\noindent (9) Letting 
$u_n := (1/n) E[(X_1 + X_2 + \dots + X_n)^2]$ for
each $n \in \N$, one has that
\begin{align}
\label{eq5.301}
&u_1 \leq u_2 \leq u_3 \leq \dots \quad {\rm and} \quad
\lim_{n \to \infty} u_n = \zeps [(2 / \theta) - 1].  
\end{align}
\noindent (10) As $n \to \infty$, the normalized partial sum
$n^{-1/2} (X_1 + X_2 + \dots + X_n)$
converges in distribution to the normal law with 
mean 0 and variance $\zeps [(2 / \theta) - 1]$. 
\end{lemma}

   {\bf Proof.}\ \ The state space $\{-1,0,1\}$ is identified in
Definition \ref{def5.2}. 
   
   To verify statement (1), note that for each of the nine 
ordered pairs $(i,j) \in \{-1,0,1\}^2$, one has that
\begin{equation}
\label{eq5.31}
P(\{X_0 = i\} \cap \{X_1 = j\})\  
=\ \pi^{(\zeps)}_i \cdot p^{(\zeps, \theta)}_{ij}\
=\ \lambda^{(\zeps, \theta)}_{i,j}
\end{equation}
by Conditions (ii) and (iii) in Definition \ref{def5.2}
followed by (\ref{eq5.142}).

     To verify statement (2), recall that the matrix $\Lambda^{(\zeps, \theta)}$ 
 in (\ref{eq5.121}) is symmetric, and accordingly apply
(\ref{eq5.31}) and Notations \ref{nt2.8}(B).

   Statements (3a) and (4) are trivial consequences of 
(\ref{eq5.131}) and condition (ii) in Definition \ref{def5.2}.
The ``first half'' of statement (3b) follows from 
(\ref{eq5.31}) and (\ref{eq5.121}).
The ``second half'' of (3b) comes from just a trivial crude 
look at a complement:  
By (say) (3a), 
$P([\{X_0 = 0\} \cap \{X_1 = 0\}]^c)
= P(\{X_0 \neq 0\} \cup \{X_1 \neq 0\}) \leq 2\zeps$.     
\medskip

   {\bf Proof of (5).}\ \ 
The random sequence $Y$ in (5) 
is of the form $Y := (Y_k, k \in \Z)$ given by
$Y_k := f(X_k)$ (for each $k \in \Z$) 
where $f: \R \to \R$ is the function
defined by $f(x) = -x$ for $x \in \R$. 
Hence the sequence $Y$ inherits strict stationarity from the sequence $X$.
Since that function $f$ is one-to-one, the sequence $Y$ also
inherits the Markov property from the sequence $X$.
Thus $Y$ is a strictly stationary Markov chain. 
\smallskip

   By a triviality and then (say) (3a), one has that
$P(-X_k = i) = P(X_k = -i) = P(X_k = i)$ for each $i \in \{-1,0,1\}$.
By  a triviality and then (\ref{eq5.141}), one has that  
$P(-X_1 = j| -X_0 = i) = P(X_1 = -j|X_0 = -i) = P(X_1 = j|X_0 = i)$ 
for any ordered pair $(i,j) \in \{-1,0,1\}^2$.
Hence the (strictly stationary) Markov chain $Y$ in (5) satisfies
Condition $\cH(\zeps, \theta)$ by Definition \ref{def5.2}.
That completes the proof of statement (5).
\medskip

   {\bf Proof of (6).}\ \
   Let $n \in \N$ be arbitrary but fixed.   
    
    For each integer $k$, the sigma-field $\sigma(X_k)$ is 
of course purely atomic, with the three atoms
$\{X_k = -1\}$, $\{X_k = 0\}$, and $\{X_k = 1\}$
(with those three atoms forming a partition of the 
sample space $\Omega$).

By e.g.\ [\cite {ref-journal-Bradley2007},
Vol.\ 1, Proposition 3.21],  
\begin{equation}
\label{eq5.361}
\beta\bigl(\sigma(X_0), \sigma(X_n)\bigl)\ =\ 
\frac {1}{2}\, \sum_{i=-1}^1 \sum_{j=-1}^1
\Bigl|P(\{X_0 = i\} \cap \{X_n = j\}) -
P(X_0 = i) P(X_n = j)\Bigl|. 
\end{equation}

   For each fixed element $j \in \{-1,0,1\}$, the three events
$\{X_0 = i\} \cap \{X_n = j\}$ for $i \in \{-1,0,1\}$ together form a partition
of the event $\{X_n = j\}$ (with the empty set allowed in such a ``partition'').
Hence for each fixed $j \in \{-1,0,1\}$,    
\begin{align}
\label{eq5.362}
\sum_{i=-1}^1&\Bigl[P(\{X_0 = i\} \cap \{X_n = j\}) -
P(X_0 = i) P(X_n = j)\Bigl] \nonumber\\
&=\ \sum_{i=-1}^1 P(\{X_0 = i\} \cap \{X_n = j\})\ 
-\ \sum_{i=-1}^1 P(X_0 = i) P(X_n = j) \nonumber\\
&=\ P(X_n = j)\ -\ 
\biggl[P(X_n =j) \cdot \sum_{i=-1}^1 P(X_0 = i)\biggl]\
=\ P(X_n = j)\ -\ P(X_n = j) \cdot 1\ =\ 0.
\end{align}
Now for any three real numbers $q$, $r$, and $s$ such that
$q + r + s = 0$, one has that
$|r| = |-q-s| \leq |q| + |s|$ and hence
$|q| + |r| +|s| \leq 2(|q| + |s|)$.
Hence for each $j \in \{-1,0,1\}$, by (\ref{eq5.362}),
\begin{align}
\label{eq5.363}
\sum_{i=-1}^1& \Bigl|P(\{X_0 = i\} \cap \{X_n = j\}) -
P(X_0 = i) P(X_n = j)\Bigl| \nonumber \\ 
&\leq\ 2\, \cdot 
\sum_{i \in \{-1,1\}}\Bigl|P(\{X_0 = i\} \cap \{X_n = j\}) -
P(X_0 = i) P(X_n = j)\Bigl|\ .  
\end{align}
Applying (\ref{eq5.363}) for each $j$ to (\ref{eq5.361})
[say after switching the order of summation in 
(\ref{eq5.361})], one has that 
\begin{equation}
\label{eq5.364}
\beta\bigl(\sigma(X_0), \sigma(X_n)\bigl)\ \leq\ 
\sum_{i \in \{-1,1\}} \sum_{j=-1}^1
\Bigl|P(\{X_0 = i\} \cap \{X_n = j\}) -
P(X_0 = i) P(X_n = j)\Bigl|. 
\end{equation}

   Next, for each of the six ordered pairs 
$(i,j) \in \{-1,1\} \times \{-1,0,1\}$, one has from (\ref{eq5.131}) 
and (\ref{eq5.141}) that
\begin{equation}
\nonumber
P(\{X_0 = i\} \cap \{X_n = j\})\ -\ P(X_0 = i) P(X_n = j)\
=\ \pi^{(\zeps)}_i  p^{(\zeps,\theta)(n)}_{ij} 
- \pi^{(\zeps)}_i \pi^{(\zeps)}_j\ 
=\ (\zeps/2) \cdot
\Bigl(p^{(\zeps,\theta)(n)}_{ij} - \pi^{(\zeps)}_j \Bigl)
\nonumber   
\end{equation}
and hence by (\ref{eq5.171}),
\begin{equation}
\nonumber
\Bigl| P(\{X_0 = i\} \cap \{X_n = j\})\ 
-\ P(X_0 = i) P(X_n = j) \Bigl|\
\leq\ \zeps \cdot (1 - \theta)^n\ .
\end{equation}
Hence by (\ref{eq5.364}), 
$\beta (\sigma(X_0), \sigma(X_n)) 
\leq 6 \zeps (1 -\theta)^n$.
Thus [recall (\ref{eq2.42}) and its entire paragraph] the first inequality 
in statement (6) holds (for the given $n$).
\smallskip

   The second inequality in statement (6) comes from the elementary 
fact that $e^u \geq 1+u$ for every real number $u$.   
Since $n \in \N$ was arbitrary, statement (6) holds.
\bigskip

   {\bf Proof of (7).}\ \ 
 The proof will be divided into the two cases $n=1$ and $n \geq 2$.
 \medskip  
       
{\bf The case $n = 1$.}\ \  One has by
(\ref{eq5.31}) and (\ref{eq5.121}) that
$P(X_0 = X_1 = 1) = \lambda^{(\zeps, \theta)}_{1,1}  
= (\zeps/2) \cdot (1 - \theta)$
and similarly $P(\{X_0 = 1\}\cap\{X_1 = -1\}) = 0$.
Hence statement (7) holds for $n=1$. 
\vskip 0.1 in

   {\bf The case $n \geq 2$.}\ \ 
To complete the proof of statement (7),
now let $n \in \N$ be arbitrary but fixed such that $n \geq 2$.
 \smallskip
 
    Define the event $D := \bigcap_{k=0}^n \{X_k \neq 0\}$ and the event
$E_0 := \{X_0= 0\}$; and for each integer $j \in \{1,2,\dots, n\}$, define the
event $E_j := (\bigcap_{k=0}^{j-1} \{X_k \neq 0\}) \bigcap \{X_j = 0\}$.
(The last pair of equalities in (\ref{eq5.141}) will come into play shortly.)\ \ 
The events $D, E_0, E_1, E_2, \dots, E_n$ together form a partition of the
sample space $\Omega$.
Hence for any event $F$, the events $F \cap D$ and $F \cap E_j$, 
$j \in \{0,1,2, \dots, n\}$ together form a partition of of the event $F$ 
(with some events in that ``partition'' possibly being empty), and hence
$P(F) = P(F \cap D) + \sum_{j=0}^n P(F \cap E_j)$.
In particular, 
\begin{align}
\label{eq5.71}
&P(X_0 = X_n = 1)\ =\ 
P(\{X_0 = X_n = 1\} \cap D)\ +\ \sum_{j=0}^n P(\{X_0 = X_n = 1\} \cap E_j) 
\indent {\rm and}\\
\nonumber
&P(\{X_0 = 1\} \cap \{X_n = -1\})\ =\ 
P(\{X_0 = 1\} \cap \{X_n = -1\} \cap D)\\
\label{eq5.72} 
&\hskip 2 in \indent \quad +\ \sum_{j=0}^n P(\{X_0 = 1\} \cap \{X_n = -1\}\cap E_j).  
\end{align}

Let us first analyze the sums $\sum_{j=0}^n(\dots)$ in the right hand sides of
both (\ref{eq5.71}) and (\ref{eq5.72}).    
Of course the events $\{X_0 = X_n = 1\} \cap E_j$ and  
$\{X_0 = 1\} \cap \{X_n = -1\}\cap E_j$ are empty for $j = 0$ and for $j=n$
(recall that $E_n \subset \{X_n = 0\}$).
Hence in those two sums, the summands for $j=0$ and for $j=n$ are 0. 
Also, for each integer $j$ such that $1 \leq j \leq n-1$, by the Markov property
and the definition of the event $E_j$, one has that
\begin{align}
\label{eq5.73}
P(\{X_0 = X_n = 1\} \cap E_j)\ &=\ P(\{X_0 = 1\} \cap E_j) \cdot P(X_n = 1 | X_j = 0)
\quad {\rm and} \\
\label{eq5.74}
P(\{X_0 = 1\} \cap \{X_n = -1\} \cap E_j)\ &=\ 
P(\{X_0 = 1\} \cap E_j) \cdot P(X_n = -1 | X_j = 0).
\end{align}
For each integer $j$ such that $1 \leq j \leq n-1$, one has that
$P(X_n = -1|X_j = 0) = P(-X_n = 1|-X_j = 0) = P(X_n = 1|X_j = 0)$ by
a triviality followed by an application of statement (5), proved above.
(Recall the third sentence after Definition \ref{def5.2}.)\ \  
Hence for each integer $j$ such that $1 \leq j\leq n-1$, the right hand sides 
of (\ref{eq5.73}) and (\ref{eq5.74}) are equal,
and hence so are the left sides of those two equations equal.
Combining  that with the third sentence of this paragraph (involving $j=0$ and $j=n$),
one has that the sums $\sum_{j=0}^n(\dots)$ in the right hand sides of
(\ref{eq5.71}) and (\ref{eq5.72}) are equal.
Hence by subtracting (\ref{eq5.72}) from (\ref{eq5.71}), one obtains
\begin{align}
\label{eq5.75}
P(X_0 = X_n = 1)\ &-\ P(\{X_0 = 1\} \cap \{X_n = -1\}) \nonumber \\
&=\ P(\{X_0 = X_n = 1\} \cap D)\ -\ P(\{X_0 = 1\} \cap \{X_n = -1\} \cap D).   
\end{align}

     Next, by a simple argument involving the definition of the event $D$ and
the last line of (\ref{eq5.141}), one has that
$P(\{X_0 = 1\} \cap \{X_n = -1\} \cap D) = 0$.
Also, by (\ref{eq5.141}) (again including its last line) and again the definition
of the event $D$, and also (\ref{eq5.131}),
\begin{equation}
\nonumber
P(\{X_0 = X_n = 1\} \cap D)\ =\ 
P(X_0 = X_1 = X_2 = \dots = X_n = 1)\ =\ 
\pi^{(\zeps)}_1 (p^{(\zeps, \theta)}_{1,1})^n\ 
=\ (\zeps/2) (1 - \theta)^n. \quad 
\end{equation}
Hence by (\ref{eq5.75}), 
$P(X_0 = X_n = 1) - P(\{X_0 = 1\} \cap \{X_n = -1\}) = (\zeps/2) (1 - \theta)^n$. 
Thus statement (7) holds for the case $n \geq 2$.
That completes the proof of statement (7).
\medskip

   {\bf Proof of (8a).}\ \ For this proof, let $n$ be an arbitrary
fixed positive integer.
By statement (4) ($EX_0 = 0$), one has that
\begin{align}
\label{eq5.381}
\cov(X_0, X_n)\ &=\ E(X_0X_n)\ =\ 
\sum_{i \in \{-1,0,1\}}\ \sum_{j \in \{-1,0,1\}}
ij\cdot P(\{X_0 = i\} \cap \{X_n = j\}) \nonumber \\
&=\ 0\ +\ \bigl[1 \cdot P(X_0 = X_n = 1)\ 
            +\ (-1) \cdot P(\{X_0 = 1\} \cap (\{X_n = -1\}) \bigl] 
            \nonumber\\
& \indent \indent +\ \bigl[1 \cdot P(X_0 = X_n = -1)\
              +\ (-1) \cdot P(\{X_0 = -1\} \cap (\{X_n = 1\})\bigl]
\end{align}
Let us represent the last expression in ({\ref{eq5.381}) in the obvious way
as $0 + [A] + [B]$.  
By statement (5), the value of $[B]$ remains unchanged if $X_0$ and $X_n$ are 
both replaced respectively by $-X_0$ and $-X_n$.
Hence trivially $[A] = [B]$.
Also, statement (7) immediately gives that $[A] = (\zeps/2) (1 - \theta)^n$. 
Hence by (\ref{eq5.381}) itself, 
$\cov(X_0, X_n) = E(X_0X_n) = \zeps (1 - \theta)^n$. 
Thus (8a) holds.
\bigskip

   {\bf Proof of (8b).}\ \ Statement (8b) follows from 
statements (4) and (8a) and an elementary calculation.
\bigskip

   {\bf Proof of (9).}\ \ Refer to statements (4), (8a), and (8b).
For each nonnegative integer $m$, define the positive 
number $\gamma_m$ by
\begin{equation}
\label{eq5.391}
\gamma_m\ :=\ E(X_0X_m)\ =\ \zeps(1 - \theta)^m.
\end{equation} 
For each pair of positive integers $n$ and $\ell$, define 
the nonnegative number $\tau_{n,\ell}$ by
\begin{equation} 
\label{eq5.392}
\tau_{n,\ell}\ :=\ 
\begin{cases}
1 - \ell/n & {\rm if}\ n > \ell \\
0 & {\rm if}\ n \leq \ell.
\end{cases}                            
\end{equation}
For each positive integer $\ell$, one trivially has that
$0 \leq \tau_{1,\ell} \leq \tau_{2,\ell} 
\leq \tau_{3,\ell} \leq \dots < 1$, and that
$\tau_{n, \ell} \to 1$ as $n \to \infty$.

    For each integer $n \geq 2$ (say), referring to the 
notation $u_n$ in statement (9), one has that 
\begin{align}
\label{eq5.393}
u_n\ &=\ (1/n) \sum_{j=1}^n \sum_{k=1}^n E(X_jX_k)\
=\ (1/n) \sum_{j=1}^n \sum_{k=1}^n \gamma_{|j-k|} 
\nonumber\\
& =\ (1/n) \cdot n \gamma_0\ 
+\ (1/n) \cdot 2 \sum_{\ell = 1}^{n-1} (n-\ell) \gamma_\ell\
=\ \gamma_0\ 
+\ 2 \sum_{\ell = 1}^{n-1} \tau_{n,\ell} \gamma_\ell\
=\ \gamma_0\ 
+\ 2 \sum_{\ell = 1}^\infty \tau_{n,\ell} \gamma_\ell.      
\end{align}
Since the numbers $\gamma_\ell$ are positive (see
(\ref{eq5.391})), one has that 
$u_1 \leq u_2$ by a trivial calculation, and 
$u_n \leq u_{n+1}$ for each $n \geq 2$ by 
(\ref{eq5.393}) and the sentence right after (\ref{eq5.392}).
Thus the ``first half'' of (\ref{eq5.301}) holds.
The ``second half'' of (\ref{eq5.301}) holds 
by (\ref{eq5.393}) and the Monotone Convergence Theorem 
and (again) the sentence right after (\ref{eq5.392}),
followed by (\ref{eq5.391}) and statement 8(b).
That completes the proof of statement (9).
\medskip

   {\bf Proof of (10).}\ \ Referring to statements (8a) and (8b), define the 
positive number $\sigma$ via the equation
$\sigma^2 := E(X_0^2) + 2 \sum_{n=1}^\infty E(X_0X_n)= \zeps[(2/ \theta) - 1]$ 
(of course the sum converges absolutely).
Note that by statement (6) and trivial arithmetic (and eq.\ (\ref{eq2.34})),
$\sum_{n=1}^\infty \alpha_X(n) < \infty$. 
Now statement (10) holds by statement (4) ($EX_0 = 0$) 
and (as a quick convenient reference here) the central limit theorem 
of Ibragimov presented in Remark \ref{rem2.6}(A).
That completes the proof of Lemma \ref{lem5.3}.

\begin{lemma}
\label{lem5.4} 
Suppose $\zeps$ and $\theta$ satisfy (\ref{eq5.111}).
Suppose $X := (X_k, k \in \Z)$ is a strictly stationary
Markov chain that satisfies Condition $\cH(\zeps, \theta)$
(see Definition \ref{def5.2}).
Refer to (\ref{eq5.112}) and (\ref{eq5.113}).
Then there exists [on the same probability space 
$(\Omega, \cF, P)$] 
a sequence $(W_1, W_2, W_3, \dots)$
of independent, identically distributed, integer-valued
random variables with the following two properties: \zhfb
(1) The probability function of each $W_m$ ($m \in \bbN$)
is ${\bf g}[\theta^*\zeps, \theta]$ (recall Notations \ref{nt4.1}(A),
including its last sentence).
\zhfb      
(2) Letting $I$ be the positive integer defined by 
\begin{equation}
\label{eq5.41}
I\ :=\ \biggl[ \frac {1} {\theta^*\zeps} \biggl] 
\end{equation}
[where in this particular equation, the tall brackets indicate
that the right hand side is the greatest integer that is
$\leq 1/(\theta^*\zeps)$], one has that
\begin{equation}
\label{eq5.42}
P \biggl(\ \sum_{k=1}^I X_k\ \neq\ \sum_{j=1}^I W_j \biggl)\
\leq\ 3\zeps.
\end{equation}
\end{lemma}

   {\bf Proof.}\ \ The proof will be divided into 14 ``steps'',
with labels such as ``Step A'' and ``Claim E''.
\medskip

   {\bf Step A.}\ \ Here, to slightly simplify technicalities, 
we shall spell out a few trivial formalities.
\medskip
   
   {\bf Sub-step A1.}\ \ As noted in the paragraph right before
Lemma \ref{lem5.3}, the Markov chain $X$ is ergodic.
Refer to (\ref{eq5.131}) and the last pair of equalities in (\ref{eq5.141}).
Without loss of generality, we assume that for {\it every\/}
$\omega \in \Omega$, the following two conditions (a) and (b) hold: 

\noindent (a) For every state $s \in \{-1,0,1\}$, one has that
$X_k(\omega) = s$ for infinitely many positive integers $k$ and
infinitely many negative integers $k$. \hfil\break
(b) There does {\it not\/} exist $k \in \Z$ such that
the vector $(X_k(\omega), X_{k+1}(\omega))$ equals either
$(-1,1)$ or $(1,-1)$.
\smallskip

   {\bf Sub-step A2.}\ \ From Definition \ref{def5.2} and 
eqs.\ (\ref{eq5.131}) and (\ref{eq5.141}), 
for any pair of integers $j$ and $\ell$ such that $j \leq \ell$, and any vector 
$(s_j, s_{j+1}, \dots, s_\ell) \in \{-1,0,1\}^{\ell - j + 1}$
that does not have the states $-1$ and $+1$ next to each other
in either order, one has that
$P((X_j, X_{j+1}, \dots, X_\ell) = (s_j, s_{j+1}, \dots, s_\ell)) > 0$.
\smallskip

   That fact will trivially imply that certain conditional
probabilities later in this proof are legitimate --- that is, the conditioning is 
being done on events of positive probability.     
\medskip   

   {\bf Sub-step A3.}\ \ For the Markov chain $X$ here, the
Markov property can actually be formulated as follows:
For every integer $m$, every state $s \in \{-1,0,1\}$, and every pair of 
events $A \in \sigma(X_k, k \leq m)$ and
$B \in \sigma(X_k, k \geq m)$ such that
$P(A \cap \{X_m = s\}) > 0$, one has that
$P(B|A \cap \{X_m = s\}) = P(B|X_m = s)$.
\smallskip

   [The Markov property is traditionally formulated in that way 
for just
$A \in \sigma(X_k, k \leq m-1)$ and 
$B \in \sigma(X_k, k \geq m+1)$ 
(such that $P(A \cap \{X_m = s\}) > 0$).
The (equivalent) slightly ``augmented'' formulation here 
is just a standard elementary by-product.
[See e.g.\ statement (III)(B2) in Background A701 in the 
Appendix of Volume 1 of \cite {ref-journal-Bradley2007}.] 
\medskip

   {\bf Sub-step A4.}\ \ Since the state space of the Markov 
chain $X$ here is $\{-1,0,1\}$, one has the following:
\smallskip
   
   For every non-empty finite set $T \subset \bbZ$ and 
every $\omega \in \Omega$, 
the sum $\sum_{k \in T} |X_k(\omega)|$      
is simply the number of indices $k \in T$ such that
$|X_k(\omega)| = 1$ [that is,
$X_k(\omega) \in \{-1,1\}$].
\medskip

   {\bf Step B.}\ \ 
Define the nonnegative integer-valued random variables
$\kappa_0, \kappa_1, \kappa_2, \dots$ (recursively) as follows:   
For every $\omega \in \Omega$,
\begin{align}
\label{eq5.4B1}
&0\ \leq\ \kappa_0(\omega) < \kappa_1(\omega) < 
\kappa_2(\omega) < \dots;
\quad {\rm and} \\
\label{eq5.4B2}
& \Bigl\{ k \in \{0,1,2,\dots\}: X_k(\omega) = 0 \Bigl\}\
=\ \bigl\{\kappa_0(\omega),\ \kappa_1(\omega),\
\kappa_2(\omega),\ \dots\bigl\}\ . 
\end{align}     
[Eqs.\ (\ref{eq5.4B1}) and (\ref{eq5.4B2}) together 
uniquely define these random variables 
$\kappa_0, \kappa_1, \kappa_2, \dots$\thinspace .]
Of course one has the equality of events 
\begin{equation}
\label{eq5.4B3}
\{\kappa_0 = 0\} = \{X_0 = 0\}.
\end{equation}  

   For each $n \in \bbN$, define the random (possibly empty) 
set $Q(n)$ as follows:
For each $\omega \in \Omega$,
\begin{equation}
\label{eq5.4B4}
Q(n)(\omega)\ :=\ \bigl\{k \in \bbN: 
\kappa_{n-1}(\omega) < k < \kappa_n(\omega)\bigl\}.
\end{equation}
 
   By (\ref{eq5.4B1}), for any given $\omega \in \Omega$, these sets 
$Q(1)(\omega),\,  Q(2)(\omega),\,  Q(3)(\omega),\, \dots$
are (pairwise) disjoint.

   Also, for a given $n \in \bbN$ and a given 
$\omega \in \Omega$, the number of elements in the set
$Q(n)$ is exactly 
$\kappa_n(\omega) - \kappa_{n-1}(\omega) - 1$.   
If $\kappa_n(\omega) = 1 + \kappa_{n-1}(\omega)$,
then the set $Q(n)(\omega)$ is empty.
\smallskip

    Now define the random variables $W_1, W_2, W_3, \dots$
as follows:
For each $n \in \bbN$ and each $\omega \in \Omega$,
\begin{equation}
\label{eq5.4B5}
W_n(\omega)\ :=\ \sum_{k \in Q(n)(\omega)} X_k(\omega).
\end{equation}
In the case where the set $Q(n)(\omega)$ is empty,
i.e.\ $\kappa_n(\omega) = 1 + \kappa_{n-1}(\omega)$,  
one defines $W_n(\omega) := 0$, with the ``empty sum''
in the right hand side of (\ref {eq5.4B5}) understood 
to take the value 0.   
\medskip

   {\bf Step C.}\ \ Suppose $m$ is a positive integer and $\omega \in \Omega$.
(These items $m$ and $\omega$ will be fixed throughout the rest of Step C here.)\ \ 
Throughout Step C here, we shall freely tacitly use (\ref{eq5.4B1})-(\ref{eq5.4B2}) and
(\ref{eq5.4B4})-(\ref{eq5.4B5}) and the sentence after (\ref{eq5.4B5}). 
\medskip

   {\bf Sub-step C1.}\ \ Consider the following three (mutually exclusive) conditions
(A), (B), and (C):
\smallskip
   
    (A) One has that 
$\kappa_m(\omega) = \kappa_{m-1}(\omega) + 1$ [and 
hence $Q_m(\omega)$ is the empty set by (\ref{eq5.4B4})].
\smallskip

   (B) For some positive integer $n$, one has that (i)
$\kappa_m(\omega) = \kappa_{m-1}(\omega) + n + 1$ [and hence
\begin{equation}
\label{eq5.4C1}
Q_m(\omega)\ =\ \{ \kappa_{m-1}(\omega) + 1,\, \kappa_{m-1}(\omega) + 2,\, \dots,\,
\kappa_{m-1}(\omega) + n\}
\end{equation}
by (\ref{eq5.4B4})], and (ii) $X_k(\omega) = 1$ for every $k \in Q_m(\omega)$.
\smallskip

   (C) For some positive integer $n$, one has that (i) [as in (B)(i)] 
$\kappa_m(\omega) = \kappa_{m-1}(\omega) + n + 1$ [and hence
(\ref{eq5.4C1}) holds by (\ref{eq5.4B4})], and 
(ii) $X_k(\omega) = -1$ for every $k \in Q_m(\omega)$.
\medskip

   {\bf Sub-step C2.}\ \ By Sub-step A1, one has [for the given $m \in \N$ and 
$\omega \in \Omega$ here in Step C] that exactly one of the 
conditions (A), (B), (C) in Sub-step C1 holds; and if either (B) or (C) holds,
it does so for exactly one positive integer $n$.
\medskip

   {\bf Sub-step C3.}\ \ Referring again to the conditions in Sub-step C1, 
one has the following: \zhfb 
(a) If condition (A) holds, then $W_m(\omega) = 0$. \zhfb
(b) If condition (B) holds for a given $n \in N$, then $W_m(\omega) = n$ for that $n$. \zhfb 
(c) If condition (C) holds for a given $n \in N$, then $W_m(\omega) = -n$ for that $n$. 
\medskip

   {\bf Sub-step C4.}\ \ Refer again to the conditions in Sub-step C1.
By Sub-steps C2 and C3 and elementary logic, the following
converses of the three assertions in Sub-step C3 hold: \zhfb     
(a) If $W_m(\omega) = 0$, then condition (A) holds. \zhfb
(b) If $n \in \N$ and $W_m(\omega) = n$, then condition (B) holds for that $n$. \zhfb
(c) If $n \in \N$ and $W_m(\omega) = -n$, then condition (C) holds for that $n$. 
\medskip

   {\bf Step D.}\ \ We shall apply
Chung [\cite {ref-journal-Chung1967}, p.\ 84, Theorem 3].
With regard to the random variables $X_k$ here, the random 
variables $\kappa_j$ here, and Chung's use of the words``optional'' and 
``entrance time of a given state'', we trivially translate here from index 
set $\bbN$ (the index set employed by Chung) to index set $\{0,1,2,\dots\}$.
Under that translation, the random variable $\kappa_0$ 
(for example) is ``optional'', in that for each nonnegative 
integer $m$, the event $\{\kappa_0 = m\}$ is a member 
of the $\sigma$-field $\sigma(X_k, 0 \leq k \leq m)$.  
By that theorem in 
[\cite {ref-journal-Chung1967}, p.\ 84, Theorem 3]
--- on p.\ 83, line 2 there, take $f(x) = x$ for
each state $x \in \{-1,0,1\}$ [thus
$f(X_{\kappa(j)(\omega)}(\omega)) = 0$ for each
$j \in \{0,1,2,\dots\}$ and each $\omega \in \Omega$ 
by (\ref{eq5.4B2})] --- one has the following: \hfil\break
(i) the random variables $W_1, W_2, W_3, \dots$ are
independent and identically distributed; and \zhfb  
(ii) that entire random sequence $(W_1, W_2, W_3, \dots)$ 
is independent of the random variable $\kappa_0$.
\medskip

   {\bf Claim E.}\ \ {\em The common probability function of the
random variables $W_1, W_2, W_3, \dots$ is the function
${\bf g}[\theta^* \zeps, \theta]$ in (\ref{eq4.1A1})
[see the last sentence of Notations \ref{nt4.1}(A)].} 
\medskip

   {\bf Proof of Claim E.}\ \
Refer to (\ref{eq5.4B4}) and (\ref{eq5.4B5}).  
From statement (i) at the end of Step D, it suffices to verify 
Claim E for the random variable $W_1$.
By statement (ii) at the end of Step D, 
$P(W_1 = \ell) = P(W_1 = \ell | \kappa_0 = 0\}$
for every integer $\ell$.
Below, we shall simply cite ``Step D(ii)'' for that fact.
With that as our starting point, we have the following three observations: 
\smallskip

First, by Step D(ii), then Sub-steps C3(a) and C4(a) together, then
eqs.\ (\ref{eq5.4B1})-(\ref{eq5.4B2}), then eq.\ (\ref{eq5.4B3}), 
then (\ref{eq5.141}), and finally (\ref{eq4.1A1})
[and the last sentence of Notations \ref{nt4.1}(A)],
one has that      
\begin{align}
\label{eq5.4E1}
P(W_1 = 0)\  &=\ P(W_1 = 0 | \kappa_0 = 0)\
=\ P(\kappa_1 = 1 | \kappa_0 = 0)\ 
=\ P(X_1 = 0| \kappa_0 = 0) \nonumber\\
&=\ P(X_1 = 0| X_0 = 0)\ 
=\ p^{(\zeps, \theta)}_{0,0}\
=\ 1 - \theta^*\zeps\ 
=\ {\bf g}[\theta^*\zeps,\, \theta](0).  
\end{align}
Second, for every positive integer $n$ (including $n=1$) one has the
following:  By (again) Step D(ii), then Sub-steps C3(b) and C4(b) together, 
then eqs.\ (\ref{eq5.4B1})-(\ref{eq5.4B2}), then eq.\ (\ref{eq5.4B3}), 
then (\ref{eq5.141}), and finally (\ref{eq4.1A1})
[and the last sentence of Notations \ref{nt4.1}(A)],    
one has that 
\begin{align}
\label{eq5.4E2}
P(W_1 = n)\ &=\ P(W_1 = n | \kappa_0 = 0)\
=\ P \bigl( \{\kappa_1 = n+1\} \cap
\{ X_1 = X_2 = \dots = X_n = 1\} \bigl| \kappa_0 = 0 \bigl) \nonumber\\
&=\ P \bigl( \{X_{n+1} = 0\} \cap
\{ X_1 = X_2 = \dots = X_n = 1\} \bigl| \kappa_0 = 0 \bigl) \nonumber\\
&=\ P \bigl( \{X_{n+1} = 0\} \cap
\{ X_1 = X_2 = \dots = X_n = 1\} \bigl| X_0 = 0 \bigl) \nonumber\\ 
&=\ p^{(\zeps, \theta)}_{0,1} \cdot (p^{(\zeps, \theta)}_{1,1})^{n-1} 
\cdot p^{(\zeps, \theta)}_{1,0}\
=\ (\theta^* \zeps/2) \cdot (1 - \theta)^{n-1} \cdot \theta\
=\ {\bf g}[\theta^*\zeps,\, \theta](n). 
\end{align} 
Third, for every positive integer $n$ (including $n=1$),
by a calculation exactly analogous to (\ref {eq5.4E2}) but using 
Sub-steps C3(c) and C4(c) together [instead of Sub-steps C3(b) and C4(b)], 
one has that for every positive integer $n$ (again including $n=1$), 
$P(W_1 = -n) = (\theta^*\zeps/2) \cdot (1 - \theta)^{n-1} \cdot \theta 
= {\bf g}[\theta^*\zeps, \, \theta](-n)$.
Combining that with the entire sentences of
(\ref {eq5.4E1}) and (\ref {eq5.4E2}), one has that
for every integer $\ell$,
$P(W_1 = \ell) = {\bf g}[\theta^*\zeps, \, \theta](\ell)$.  
Thus Claim E holds.
\medskip

   {\bf Step F.}\ \ Refer again to the statement of Lemma \ref{lem5.4}.    
Refer to Claim E and line (i) at the end of Step D.
To complete the proof of Lemma \ref{lem5.4}, all that 
remains is to prove Statement (2) in that lemma.
\bigskip

   {\bf Proof of Statement (2).}\ \ 
The rest of the ``Steps'' here (starting with ``Step G'' below)
will be devoted to the proof of Statement (2).
The integer $I$ will of course be as defined in 
(\ref {eq5.41}).
For a given nonnegative integer $m$ and a given 
$\omega \in \Omega$, whenever the
nonnegative integer $\kappa_m(\omega)$
(see Step B) appears in a subscript or superscript, it will 
be written as $\kappa(m)(\omega)$ for typographical convenience.
\bigskip

   {\bf Step G.}\ \ By (\ref{eq5.4B1})-(\ref{eq5.4B2})   
and (\ref{eq5.4B4})-(\ref{eq5.4B5}), one has that    
for each positive integer $j$ and each $\omega \in \Omega$,
$X_{\kappa(j)(\omega)}(\omega) = 0$  
and hence
$W_j(\omega) = 
\sum_{k = \kappa(j-1)(\omega)+1}^{\kappa(j)(\omega)}X_k(\omega)$.
[That is, in that sum, the upper index can trivially be raised 
from its ``original natural value'' $\kappa_j(\omega) - 1$ to
$\kappa_j(\omega)$.]\ \  
Hence referring again to eq.\ (\ref{eq5.41}), 
one has that for every 
$\omega \in \Omega$, 
$\sum_{j=1}^I W_j(\omega) = 
\sum_{k = \kappa(0)(\omega) + 1} ^{\kappa(I)(\omega)}
X_k(\omega)$.   
Hence by (\ref{eq5.4B3}), the following holds:
\begin{equation}
\label{eq5.4G1}
{\rm For\ every}\ \omega \in \{X_0 = 0\},\ \ \ 
{\rm one\ has\ that}\ \ \ 
\sum_{j=1}^I W_j(\omega) 
= \sum_{k=1}^{\kappa(I)(\omega)}X_k(\omega).  
\end{equation}

   {\bf Step H.}\ \ Define the event $F$ (for ``fine'') as follows:
\begin{equation}
\label{eq5.4H1}
F\ :=\ \{X_0 = 0\} \cap \{X_I = 0\}.    
\end{equation}
For each integer $n \in \{0,1,2,\dots,I-1\}$, define the 
event $F_n$ as follows:
\begin{equation}
\label{eq5.4H2}
F_n\ :=\ F \bigcap\Biggl\{\sum_{k=1}^{I-1}|X_k| = n \Biggl\}\ . 
\end{equation}
For each integer $n \in \{1,2,\dots,I-1\}$ (note that $n=0$
is not included here), define the events $G_n$ (for ''good'')
and $H_n$ [for ``(slight) hindrance''] as follows:
\begin{align}
\label{eq5.4H3}
G_n\ &:=\ F_n\, \bigcap\, 
\Biggl\{\sum_{k=I+1}^{I+n}|X_k| = 0 \Biggl\} 
\quad {\rm and} \\
\label{eq5.4H4}
H_n\ &:=\ F_n\, \bigcap\,
\Biggl\{\sum_{k=I+1}^{I+n}|X_k| > 0 \Biggl\}\ .  
\end{align}
Since the state space of the Markov chain $X$ is $\{-1,0,1\}$, 
the events $F_0,\, F_1,\, \dots,\, F_{I-1}$ together     
form a partition of the event $F$.
For each $n \in \{1,2,\dots,I-1\}$, the events
$G_n$ and $H_n$ together form a partition of the 
event $F_n$.
Hence the events 
$F_0, G_1, H_1, G_2, H_2, \dots, G_{I-1}, H_{I-1}$
together form a partition of $F$.
Also, of course the events $F$ and $F^c$ together form a 
partition of the sample space $\Omega$.
Hence, the events
\begin{equation}
\label{eq5.4H5}
F^c,\ F_0,\ G_1,\ H_1,\ G_2,\ H_2,\ \dots,\ 
G_{I-1},\ H_{I-1}
\end{equation}  
together form a partition of the sample space $\Omega$.
\medskip

   {\bf Claim I.}\ \ {\em For any $\omega \in F_0$, one has that
$\kappa_I(\omega) = I$ and hence}
\begin{equation}
\label{eq5.4I1}
\sum_{k=1}^{\kappa(I)(\omega)} X_k(\omega)\ 
=\ \sum_{k=1}^I X_k(\omega)\ .
\end{equation}

   {\bf Proof of Claim I.}\ \ 
By (\ref{eq5.4H1}) and (\ref{eq5.4H2}), one has that 
$F_0 = \bigcap_{k=0}^I\{X_k = 0\}$.
Hence by (\ref{eq5.4B1})--(\ref{eq5.4B2})
and trivial induction, for every $\omega \in F_0$,
one has that $\kappa_j(\omega) = j$ for every
$j \in \{0,1,\dots, I\}$
and in particular, $\kappa_I(\omega) = I$. 
Eq.\ (\ref{eq5.4I1}) is (for $\omega \in F_0$) a trivial by-product.
Thus Claim I holds.
\medskip

   {\bf Claim J.}\ \ 
{\em Suppose $n \in \{1,2,\dots, I-1\}$ and $\omega \in G_n$.
Then $\kappa_{I-n}(\omega) = I$, and (\ref{eq5.4I1}) holds.}
\medskip

   {\bf Proof of Claim J.}\ \ 
Now $\omega \in G_n \subset F_n \subset F$ 
by (\ref{eq5.4H3}) and (\ref{eq5.4H2}). 
Now by (\ref{eq5.4H2}), there are exactly $n$ integers 
$k \in \{1,2,\dots, I-1\}$ such that $X_k(\omega) \neq 0$; 
and by (\ref{eq5.4H1}), all other integers $k \in \{0,1,2, \dots,I\}$
(which has exactly $I+1$ elements), including $k=0$ and $k=I$, 
satisfy $X_k(\omega) = 0$.
Hence the set $\{0,1,2,\dots, I\}$ has exactly $I+1-n$
integers $k$ (including $k=0$ and $k=I$) such that
$X_k(\omega) = 0$.
By (\ref{eq5.4B1})--(\ref{eq5.4B2}), those elements must be, 
in increasing order, the $I+1 - n$ integers
$\kappa_0(\omega), \kappa_1(\omega), \dots,
\kappa_{I-n}(\omega)$, 
with $\kappa_0(\omega) = 0$ and 
$\kappa_{I-n}(\omega) = I$.
\smallskip
   
     Since $\omega \in G_n$, one has by (\ref{eq5.4H3}) 
that $X_k(\omega) = 0$ for every 
$k \in \{I+1, I+2, \dots, I+n\}$.    
Hence by the last equality in the preceding paragraph, 
together with (\ref{eq5.4B1})--(\ref{eq5.4B2}) and trivial induction, 
one has that $\kappa_{I-n+j}(\omega) = I+j$ for every 
$j \in \{1,2, \dots, n\}$, and in particular (take $j=n$)
$\kappa_I(\omega) = I+n$.
And now one trivially has
$\sum_{k=I+1}^{\kappa(I)(\omega)} X_k(\omega)
= \sum_{k=I+1}^{I+n} X_k(\omega) = 0$.
Hence (\ref{eq5.4I1}) holds.  
That completes the proof of Claim J.
\medskip

   {\bf Step K.}\ \ Again recall from 
(\ref{eq5.4H1}), (\ref{eq5.4H2}), and
(\ref{eq5.4H3}) that
$F_0 \subset F \subset \{X_0 = 0\}$ 
and that for each $n \in \{1,2,\dots, I-1\}$,
$G_n \subset F_n \subset F \subset \{X_0 = 0\}$.
For each 
$\omega \in F_0 \bigcup(\bigcup_{n=1}^{I-1} G_n)$,
one now has by (\ref{eq5.4G1}) and Claims I and J that
$\sum_{j=1}^I W_j(\omega) 
=\sum_{k=1}^{\kappa(I)(\omega)} X_k(\omega)
= \sum_{k=1}^I X_k(\omega)$.
Hence by the entire sentence of (\ref{eq5.4H5}), 
one now has that
\begin{equation}
\label{eq5.4K1}
\Biggl\{\, \sum_{k=1}^I X_k \neq \sum_{j=1}^I W_j \Biggl\}\
\subset\
F^c \bigcup\biggl(\, \bigcup_{n=1}^{I-1} H_n \biggl)\ .
\end{equation}
 
   To complete the proof of eq.\ (\ref{eq5.42}), 
the remaining task is to get an appropriate bound on the 
probability of the right hand side of (\ref{eq5.4K1}).
\medskip

   {\bf Step L.}\ \ 
Of course $F^c \subset \{X_0 \neq 0\} \cup \{X_I \neq 0\}$
by (\ref{eq5.4H1}); and hence by Lemma 5.3(3a) (and stationarity),    
\begin{equation}
\label{eq5.4L1}
P(F^c)\ \leq\ 2 \cdot P(X_0 \neq 0)\ =\ 2\zeps.   
\end{equation}
We now turn our attention to the events $H_n$,
$1 \leq n \leq I-1$.
\medskip

  {\bf Claim M.}\ \ 
{\em For each $n \in \{1,2,\dots, I-1\}$, one has that}
\begin{equation}
\label{eq5.4M1}
P(H_n)\ \leq\ n \theta^* \zeps \cdot P(F_n).     
\end{equation}

    {\bf Proof of Claim M.}\ \ 
Let $n \in \{1,2,\dots, I-1\}$ be arbitrary but fixed.
Define the event $B$ (which will depend on the fixed $n$) by
\begin{equation}
\label{eq5.4M2}
B\ :=\ \Biggl\{\, \sum_{k=I+1}^{I+n} |X_k| > 0 \Biggl\}\ .
\end{equation}
Then by (\ref{eq5.4H4}), 
\begin{equation}
\label{eq5.4M3}
   H_n\ =\ F_n \cap B\ .    
\end{equation}

   Now by (\ref{eq5.4H1}) and (\ref{eq5.4H2}), 
$F_n \in \sigma(X_k, k \leq I)$ and also
$F_n \subset \{X_I = 0\}$
(and hence $F_n = F_n \cap \{X_I = 0\}$);
and by (\ref{eq5.4H1}) and (\ref{eq5.4H2}) 
and Sub-step A2 in Step A, $P(F_n) > 0$. 
Also, by (\ref{eq5.4M2}), 
$B \in \sigma(X_k, k \geq I+1)$. 
Hence by the Markov property,
\begin{equation}
\label{eq5.4M4}
P(B|F_n)\ =\ P(B| F_n \cap \{X_I = 0\}) = P(B|X_I = 0).
\end{equation}

    Now define the event
\begin{equation}
\label{eq5.4M5}
B_1\ :=\ \{X_{I+1} \neq 0\}\ ;    
\end{equation}
and for each $j \in \{2,3, \dots, n\}$ (if $n \geq 2$), 
define the event
\begin{equation}
\label{eq5.4M6}
B_j\ :=\ 
\biggl[\, \bigcap_{k = I+1}^{I + j -1} \{X_k = 0\} \biggl]
\bigcap \{X_{I + j} \neq 0\}\ . 
\end{equation}
Then by (\ref{eq5.4M2}), (\ref{eq5.4M5}), 
and (\ref{eq5.4M6}), 
\begin{equation}
\label{eq5.4M7}
B\ =\ \bigcup_{j=1}^n B_j\ .    
\end{equation}

   By (\ref{eq5.4M5}), (\ref{eq5.4M6}), and (\ref{eq5.141}), 
\begin{equation}
\label{eq5.4M8}
P(B_1 | X_I = 0)\ =\ 
p^{(\zeps, \theta)}_{0,-1}\ +\ p^{(\zeps, \theta)}_{0,1}\
=\ \theta^*\zeps\ ; 
\end{equation}
and for each $j \in \{2,3, \dots, n\}$ (if $n \geq 2$), 
\begin{equation}
\label{eq5.4M9}
P(B_j| X_I = 0)\ =\ (p^{(\zeps, \theta)}_{0,0})^{j-1}
(p^{(\zeps, \theta)}_{0,-1}\ +\ p^{(\zeps, \theta)}_{0,1})\
\leq\ 
p^{(\zeps, \theta)}_{0,-1}\ +\ p^{(\zeps, \theta)}_{0,1}\
=\ \theta^*\zeps\ .  
\end{equation}

   It is unnecessary but trivial to note here that (if $n \geq 2$)
the events $B_1, B_2, \dots, B_n$ are (pairwise) disjoint.   
By (\ref{eq5.4M4}), (\ref{eq5.4M7}), 
(\ref{eq5.4M8}), and (\ref{eq5.4M9}),
\begin{equation}
\label{eq5.4M10}
P(B|F_n)\ =\ P(B|X_I = 0)\ =\ \sum_{j=1}^n P(B_j|X_I = 0)\
\leq\ n \theta^* \zeps\ .
\end{equation}
Hence by (\ref{eq5.4M3}) and (\ref{eq5.4M10}),
\begin{equation}
\nonumber
P(H_n)\ =\ P(F_n \cap B)\ =\ P(F_n) \cdot P(B|F_n)\
\leq\ P(F_n) \cdot n\theta^*\zeps\ .
\end{equation}
That is, (\ref{eq5.4M1}) holds.
That completes the proof of Claim M.
\medskip

  {\bf Step N.}\ \ 
Now by Claim M, eq.\ (\ref{eq5.4H2}), Lemma \ref{lem5.3}(3a)(4)
[$E|X_0| = \zeps$], and eq.\ (\ref{eq5.41})
(and including one trivial extra term in the third sum here below),
one has that
\begin{align} 
\label{eq5.4N1}
\sum_{n=1}^{I-1} P(H_n)\ 
&\leq\ \sum_{n=1}^{I-1} \theta^*\zeps n P(F_n)\
\leq\ 
\theta^*\zeps \sum_{n=0}^{I-1}
\Biggl[\, n \cdot P\Biggl(\, \sum_{k=1}^{I-1} |X_k|\, =\, n
\Biggl) \Biggl]\ 
=\ \theta^* \zeps E \Biggl(\, \sum_{k=1}^{I-1} |X_k| \Biggl) 
\nonumber\\
&=\ \theta^*\zeps (I-1) \zeps\
<\ \theta^*\zeps  \frac {1} {\theta^*\zeps} \zeps\
=\ \zeps.    
\end{align}

   Now by (\ref{eq5.4N1}) and (\ref{eq5.4L1}),
$P(F^c \bigcup(\bigcup_{n=1}^{I-1}H_n)) \leq 3\zeps$.
Hence by (\ref{eq5.4K1}), eq.\ (\ref{eq5.42}) holds.
Thus Statement (2) in Lemma \ref{lem5.4} holds.
That completes the proof of Lemma \ref{lem5.4}. 

\section{Other Preliminaries}
\label{sc6}

   In this section, we shall collect together several technical 
background facts in analysis and probability theory that will be employed in 
Section \ref{sc7} in the proof of Theorem \ref{thm3.3}.
For a couple of elementary statements here of a somewhat ad hoc nature,
proofs are included for convenience. 
\smallskip    

  We start with a standard technical lemma. 
A proof for it can be found e.g.\ in
[\cite {ref-journal-Bradley2007}, Vol.\ 1, Theorem 6.2(III)]. 

\begin{lemma}
\label{lem6.1}
Suppose $\cA_1, \cA_2, \cA_3, \dots$ and
$\cB_1, \cB_2, \cB_3, \dots$ are 
$\sigma$-fields [on the given probability space
$(\Omega, \cF, P)$], and the
$\sigma$-fields $\cA_j \vee \cB_j$, $j \in \N$
are independent.
Then
\begin{equation}
\nonumber
\beta\biggl(\, \bigvee_{j=1}^\infty \cA_j,\ 
\bigvee_{j=1}^\infty \cB_j \biggl)\
\leq\ \sum_{j=1}^\infty \beta(\cA_j, \cB_j). 
\end{equation}  
\end{lemma}

   The following lemma is of course well known. 
It is stated here for convenient reference.

\begin{lemma}
\label{lem6.2}
     Suppose $V_1, V_2, V_3, \dots$ is a sequence
of independent square-integrable random variables such that
(i) $EV_j = 0$ for each $j \in \N$, and
(ii) $\sum_{j=1}^\infty E(V_j^2) < \infty$.
     Then the sum $\sum_{j=1}^\infty V_j$ converges almost surely and 
in $\cL^2$ to a square-integrable random variable $Y$ such that $EY = 0$ 
and $E(Y^2) = \sum_{j=1}^\infty E(V_j^2)$, 
\end{lemma}

   The next lemma gives just an elementary technical fact that 
will be employed in order to show that the construction in 
Sections \ref{sc7} is in fact a Markov chain ---
not just a function of a Markov chains. 

\begin{lemma}
\label{lem6.3}
    Suppose $(y_1, y_2, y_3, \dots)$ and $(z_1, z_2, z_3, \dots)$
are each a sequence of elements of the set $\{-1,0,1\}$.
Suppose that (i) there exists $m \in \N$ such that
$y_k = z_k = 0$ for all $k > m$, and
(ii) there exists $k \in \N$ such that $y_k \neq z_k$.

   Suppose $r_1, r_2, r_3, \dots$ is a sequence of positive numbers 
such that $r_{k+1}/r_k \geq 3$ for every $k \in \N$.

   Then each of the sums $\sum_{k=1}^\infty r_ky_k$ and 
$\sum_{k=1}^\infty r_kz_k$ has at most finitely many non-zero summands 
(and therefore trivially converges absolutely); and these sums satisfy 
$\sum_{k=1}^\infty r_ky_k \neq \sum_{k=1}^\infty r_kz_k$.
\end{lemma}
   
    {\bf Proof.}\ \ 
Let $m$ be as in assumption (i) in the statement of the lemma.
Let $j$ denote the greatest element of $\{1, 2, \dots, m\}$
such that $y_j \neq z_j$.
If $j=1$, then the final inequality in Lemma \ref{lem6.3} holds trivially.
Therefore, suppose instead that $j \in \{2,3,\dots,m\}$.
Obviously it suffices to show that 
$\sum_{k=1}^j r_ky_k \neq \sum_{k=1}^j r_kz_k$.
\smallskip

   For each $k \in \{1,2,\dots,j-1\}$,
$r_k\leq 3^{k-j} r_j$ and $|y_k - z_k| \leq 2$ and hence
$|r_ky_k - r_kz_k| = r_k |y_k - z_k| \leq 3^{k-j}r_j \cdot 2$. 
Hence 
\begin{equation}
\label{eq6.301}
\biggl|\sum_{k=1}^{j-1} r_ky_k\ 
-\ \sum_{k=1}^{j-1} r_kz_k\biggl|\ 
\leq\ 2 r_j \cdot 3^{-j}\sum_{k=1}^{j-1} 3^k\
=\ (2 r_j \cdot 3^{-j}) \cdot (1/2) \cdot (3^j-3)\
<\ r_j. 
\end{equation}
Also, 
$|r_jy_j - r_j z_j| = r_j |y_j - z_j| \geq r_j \cdot 1$.
From that and (\ref{eq6.301}) and simple arithmetic, the
desired inequality 
$\sum_{k=1}^j r_ky_k \neq \sum_{k=1}^j r_kz_k$
follows.
That completes the proof.

\begin{notations}
\label{nt6.4}
     (A) In Lemma \ref{lem6.5} below, and in arguments in 
Section \ref{sc7}, the term ``affine function $L: \R \to \R$'' 
of course just means a function (``line'') given by 
$L(x) = q + rx$ for $x \in \R$, 
where $q$ and $r$ are each a real number.
For such an affine function $L$, the notation
``$[{\rm slope\ of}\ L]$'' will simply mean its (constant) 
derivative (the parameter $r$ in the preceding sentence).
\smallskip

     (B) In fact Lemma \ref{lem6.5} will involve lines that are tangent to
convex functions at certain points. 
In that context, the following elementary facts from analysis will
be tacitly kept in mind:

     If $I$ is an open interval or open half line $\subset \R$, or $I$ is the whole
real number line $\R$ itself, and $\phi: I \to \R$ is a convex function, then \zhfb
(i) $\phi$ is continuous on $I$, and \zhfb
(ii) if $\Gamma$ is the set of all points $x \in I$ such that $\phi$ is
{\it not\/} differentiable at $x$, then that set $\Gamma$ is at most countable,
and in particular the set $I - \Gamma$ 
(of points in $I$ where $\phi$ {\it is} differentiable) is dense in $I$.
\smallskip

   (In arguments below and in Section \ref{sc7} involving convex functions, 
it will not be really necessary to exclude that set $\Gamma$, but doing so 
will perhaps make the arguments a little tidier.) 
\end{notations}

\begin{lemma}
\label{lem6.5}
Suppose $\phi: [1,\infty) \to (-\infty, 0]$ is a function
that satisfies the following conditions: \zhfb
(a) $\phi$ is strictly decreasing and convex on $[1, \infty)$; and \zhfb
(b) one has that
\begin{equation}
\label{eq6.501}
\lim_{x \to \infty} \phi(x)\ =\ -\infty 
\indent {\rm and} \indent
\lim_{x \to \infty} \frac {\phi(x)} {x}\ =\ 0.
\end{equation}

   Let $\Gamma$ denote the set of all points $x \in (1, \infty)$ (the open half line) 
such that $\phi$ is {\rm not} differentiable at $x$.
Refer to Notations \ref{nt6.4}(A)(B).   
For each $y \in (1, \infty) - \Gamma$, let 
$L^{(y)}: \R \to \R$ denote the affine function whose graph is tangent to that 
of $\phi$ at the point $(y, \phi(y))$ --- that is, the affine function that satisfies
\begin{equation}
\label{eq6.502}
L^{(y)}(y)\ =\ \phi(y) \indent {\rm and} \indent 
L^{(y)}(x)\ \leq\ \phi(x)\ \ {\rm for\ all}\ x \in [1, \infty). 
\end{equation}
 
   Then the following statements (I) and (II) hold:
\medskip

   (I) The following statements hold: \zhfb 
(i) For each $x \in (1, \infty)$, $\phi(x) < \phi(1) \leq 0$. \zhfb 
(ii) For each $y \in (1,\infty) - \Gamma$,  
$[{\rm slope\ of}\ L^{(y)}] = \phi'(y) < 0$. \zhfb   
(iii) The function $x \mapsto \phi'(x)$ for $x \in (1,\infty) - \Gamma$, 
is nondecreasing. \zhfb  
(iv) One has that $\lim_{x \to \infty,\, x \notin \Gamma} \phi'(x) = 0$. 
\medskip

   (II) For any negative number $D$ and any positive number 
$s$, there exists a number    
$T = T(\phi, D, s) > 1$ with the following property:

   For each $y \in [T,\infty)  - \Gamma$, the affine function 
$L^{(y)}$ in (the entire sentence of) eq.\ (\ref{eq6.502})
satisfies 
\begin{equation}
\label{eq6.521}
-s\ \leq\ [{\rm slope\ of}\ L^{(y)}]\ =\ \phi'(y)\ <\ 0 
\indent {\rm and} \indent
L^{(y)}(0)\ \leq\ D. 
\end{equation}
\end{lemma}

   {\bf Proof.}\ \ {\bf Proof of (I).}\ \ 
In Statement (I), item (i) and the inequality $\phi'(y)< 0$ in (ii)
follow from assumption (a) and the notations in the first line
of the lemma.  
The equality in (ii) follows from (\ref{eq6.502}) and convexity.  
Item (iii) follows from convexity.
Item (iv) follows from (iii) and the ``second half'' of (\ref{eq6.501}).
Thus Statement (I) holds.
\medskip

   {\bf Proof of (II).}\ \
This statement and its geometry are elementary.
A detailed proof seems a little awkward to write out.
To facilitate that process, we shall label the paragraphs 
in the proof as (P1), (P2), etc.
\smallskip
    
     (P1) Suppose $D< 0$ and $s>0$, as assumed in Statement (II).
Referring to (\ref{eq6.501}) (its ``first half'') and to Statement (I)(iv),       
let $w \in (1, \infty) - \Gamma$ be chosen sufficiently large that
$\phi(w) < D$ and $-s < \phi'(w)$. 

     (P2) Referring to (P1), let $q$ be a number such that $\phi(w) < q < D$.

     (P3) Referring to (P1) and (P2), define the number $\eta$ by $\eta := (D-q)/w$.
Then $\eta > 0$.
Also, $q = D - \eta w$.

     (P4) Define the affine function $L^*: \R \to \R$ as follows: 
$L^*(x) := D - \eta x$ for all $x \in \R$.
Then $L^*(w) = D - \eta w = q > \phi(w)$ by (P3) and (P2).

     (P5) Recall that $\lim_{x \to \infty} \phi(x)/x = 0$ by (\ref{eq6.501});
and note that $\lim_{x \to \infty} L^*(x)/x = - \eta < 0$ by (P4) and (P3). 
It follows that for all $x \in (1,\infty)$ sufficiently large, one has that
$\phi(x)/x > L^*(x)/x$. 
Hence $\phi(x) > L^*(x)$ for all $x \in (1, \infty)$ sufficiently large.
Also, $\phi(w) < L^*(w)$ by (P4).
Also, the function $\phi$ is continuous on $(1, \infty)$ (since it is 
convex on that open half line).
Hence there exists $z > w$ ($z$ henceforth fixed) such that
$\phi(z) = L^*(z)$.

     (P6) By (P1) and (P5), $1 < w < z$. 
Let $T := T(\phi, D, s) \in (1, \infty)$ be a number such that $T > z$.

     (P7) Now suppose $y$ is any number such that $y \in [T ,\infty) - \Gamma$. 
To complete the proof of (II), our task is to verify (\ref{eq6.521}) for this $y$.

     (P8) By (P6) and (P7), $1 < w < z < T \leq y$.
Hence by (P1) and Statement (I)(ii)(iii),
$-s < \phi'(w) \leq \phi'(y) = [{\rm slope\ of}\ L^{(y)}] < 0$.
Thus the ``first half''  of (\ref{eq6.521}) holds.  

     (P9) Recall that $\phi(w) < L^*(w)$ by (P4), and $\phi(z) = L^*(z)$ by (P5),
and that $w < z$ (see e.g.\ (P8)).
It follows from the convexity of $\phi$ (and again the first sentence of (P8))
that $\phi(y) > L^*(y)$.
Substituting from the ``first half'' of 
(\ref{eq6.502}), one thus has that $L^{(y)}(y) > L^*(y)$.

   (P10) Refer again to all of  the first sentence of (P8).   
By the ``second half'' of (\ref{eq6.502}), and then (P4), 
$L^{(y)}(w) \leq \phi(w) < L^*(w)$.
From that and the last sentence of (P9), it follows that
$L^{(y)}(x) < L^*(x)$ for all $x \leq w$
(since the function $x \mapsto L^{(y)}(x) - L^*(x)$ for $x \in \R$ is affine).      
In particular, $L^{(y)}(0) < L^*(0) = D - \eta \cdot 0 = D$ by (P4).
Thus the ``second half'' of (\ref{eq6.521}) holds.

   (P11) By the last sentences of (P7), (P8), and (P10), the
proof of Statement (II) is complete.
That completes the proof of Lemma \ref{lem6.5}.

\begin{lemma}
\label{lem6.6}
Suppose $(E_1, E_2, E_3, \dots)$ is a sequence of independent events
on a probability space $(\Omega, \cF, P)$, such that 
$\sum_{j=1}^\infty P(E_j^c) < \infty$ (where the superscript 
$c$ denotes complement).

   Suppose that on the same probability space, $(F_1, F_2, F_3, \dots)$ is a
sequence of independent events such that
(i) $P(F_j) > 0$ for each $j\in \N$, and 
(ii) $F_j = E_j$ (an equality of events) for all except at most finitely
many indices $j \in \N$.
Then $P(\bigcap_{j=1}^\infty F_j) = \prod_{j=1}^\infty P(F_j) > 0$.  
\end{lemma}

   This formulation is rather contrived, but will fit our applications later on.  
Of course the final equality holds by independence.
The subsequent final inequality ($\dots > 0$) is simply an application of the 
basic fact in analysis that if $(s_1, s_2, s_3, \dots)$ is a sequence of 
numbers in $[0,1)$ such that $\sum_{j=1}^\infty s_i < \infty$,
then $\prod_{j=1}^\infty (1 - s_j) > 0$.
[Simply take $s_j = P(F_j^c)$.]

\section{Proof of Theorem \ref{thm3.3}}
\label{sc7}
In this proof, in order to simplify certain mathematical arguments, 
we shall sometimes build unnecessary redundancies into the
choices of parameters, and also we shall sometimes settle on adjectives of
``less than full strength'' (such as stating that a sequence of real numbers
is ``monotonically decreasing'' when in fact it is strictly decreasing) when
``full strength'' is unnecessary. 
\medskip      

   The proof will be divided into seven ``steps''.
\medskip

     {\bf Step 1.}\ \ This step will be divided into several ``sub-steps''.
\medskip

     {\bf Sub-step 1A.}\ \ Refer to assumptions (a) and (b) in the statement of 
Theorem \ref{thm3.3}.
To make the presentation of the argument superficially a little tidier, 
reducing $\zeta_n$ for just finitely many positive integers $n$ if necessary,
we also assume without loss of generality that 
\begin{equation}
\label{eq7.1A1}
0 < \zeta_n < 1\ \ \ {\rm for\ every}\ {n \in \N}.    
\end{equation}
     
     For each $n \in \N$, define the real number $\xi_n$ by
\begin{equation}
\label{eq7.1A2}
\xi_n\ :=\ \log \zeta_n\ .
\end{equation}  
By eqs.\ (\ref{eq7.1A1}) and (\ref{eq7.1A2}) and assumptions (a) and (b) in 
the statement of Theorem \ref{thm3.3}, one has that
\begin{equation}
\label{eq7.1A3}
\xi_n < 0\ \ {\rm for\ every}\ n \in \N;  \indent {\rm and} \indent
\lim_{n \to \infty} \xi_n = -\infty \quad {\rm and} \quad
\lim_{n \to \infty}  \xi_n / n\ =\ 0.
\end{equation}
[Note that for every $c > 0$, by assumption (b) in  
the statement of Theorem \ref{thm3.3},
one has that $\xi_n > -cn$ for all $n \in \N$ sufficiently large (depending on $c$),
forcing $\liminf_{n \to \infty} \xi_n/n \geq -c$. 
That and the (``first'' or) ``second third'' of (\ref{eq7.1A3}) together imply the
``final third'' of (\ref{eq7.1A3}).]
\medskip

     {\bf Sub-step 1B.}\ \ Some material in this section will have aspects that are
redundant or superfluous, in the hope of allowing a (trivially) slightly 
less ``cluttered'' presentation of arguments later on.  
\smallskip     
     
     Define the quantity
\begin{equation}
\label{eq7.1B1} 
{\bf r}\ :=\ \inf_{n \in \N} \xi_n / n.
\end{equation}
By (the ``first and last thirds'' of) (\ref{eq7.1A3}), 
\begin{equation}
\label{eq7.1B2}
-\infty\ <\ {\bf r}\ <\ 0.
\end{equation}
That is, {\bf r} is a (finite) negative real number.
\smallskip     
     
     For each ordered pair $(q,r)$ of real numbers such that
$q \leq 0$ and ${\bf r} \leq r < 0$, let $L_{(q,r)}: \R \to \R$ denote the 
affine function (``line'') defined as follows: 
\begin{equation}
\label{eq7.1B3}
{\rm for\ every}\ x \in \R,\ \ \ L_{(q,r)}(x)\ :=\ q + rx. 
\end{equation}

   Let ${\bf S}$ denote the set of all ordered pairs $(q,r)$ of real numbers with the 
following properties:
\begin{equation}
\label{eq7.1B4}
q \leq 0; \indent {\bf r} \leq r < 0; \indent 
{\rm and\ for\ every}\ n \in \N,\ \ L_{(q,r)}(n)\ \leq\ \xi_n.
\end{equation}
Note that $L_{(0,{\bf r})}(n) = {\bf r}n \leq \xi_n$ for every $n \in \N$
by (\ref{eq7.1B1}); and hence $(0, {\bf r}) \in {\bf S}$ by 
the entire sentence of (\ref{eq7.1B4}).
Hence the set ${\bf S}$ is nonempty.
\medskip

     {\bf Sub-step 1C.}\ \ Refer again to (\ref{eq7.1B1}) and (\ref{eq7.1B2}). 
For any given ordered pair $(q,r) \in {\bf S}$, one has by
(\ref{eq7.1B3}) and (\ref{eq7.1B4}) that
(i) if $x \geq 0$, then $rx \leq 0$ and $L_{(q,r)} = q + rx \leq 0 + 0 = 0$, and 
(ii) if instead $x < 0$, then $rx > 0$ and hence
$L_{(q,r)} = q + rx \leq 0 + rx = |r| \cdot |x| \leq (-{\bf r}) \cdot (-x)$.
Hence for each $x \in \R$, one has that $\sup_{(q,r) \in {\bf S}} L_{(q,r)}(x) < \infty$.       
[That supremum is $\leq 0$ if $x \geq 0$, and $\leq (-x) \cdot (-{\bf r})$ if 
$x < 0$.]
\smallskip     
     
Accordingly, define the function $\phi: \R \to \R$ as follows:
\begin{equation}
\label{eq7.1C1}
{\rm For\ every}\ x \in \R,\ \ \ \phi(x)\ :=\ \sup_{(q,r) \in {\bf S}} L_{(q,r)}(x).
\end{equation}
From (\ref{eq7.1C1}) and the sentence in brackets at the end of the 
preceding paragraph, one has that $\phi(0) \leq 0$.
Also, by (\ref{eq7.1C1}) and the sentence right after eq.\ (\ref{eq7.1B4}), 
one has that $\phi(0) \geq L_{(0, {\bf r})}(0) = 0 + {\bf r}0 = 0$.
Hence 
\begin{equation}
\label{eq7.1C2}
\phi(0)\ =\ 0.
\end{equation}
More information about this function $\phi$
will be developed below. 
\medskip  

     {\bf Claim 1D.}\ \ {\it For every $r \in [{\bf r}, 0)$ there exists $q < 0$ such that
$(q,r) \in {\bf S}$.}
\medskip

     {\bf Proof.}\ \  This is elementary.  Suppose $r \in [{\bf r}, 0)$.
By (\ref{eq7.1A3}) (its ``final third''), there can exist at most finitely 
many positive integers $n$ (if any) such that $L_{(0,r)}(n) = rn > \xi_n$.
Hence for $q < 0$ sufficiently far below 0, one has that
$L_{(q,r)}(n) = q + rn \leq \xi_n$ for all $n \in \N$
and hence $(q,r) \in {\bf S}$ by the entire sentence of (\ref{eq7.1B4}).
Thus Claim 1D holds.
\medskip

     {\bf Claim 1E.}\ \ {\it The function $\phi$ has the following properties: \zhfb
(a) For every positive integer $n$, one has that $\phi(n) \leq \xi_n$
and (hence) $\exp(\phi(n)) \leq \zeta_n$. \zhfb
(b) The function $\phi$ is continuous and strictly decreasing on $\R$. \zhfb
(c) The function $\phi$ is convex on $\R$. \zhfb 
(d) One has that $\lim_{x \to \infty} \phi(x) = -\infty$. \zhfb
(e) One has that $\lim_{x \to \infty} [\phi(x)]/x = 0$.}
\medskip

     This statement is elementary (and has considerable redundancy).
We shall go ahead and give the proofs of these statements, in the order 
(a), (d), (e), (c), (b).
\medskip

     {\bf Proof.}\ \ {\bf Proof of (a).}\ \  
In statement (a), the first inequality  holds by (\ref{eq7.1C1}) and the entire sentence of (\ref{eq7.1B4}), and the second inequality then follows from (\ref{eq7.1A2}).
\smallskip

   {\bf Proof of (d).}\ \ By property (a) here and the ``middle third'' 
of (\ref{eq7.1A3}), one has that 
$\phi(n) \to -\infty$ as $n \to \infty$ (that is, along the positive {\it integers\/}).
\smallskip

     Also, for each ordered pair $(q,r) \in {\bf S}$, the affine function $L_{(q,r)}$ is
strictly decreasing (on $\R$), by (\ref{eq7.1B3}) and the entire sentence of (\ref{eq7.1B4}).
Hence trivially by (\ref{eq7.1C1}), the function $\phi$ is nonincreasing, thus monotonic
(on $\R$).
Combining that monotonicity with the observation in the preceding paragraph above,
one has that  
$\phi(x) \to -\infty$ as $x \to \infty$ (that is, along the positive {\it real numbers\/}).
Thus statement (d) holds.
\smallskip

     {\bf Proof of (e).}\ \ Suppose $r \in [{\bf r},0)$.
Applying Claim 1D, let $q < 0$ be such that $(q,r) \in {\bf S}$.
Then for each real number $x$, 
by (\ref{eq7.1C1}) and (\ref{eq7.1B3}), one has that
$\phi(x) \geq L_{(q,r)}(x) = q +rx$.
Hence $\liminf_{x \to \infty} [\phi(x)]/x \geq r$.
Since  $r \in [{\bf r},0)$ was arbitrary, one has that
$\liminf_{x \to \infty} [\phi(x)]/x \geq 0$. 
Hence by property (d), statement (e) holds. 
\smallskip

   {\bf Proof of (c).}\ \ Suppose $y$ and $z$ are real numbers such that $y < z$.
Suppose $0 <  t < 1$.
Our task is to show that
\begin{equation}
\label {eq7.1E1}
\phi\bigl((1-t)y + tz\bigl)\ \leq\ (1-t)\phi(y) + t\phi(z).
\end{equation}

   Suppose $\zeps > 0$.
By (\ref{eq7.1C1}), there exists an ordered pair $(q,r) \in {\bf S}$ such that 
the following holds:
\begin{align}
\nonumber
\phi \bigl((1-t)y + tz\bigl)\ \leq\ \zeps\ +\  L_{(q,r)} \bigl((1-t)y + tz\bigl)\
&=\ \zeps\ +\ (1 - t)\, L_{(q,r)} (y)\ +\ t\, L_{(q,r)} (z)\\
\nonumber 
&\leq\ \zeps\ + (1 - t)\, \phi (y)\ +\ t\, \phi (z).    
\end{align}
Since $\zeps > 0$ was arbitrary, (\ref{eq7.1E1}) holds.
Thus statement (c) holds.
\smallskip

   {\bf Proof of (b).}\ \ The continuity of the function $\phi$ on $\R$ is a standard 
elementary consequence of convexity (property (c)).
The fact that the function $\phi$ is strictly decreasing on $\R$ is an elementary
consequence of properties (c) and (d).
Thus statement (b) holds.
That completes the proof of Claim 1E. 
\medskip

   {\bf Sub-step 1F.}\ \
By (\ref{eq7.1C2}) and Claim 1E, the function $\phi$ [say restricted to the 
closed half line $[1,\infty)$] satisfies the assumptions in Lemma \ref{lem6.5}. 
Refer to Notations \ref{nt6.4}(A)(B). 
Let $\Gamma$ denote the set of all numbers $x \in (1, \infty)$ (the open half line)
such that $\phi$ fails to be differentiable at the point $x$.  
As noted in Lemma \ref{lem6.5}(I)(iv), $\phi'(y) \to 0$ as 
$y \to \infty,\ y \notin \Gamma$.
\smallskip
       
For each $y \in (1, \infty) - \Gamma$ let $L^{(y)}: \R \to \R$ 
denote the affine function whose graph is tangent to
the graph of the (convex) function $\phi$ at the point
$(y, \phi(y))$ (as in the statement of Lemma \ref{lem6.5}).
Then for each $y \in (1, \infty)-\Gamma$, one has that 
[also repeating Lemma 6.5(I)(ii)]
\begin{align}
\label{eq7.1F1}
&L^{(y)}(y) = \phi(y) 
\indent {\rm and} \indent
L^{(y)}(x)\ \leq\ \phi(x)\ \ {\rm for\ all}\ x \in [1, \infty)\ 
[{\rm in\ fact\ for\ all}\ x \in \R];\ \quad {\rm and} \\
\label{eq7.1F2}
&[{\rm slope\ of}\ L^{(y)}]\ =\ \phi'(y)\ <\ 0.
\end{align}

   {\bf Step 2.  The parameters.}\ \
In this step, we shall recursively define for each 
positive integer $j$ nine parameters
(one of those ``parameters'' being an affine function). 
They are listed here, along with their most basic restrictions: 
\medskip

    {\bf Key List:} \zhfb
(i) a positive number $B_j$; \zhfb
(ii) a number $\zeps_j^* \in (0, 1/9]$; \zhfb 
(iii) a number $t_j \in (1, \infty)$; \zhfb
(iv) an affine function ${\bf L}_j: \R \to \R$; \zhfb
(v) a number $\zeps_j \in (0, 1/9]$; \zhfb
(vi) a number $\theta_j \in (0, 1/9]$; \zhfb
(vii) a number $\theta^*_j \in (0, 1/8]$; \zhfb  
(viii) a positive integer $I_j$; and \zhfb
(ix) a positive number $h_j$.
\medskip

   Such parameters $\zeps_j$, $\theta_j$, $I_j$, and $h_j$ 
in (v)-(vi) and (viii)-(ix) will also be defined for $j=0$.
Those four parameters for $j=0$ will provide a convenient ``springboard'' 
to start the recursion (but they will not otherwise play a significant
role later on in the construction itself for Theorem \ref{thm3.3}). 
\smallskip

     In superscripts, the numbers $t_j$ in (iii) for $j \geq 1$ 
will be denoted $t(j)$ for typographical convenience.
 \medskip

   {\bf Sub-step 2A.  The initial step.}\ \
We start with $j=0$, for just four parameters as mentioned above.
Define the following four parameters:
\begin{equation}
\label{eq7.2A1}
\zeps_0\ :=\ \theta_0\ :=\ 1/9  \quad {\rm and} \quad
I_0\ :=\ h_0:\ =\ 1.
\end{equation}
Obviously in the Key List above, the relevant basic restrictions 
in (v)-(iv) and (xiii)-(ix) for the four parameters defined in 
(\ref{eq7.2A1}) are satisfied.
That completes the initial step.
\medskip

   {\bf Sub-step 2B.  The recursion step.}\ \ 
Now suppose $j \geq 1$ is an integer, and that for each 
integer $u$ such that $0 \leq u \leq j-1$, the four parameters 
$\zeps_u$, $\theta_u$, $I_u$, and $h_u$ 
have been defined, meeting the basic requirements in (v)-(vi) and (xiii)-(ix) 
in the Key List:
$\zeps_u \in (0, 1/9]$, $\theta_u \in (0, 1/9]$, 
$I_u$ is a positive integer, and $h_u > 0$.
Those properties will repeatedly be tacitly applied below.  
Our task now is to define, for the given $j$ here, all nine parameters in
the Key List, meeting the basic requirements there.
\medskip   

   To start off,  define the number $B_j$ by
\begin{equation}
\label{eq7.2B1}
B_j\ :=\  j \cdot \sum_{u=0}^{j-1} h_u^2 \zeps_u / \theta_u.
\end{equation}
Of course $B_j > 0$ by the preceding paragraph.
 \smallskip 
 
   Next, let $\zeps_j^* \in (0, 1/9]$ be such that
\begin{equation}
\label{eq7.2B2} 
\zeps_j^*\ \leq\ \min \biggl\{\frac {\zeps_{j-1}} {2}\, ,\ 
\frac {9^{-j}} {I_{j-1}}\, ,\ \frac {2^{-j}} {9h_{j-1}^2} \biggl\}.
\end{equation}

   Next, we apply Lemma \ref{lem6.5}(II) to the (convex)
function $\phi$ here.
Of course $\log \zeps_j^* < 0$.   
Referring to Sub-step 1F and the notation $T(\dots)$ in Lemma \ref{lem6.5}(II),  
let $t_j \in (1,\infty) - \Gamma$ be a positive number such that 
\begin{equation}
\label{eq7.2B3}
t_j\ \geq\ T\bigl(\phi,\ \log \zeps_j^*,\ \min\{1/9,\, \theta_{j-1}/2,\, 2^{-j}/B_j \}\bigl).
\end{equation}

   Now referring to the entire paragraph of  (\ref{eq7.1F1})-(\ref{eq7.1F2}), 
define the affine function ${\bf L}_j: \R \to \R$ as follows:
\begin{equation}
\label{eq7.2B4}
{\rm For\ all}\ x \in \R,\ \ \ {\bf L}_j(x)\ :=\ -(j+2)\ +\  L^{(t(j))}(x). 
\end{equation} 

   Note that ${\bf L}_j(0) < L^{(t(j))}(0) \leq \log \zeps^*_j$ from (\ref{eq7.2B4}),
(\ref{eq7.2B3}), and Lemma \ref{lem6.5}(II).
It follows that  
$\exp({\bf L}_j(0)) < \zeps^*_j \leq 9^{-j}/I_{j-1} \leq 9^{-j} \leq 1/9$ by (\ref{eq7.2B2}).
Note also that from (\ref{eq7.2B3}) and (again) Lemma 6.5(II),
$-1/9 \leq  \phi'(t_j) < 0$.
Accordingly, define the positive numbers 
$\zeps_j \in (0, 1/9]$ and $\theta_j \in (0, 1/9]$; and then define the positive
number $\theta_j^*$ [\`a la (\ref{eq5.112})], the positive integer $I_j$, 
and [see also (\ref{eq7.2B1}) again] the positive number $h_j$; as follows:    
\begin{align}
\label{eq7.2B5}    
&\zeps_j\ :=\ \exp\bigl({\bf L}_j(0)\bigl) \quad {\rm and} \quad
\theta_j\ := -\phi'(t_j) \quad {\rm and} \quad
\theta_j^*\ :=\ \frac {\theta_j} {1 - \zeps_j}\ ; 
\quad {\rm and}\\
\label{eq7.2B6}
&I_j\ :=\ \Biggl[ \frac {1} {\theta_j^* \zeps_j} \Biggl] 
\indent {\rm and} \indent
h_j\ :=\ \max \biggl\{ 3h_{j-1}\, ,\  \Bigl( B_j \theta_j / \zeps_j \Bigl)^{1/2}\, \biggl\},     
\end{align}
where in the first equality in (\ref{eq7.2B6}), the big brackets indicate the
greatest integer that is $\leq 1/(\theta_j^* \zeps_j)$.
By simple arithmetic, $\theta^*_j \in (0, 1/8]$ and $I_j \geq 72$. 
Obviously $h_j > 0$.
[In fact both terms in the big braces in (\ref{eq7.2B6}) are positive.]\ \  
All nine parameters defined here in Sub-step 2B for the given integer $j$ 
meet the requirements specified in the Key List.
That completes the recursion step.
\smallskip

   That completes the recursive definition of the parameters. 
\medskip

   {\bf Sub-step 2C.  Miscellaneous observations.}\ \ 
Here are several observations that arise from, or for convenient later reference are
repeated from, that recursive definition. 
They will be employed later on.
\medskip

   First, from (\ref{eq7.2A1}), followed by (\ref{eq7.2B5}), 
the (first or) second sentence after (\ref{eq7.2B4}), and then 
(\ref{eq7.2B2}), one has that  
\begin{equation}
\label{eq7.2C1}
\zeps_0 = 1/9 \indent 
{\rm and\ for\ every}\ j  \in \N,\ \ 
\zeps_j\, =\, \exp({\bf L}_j(0))\, <\, \zeps^*_j\, \leq\, 
\min \biggl\{\frac {\zeps_{j-1}} {2}\, ,\ 
\frac {9^{-j}} {I_{j-1}}\, ,\ \frac {2^{-j}} {9h_{j-1}^2} \biggl\}.
\end{equation}

   Next, from (\ref{eq7.2B3}) and Lemma \ref{lem6.5}(II), for each $j \in \N$, 
$-\min\{ \theta_{j-1}/2,\, 2^{-j}/B_j\} \leq \phi'(t_j) < 0$.   Hence by (\ref{eq7.2B5}), 
\begin{equation}
\label{eq7.2C2}
{\rm for\ every}\ j \in \N,\ \ \ \theta_j\ \leq\ \min\{ \theta_{j-1}/2,\ 2^{-j}/B_j \}.
\end{equation}

   Next, by (\ref{eq7.2B6}) and (\ref{eq7.2B1}) and trivial arithmetic,
\begin{equation}
\label{eq7.2C3}
{\rm for\ every}\ j \in \N, \quad 
h_j^2 \zeps_j/\theta_j\ \geq\ B_j\ =
j \cdot \sum_{u=0}^{j-1} h_u^2 \zeps_u/\theta_u.
\end{equation}

    Next, from (\ref{eq7.2A1}) and (\ref{eq7.2B6}),
\begin{equation}
\label{eq7.2C4}
h_0 = 1\ \ \ 
{\rm and\ for\ every}\ j \in \N,\ \ \
h_j\ \geq\ 3h_{j-1}\ .
\end{equation}
        
Our final task here in Sub-step 2C is to show that
\begin{equation}
\label{eq7.2C5}
{\rm for\ every}\ j \in \N, \quad h_j^2 \zeps_j\ \leq\ 2^{-j}.  
\end{equation}
Suppose $j \in \N$.
Refer again to the ``second half'' of (\ref{eq7.2B6}).
If $h_j = 3 h_{j-1}$, then by (\ref{eq7.2C1}), 
\begin{equation}
\nonumber
h_j^2 \zeps_j\ =\ 9h_{j-1}^2 \zeps_j\  
\leq\ 9h_{j-1}^2 \cdot \frac {2^{-j}} {9h_{j-1}^2}\ =\ 2^{-j}  
\end{equation}
and thus (\ref{eq7.2C5}) is satisfied.
If (instead) $h_j = (B_j \theta_j / \zeps_j)^{1/2}$, then by (\ref{eq7.2C2}),
$h_j^2 \zeps_j = B_j \theta_j \leq 2^{-j}$ and thus (again) (\ref{eq7.2C5}) is satisfied.
Since $j \in \N$ was arbitrary, eq.\ (\ref{eq7.2C5}) holds (in either case) for every $j \in \N$.
\medskip

    {\bf Sub-step 2D.\ \ A summary of key properties of some of the parameters:}
 
      Justifications for the following comments will be briefly reiterated below.
\smallskip
  
\noindent (a) For each $j \in \N$, $\zeps_j \in (0, 1/9]$, $\theta_j \in (0, 1/9]$,
$\theta^*_j \in (0, 1/8]$, and $I_j \geq 72$. \zhfb 
(b) As $j \to \infty$, $\zeps_j \to 0$ monotonically
and $\theta_j \to 0$ monotonically. \zhfb
(c)  As $j \to \infty$, $\theta^*_j \to 0$ monotonically
and $I_j \to \infty$ monotonically. \zhfb  
(d) For each $j \in \N$, one has that $\theta_j/\theta^*_j < 1$ and
$\theta_j^*\zeps_j I_j \leq 1$.  
\zhfb
(e) As $j \to \infty$, one has that $\theta_j/\theta^*_j \to 1$ and 
$\theta_j^*\zeps_j I_j \to 1$ and that $\theta_j/h_j \to 0$. \zhfb
(f) One has that $\sum_{j=1}^\infty \zeps_j \leq 1/8$ and that
$\sum_{j=1}^\infty h_j^2 \zeps_j \leq 1$. \zhfb 
(g) For each $j \geq 2$, one has that $\zeps_j \leq 9^{-j}/I_{j-1}$ and that
$h_j^2 \zeps_j / \theta_j \geq B_j 
= j \cdot [1 + \sum_{u=1}^{j-1} (h_u^2 \zeps_u / \theta_u)] > j$. 
\medskip

   Property (a) holds by the entire sentence of (\ref{eq7.2B5})--(\ref{eq7.2B6})
and the sentence after it. 
Property (b) holds by (\ref{eq7.2C1}) and (\ref{eq7.2C2}). 
Property (c) holds by properties (a) and (b) 
and (\ref{eq7.2B5})--(\ref{eq7.2B6}) and simple arithmetic.
The two inequalities in (d) hold by (\ref{eq7.2B5})--(\ref{eq7.2B6})
and property (a). 
In property (e), the first two limits hold by properties (b) and (c) and
eqs.\ (\ref{eq7.2B5})--(\ref{eq7.2B6}); and the third limit holds by
property (a) [$0 < \theta_j \leq 1/9$ for all $j \in \N$] and the fact from 
(\ref{eq7.2C4}) that $h_j \to \infty$ as $j \to \infty$.
In (f), the first inequality holds since $\zeps_j \leq 9^{-j}$ for each 
$j \in \N$ by (\ref{eq7.2C1}); and the second inequality holds by 
(\ref{eq7.2C5}).
In (g), the first inequality holds by (\ref{eq7.2C1}), and the rest 
holds by (\ref{eq7.2C3}) together with the fact from (\ref{eq7.2A1})
that $h_0^2 \zeps_0/\theta_0 = 1$.
\medskip

   {\bf Step 3.  The construction.}\ \
This ``Step'' will be divided into several ``sub-steps''.
\medskip

   {\bf Sub-step 3A.}\ \ For each positive integer $j$, 
referring to the parameters $\zeps_j$ and $\theta_j$ defined in 
Sub-step 2B [see Sub-step 2D(a)], let
$X^{(j)} := (X^{(j)}_k, k \in \Z)$ be a strictly stationary
Markov chain that satisfies Condition 
${\cal H}(\zeps_j, \theta_j)$
of Definition \ref{def5.2}.
Let this construction be carried out in such a way that these
(strictly stationary) Markov chains 
$X^{(1)},\ X^{(2)},\ X^{(3)},\ \dots$
are independent of each other.
\medskip

  {\bf Sub-step 3B.}\ \ By Definition \ref{def5.2} and
Lemma \ref{lem5.3}, for each $j \in \N$ and each $k \in \Z$,
the random variable $X^{(j)}_k$ takes its values in the set
$\{-1,0,1\}$, with probabilities
$P(X^{(j)}_k = 0) = 1 - \zeps_j$ and
$P(X^{(j)}_k = -1) = P(X^{(j)}_k = 1) = \zeps_j/2$, and 
hence $P(X^{(j)}_k \neq 0) = \zeps_j$.
Hence by Sub-step 2D(f), one has that
\begin{equation}
\label{eq7.3B1}
{\rm for\ each}\ k \in \Z, \quad \sum_{j=1}^\infty P(X^{(j)}_k \neq 0)\ 
\leq\ 1/8\ <\ \infty.
\end{equation}
Hence for each $k \in \Z$, by the Borel-Cantelli Lemma, one has that
\begin{equation}
\label{eq7.3B2}
P\Bigl(X^{(j)}_k \neq 0\ {\rm for\ infinitely\ many}\
j \in \N \Bigl)\ =\ 0.
\end{equation}
Deleting a set of probability 0 from the probability space
$\Omega$ is necessary, we assume without loss of generality
that for each $k \in \Z$, the event in the left side of
(\ref{eq7.3B2}) is empty.
\medskip     

 {\bf Sub-Step 3C.}\ \ 
Referring to the second equality in (\ref{eq7.2B6}), define the sequence 
$X := (X_k, k \in \Z)$ of random variables as follows:   
For each $k \in \Z$ and each $\omega \in \Omega$, 
\begin{equation}
\label{eq7.3C1}
X_k(\omega)\ :=\ \sum_{j=1}^\infty h_j X^{(j)}_k(\omega)\ . 
\end{equation}  
For each $k \in \Z$ and each $\omega \in \Omega$, 
by the last sentence of Sub-step 3B, 
the sum in (\ref{eq7.3C1}) has at most finitely
many non-zero terms.
\medskip

  {\bf Sub-step 3D.}\ \ 
Refer again to the first and last sentences of Sub-step 3B.
\smallskip

Let $\Upsilon$ denote the set of all
sequences $y := (y_1, y_2, y_3, \dots)$ of elements of 
the set $\{-1,0,1\}$ such that $y_j \neq 0$ for at most
finitely many indices $j$.
Let $\eta: \Upsilon \to \R$ denote the function defined by
$\eta(y) := \sum_{j=1}^\infty h_jy_j$
for $y := (y_1, y_2, y_3, \dots) \in \Upsilon$.
We use the notation ${\bf 0} := (0,0,0,\dots)$ for the 
``zero-sequence'' in $\Upsilon$.
Then $\eta({\bf 0}) = 0$.
\vskip 0.1 in

    Let ${\bf S}$ denote the range of the function $\eta$.
One has that (i) the set $\Upsilon$ is countably infinite,
(ii) the function $\eta$ is one-to-one by
(\ref{eq7.2C4}) and Lemma \ref{lem6.3}, hence
(iii) the set ${\bf S}$ is countably infinite, and finally
(iv) $0 \in {\bf S}$ (since $\eta({\bf 0}) = 0$).  
\vskip 0.1 in   

   By (\ref {eq7.3C1}), for each $k \in \Z$, 
\begin{equation}
\label{eq7.3D1}
X_k\ =\ 
\eta\Bigl( (X^{(1)}_k, X^{(2)}_k, X^{(3)}_k, \dots) \Bigl).
\end{equation}
That is, for each $k \in \Z$,  the random variable $X_k$ takes its values in the
set ${\bf S}$.
Also, as a consequence of property (ii) at the end of the preceding paragraph, 
one has, for each $k \in \Z$ and each element 
$y := (y_1, y_2, y_3, \dots) \in \Upsilon$, the following equality of events: 
\begin{equation}
\label{eq7.3D2}
\bigl\{X_k = \eta(y)\bigl\}\ =\ 
\bigcap_{j=1}^\infty \bigl\{X^{(j)}_k = y_j\bigl\}. 
\end{equation}
 
     Suppose $k \in \Z$ and suppose $y := (y_1, y_2, y_3, \dots)$
is any element of $\Upsilon$.
Let us apply Lemma \ref{lem6.6} with the events there being 
$E_j := \{X^{(j)}_k = 0\}$ and $F_j := \{X^{(j)}_k = y_j\}$.
By (\ref{eq7.3B1}), the first sentence of Sub-step 3B, the independence
condition in Sub-step 3A, and the definition of the set $\Upsilon$,
all of the hypothesis of Lemma \ref{lem6.6} is satisfied.  
By that lemma,      
\begin{equation}
\label{eq7.3D3}
P\Bigl(X_k = \eta(y)\Bigl)\ 
=\ \prod_{j=1}^\infty P\Bigl(X^{(j)}_k = y_j\Bigl) \ >\ 0. 
\end{equation}
What has been shown here is that for every 
$y \in \Upsilon$ (and every $k \in \Z$), $P(X_k = \eta(y)) > 0$.
\medskip 

   Now recall again both sentences of Sub-step 3A, 
eq.\ (\ref{eq7.3C1}), the fact that 
the function $\eta$ is one-to-one,
and eqs.\ (\ref{eq7.3D1}), (\ref{eq7.3D2}), and (\ref{eq7.3D3}).
Applying these (considerably redundant) pieces together 
in an appropriate elementary argument, one has 
that the the random sequence $X$ is a 
strictly stationary Markov chain, with (countably infinite) 
state space $\bf S$, such that $P(X_0 = s) > 0$
for every $s \in {\bf S}$.
\medskip

   {\bf Sub-step 3E.}\ \
For any $j \in \N$ and any pair of (not necessarily distinct) 
elements $u, v \in \{-1,0,1\}$, one has, 
by Lemma \ref{lem5.3}(2) and Notations \ref{nt2.8}(B), that
\begin{equation}
\label{eq7.3E1}
P(\{X^{(j)}_0 = u\} \cap \{X^{(j)}_1 = v\})\
=\ P(\{X^{(j)}_1 = u\} \cap \{X^{(j)}_0 = v\}).
\end{equation}

   Now recall again the second sentence of Sub-step 3A.
By (\ref{eq7.3D2}) [its entire sentence], 
for any pair of (not necessarily distinct) elements $y := (y_1, y_2, y_3, \dots)$ 
and $z := (z_1, z_2, z_3, \dots)$ of $\Upsilon$,
one has that
\begin{align}
\label{eq7.3E2}
P\Bigl(\{X_0 = \eta(y)\} \cap \{X_1 = \eta(z)\}\Bigl)\ 
&=\ P\biggl(\,
\biggl[\, \bigcap_{j=1}^\infty \{X^{(j)}_0 = y_j\} \biggl]\,
\bigcap\, 
\biggl[\, \bigcap_{j=1}^\infty \{X^{(j)}_1 = z_j\} \biggl]\, 
 \biggl) \nonumber \\
&=\ \prod_{j=1}^\infty
P\Bigl(\{X^{(j)}_0 = y_j\} \cap \{X^{(j)}_1 = z_j\}\Bigl)   
\end{align}
By (\ref{eq7.3E1}), the final product in (\ref{eq7.3E2})
remains unchanged if for each $j \in \N$ the numbers 
$y_j$ and $z_j$ are interchanged with each other.
Hence as a consequence of the entire sentence of (\ref{eq7.3E2}) itself,  
it follows that the first term in (\ref{eq7.3E2}) remains unchanged
if the elements $y$ and $z$ are interchanged.
That is, for any two elements $q$ and $s$ of the state 
space ${\bf S}$, one has that
$P(\{X_0 = q\} \cap \{X_1 = s\})
= P(\{X_0 = s\} \cap \{X_1 = q\})$.
Since the state space ${\bf S}$ of the (strictly stationary) 
Markov chain $X$ is discrete, it follows, by Notations \ref{nt2.8}(B) and 
a trivial argument, that this Markov chain $X$ is reversible.
\medskip

   {\bf Sub-step 3F.}\ \
Refer again to the first three paragraphs of Sub-step 3D, including in particular
the equality $\eta({\bf 0}) = 0$.
\smallskip

   Suppose $s$ is any element of the state space ${\bf S}$.     
Let $y := (y_1,y_2,y_3,\dots)$ be the element of $\Upsilon$
such that $\eta(y) = s$.
By (\ref{eq7.3E2}), 
\begin{equation}
\label {eq7.3F1}
P\Bigl(\{X_0 = 0\} \cap \{X_1 = s\}\Bigl)\ =\ 
\prod_{j=1}^\infty P\Bigl(\{X^{(j)}_0 = 0\} \cap P(\{X^{(j)}_1 = y_j\}\Bigl). 
\end{equation}
We shall apply Lemma \ref{lem6.6} with the events
$E_j$ and $F_j$, $j \in \N$ there given by
\begin{equation}
\label{eq7.3F2}
E_j\ :=\ \{X^{(j)}_0 = 0\} \cap P(\{X^{(j)}_1 = 0\}
\indent {\rm and} \indent
F_j\ :=\ \{X^{(j)}_0 = 0\} \cap P\{X^{(j)}_1 = y_j\}.
\end{equation}
For each $j \in \N$, 
$P(E_j^c) \leq P(X^{(j)}_0 \neq 0) + P(X^{(j)}_1 \neq 0) \leq 2\zeps_j$
by (again) the first sentence of Sub-step 3B.
Hence $\sum_{j=1}^\infty P(E_j^c) \leq 1/4 < \infty$ by Sub-step 2D(f).
Next, for each $j \in \N$, $y_j \in \{-1,0,1\}$ (by the definition of $\Upsilon$) 
and hence $P(F_j )  > 0$ by Lemma \ref{lem5.3}(3b)  
and the first sentence of Sub-step 3A.
Next, by the definition of $\Upsilon$, one has that $E_j = F_j$ for all 
except at most finitely many indices $j \in \N$.
Finally, refer again to the second sentence (independence) of Sub-step 3A. 
It now follows from Lemma \ref{lem6.6} that $\prod_{j=1}^\infty P(F_j) > 0$.
Thus by (\ref{eq7.3F2}) and (\ref{eq7.3F1}), 
$P(\{X_0 = 0\} \cap \{X_1 = s\}) > 0$.
\vskip 0.1 in

   What has been shown is that for every $s \in {\bf S}$ (the state space),
$P(\{X_0 = 0\} \cap \{X_1 = s\}) > 0$.
Hence by reversibility (proved in Sub-step 3E above), one has that
for every $s \in {\bf S}$,
$P(\{X_0 = s\} \cap \{X_1 = 0\}) = P(\{X_0 = 0\} \cap \{X_1 = s\}) > 0$.
As a consequence, for every $s \in {\bf S}$, 
$P(X_1 = s | X_0 = 0) > 0$ and $P(X_1 = 0 | X_0 = s) > 0$. 
The last two inequalities have two key consequences:
\smallskip

     First (by strict stationarity), for any two (not necessarily distinct) elements
$s, t \in {\bf S}$,
$P(X_2 = t | X_0 = s) \geq P(X_1 = 0 | X_0 = s) \cdot P(X_2 = t | X_1 = 0) > 0$.
Hence the Markov chain $X$ is irreducible.
\smallskip

     Second (take $s = 0$), $P(X_1 = 0 | X_0 = 0) > 0$.
Hence the Markov chain $X$ is aperiodic.
\medskip 

   {\bf Sub-step 3G.}\ \ To conclude Step 3, let us compile here 
for convenience reference the 
last sentence from each of Sub-step 3D, Sub-step 3E, 
and each of the last two paragraphs of Sub-step 3F:
\smallskip

   {\it The random sequence $X$ defined in (\ref {eq7.3C1})
is a strictly stationary, countable-state, irreducible, aperiodic,
reversible Markov chain, with (countably infinite) state 
space {\bf S} [with $P(X_0 = s) > 0$ for every $s \in {\bf S}$]\/.} 
\medskip

   {\bf Step 4.  Proof of statement (i)} [in Theorem \ref{thm3.3}]. 
Refer to Sub-step 3A and the first sentence of Sub-step 3B.
Of course for any given $j \in \N$ and any given 
$k \in \Z$, one has that $EX^{(j)}_k = 0$ and $\var X^{(j)}_k = \zeps_j$
and hence also $E(h_jX^{(j)}_k) = 0$ and 
$E((h_jX^{(j)}_k)^2) = \var(h_jX^{(j)}_k) = h_j^2 \zeps_j$.
By Sub-step 2D(f), $\sum_{j=1}^\infty h_j^2\zeps_j < \infty$.
Hence for each $k \in \N$, by Lemma \ref{lem6.2} and the
second sentence of Sub-step 3A, the random variable  $X_k$ in
(\ref{eq7.3C1}) is square-integrable and has mean 0 and 
finite second moment (variance) $\sum_{j=1}^\infty h_j^2 \zeps_j$.
Thus statement (i) in Theorem \ref{thm3.3} holds.
\medskip

    {\bf Step 5.  Proof of statement (ii)} [in Theorem \ref{thm3.3}].\ \
For each $j  \in \N$, one has the following:
First, of course $t_j \in (1,\infty) - \Gamma$ by the entire sentence of
(\ref{eq7.2B3}).  
Next, by (\ref{eq7.2B4}), (\ref{eq7.1F2}), and (\ref{eq7.2B5}), 
$[{\rm slope\ of}\ {\bf L}_j] = [{\rm slope\ of}\ L^{(t(j))}] = \phi'(t_j) = -\theta_j$,  
and hence for every $x \in \R$, ${\bf L}_j(x) = {\bf L}_j(0) - \theta_jx$.
Hence by (\ref{eq7.2B5}), for every $x \in \R$, 
\begin{equation}
\label{eq7.51}
\exp\bigl({\bf L}_j(x)\bigl)\ =\ 
\exp\bigl({\bf L}_j(0)\bigl)\, \cdot\, \exp(-\theta_j x)\
=\ \zeps_j \cdot \exp(-\theta_j x).  
\end{equation}
Also, by (\ref{eq7.2B4}) again, for every $x \in \R$,
$\exp\bigl({\bf L}_j(x)\bigl)\ 
=\ \exp\bigl( - (j+2) \bigl)\, \cdot\, \exp \bigl(L^{(t(j))}(x) \bigl)$;
and hence by (\ref{eq7.1F1}), one has that 
for every $x \in [1, \infty)$,
\begin{equation}
\label{eq7.52}
\exp\bigl({\bf L}_j(x)\bigl)\ 
\leq\ \exp\bigl(-(j+2) \bigl)\, \cdot\, \exp \bigl(\phi(x) \bigl).
\end{equation}
Finally, substituting (\ref{eq7.51}) into (\ref{eq7.52}), 
one has that every $x \in [1, \infty)$,
\begin{equation}
\label{eq7.53}
\zeps_j \cdot \exp(-\theta_j x)\ \leq\ e^{-j} \cdot e^{-2} \cdot \exp \bigl(\phi(x) \bigl).
\end{equation}

   What has been shown is that (\ref{eq7.53}) holds for 
every $j \in \N$ and every $x \in [1, \infty)$.
Now by Lemma \ref{lem5.3}(6) and Sub-step 3A, and then (\ref{eq7.53}),
and then Claim 1E(a) in Step 1, one has that
for every $j \in \N$ and every $n \in \N$,
\begin{equation}
\nonumber
\beta_{X(j)}(n)\ \leq\ 6\zeps_j \cdot \exp(-\theta_j n)\
<\ e^{-j} \cdot \exp \bigl(\phi(n) \bigl)\ \leq\ e^{-j} \zeta(n).
\end{equation}
Hence by Lemma \ref{lem6.1}, for every $n \in \N$,
$\beta_X(n) \leq \sum_{j=1}^\infty \beta_{(X(j)}(n) 
\leq \sum_{j=1}^\infty e^{-j} \zeta(n) < \zeta(n)$. 
Referring to (\ref{eq2.34}), one now has that 
statement (ii) in Theorem  \ref{thm3.3} holds.
\medskip

   {\bf Step 6.  Proof of Statement (iii)} [in Theorem \ref{thm3.3}].\ \ 
Eq.\ (\ref{eq3.11}) [which is Statement (i) in the theorem] 
was verified in Step 4 above.
As noted near the end of Remark \ref {rem3.1},
(\ref {eq3.12}) then follows from (\ref {eq3.13}).
To complete the proof of Statement (iii), it now suffices 
to verify (\ref {eq3.13}).
The rest of Step 6 here will be devoted to that task.
\medskip

     {\bf Proof of (\ref {eq3.13}).}\ \ Refer to the final
paragraph of Remark \ref{rem3.1}.
It suffices to show that the equality in (\ref{eq3.13})
holds for the case $b=1$. 
Let $\gamma > 0$ be arbitrary bit fixed.
It suffices to show that
\begin{equation}
\label{eq7.61}
\limsup_{n \to \infty}\, \biggl[\sup_{r \in \bbR}\, 
P\biggl(r-1\ <\ n^{-1/2}\sum_{k=1}^nX_k\ <\ r+1\biggl)
\biggl]\ 
\leq \gamma.    
\end{equation}

   Now let $j$ be a positive integer fixed sufficiently large
that $4/\gamma^2 \leq j$.
By Sub-step 2D(g), $j < h_j^2 \zeps_j/\theta_j$.
Thus (for later reference) 
\begin{equation}
\label{eq7.62}
2/\gamma\ \leq\ h_j (\zeps_j/\theta_j)^{1/2}.
\end{equation}
Define (indirectly) the positive number $\sigma$ (which will depend on the fixed
positive integer $j$) by
\begin{equation}
\label{eq7.63}
\sigma^2\ :=\ \zeps_j [(2/\theta_j) -1].
\end{equation}
Note for later reference the trivial fact that $1 < 1/\theta_j$,
hence $(2/\theta_j) - 1 > 1/\theta_j$, and hence
$\sigma^2 > \zeps_j/\theta_j$ by (\ref{eq7.63});
and hence by (\ref{eq7.62}), 
\begin{equation}
\label{eq7.64}
\frac {1} {h_j \sigma}\ <\ 
\frac {1} {h_j}\, (\theta_j/\zeps_j)^{1/2}\ 
\leq\ \gamma/2 
\end{equation}

   For convenience, let $Z$ be a $N(0,1)$ random variable.
Then of course $aZ$ is a $N(0, a^2)$ random variable for
any $a \neq 0$.
By Lemma \ref{lem5.3}(10), eq.\ (\ref{eq7.63}), and Sub-step 3A, 
$n^{-1/2} \sum_{k=1}^n X^{(j)}_k$ converges 
in distribution to $\sigma Z$ as $n \to \infty$.
Hence by Slutsky's Theorem,
$n^{-1/2} \sum_{k=1}^n h_jX^{(j)}_k$ converges in distribution 
to $h_j \sigma Z$ as $n \to \infty$.
Since the (cumulative) distribution function of a 
(non-degenerate) normal distribution is continuous, it of course follows 
(from e.g.\ [\cite {ref-journal-Bill}, p.\ 198, Exercise 14.8(a) and 
Exercise 14.8(c) (its first sentence)])
that as $n \to \infty$, \break
(i) $P(n^{-1/2} \sum_{k=1}^n h_j X^{(j)}_k = x)$ 
converges to 0 uniformly over $x \in \R$, and
(ii) $P(n^{-1/2} \sum_{k=1}^n h_j X^{(j)}_k \leq x)$ 
converges to $P(h_j \sigma Z \leq x)$ uniformly
over $x \in \R$.
It follows from just a little bit of simple standard arithmetic 
that as $n \to \infty$,
$P(r -1 < n^{-1/2} \sum_{k=1}^n h_j X^{(j)}_k < r+1)$
converges to $P(r - 1 <  h_j \sigma Z < r+1)$
uniformly over $r \in \R$.
Accordingly, let $N$ be a positive integer such that
\begin{equation}
\label{eq7.65}
{\rm for\ all}\ n \geq N,\ \ 
\sup_{r \in \R} \biggl|
P\biggl(r -1 < 
n^{-1/2} \sum_{k=1}^n h_j X^{(j)}_k < r+1\biggl)\
-\ P\Bigl(r - 1 <  h_j \sigma Z < r+1\Bigl) 
\biggl|\
\leq\ \gamma/2.
\end{equation}

    Now let $m$ be an arbitrary fixed integer such that 
$m \geq N$.
To complete the proof of (\ref{eq7.61})
[and hence that of (\ref{eq3.13})], 
it suffices to prove for this $m$ that
\begin{equation}
\label{eq7.66}
\sup_{r \in \R} 
P\biggl(r-1\ <\ m^{-1/2}\sum_{k=1}^m X_k\ <\ r+1\biggl)\
\leq\ \gamma.
\end{equation} 

   Now by (\ref{eq7.65}), for each $r \in \R$,
\begin{equation}
\label{eq7.67}
\biggl| P\biggl(r -1 < 
m^{-1/2} \sum_{k=1}^m h_j X^{(j)}_k < r+1\biggl)\
-\ P\Bigl(r-1 <  h_j \sigma Z <  r+1\Bigl) 
\biggl|\ \leq\ \gamma/2.
\end{equation}
Since the random variable $h_j \sigma Z$ is normal with 
(mean 0 and) standard deviation $h_j \sigma$,
its probability density function is bounded above by
the positive number $1/(\sqrt{2 \pi} \cdot h_j \sigma)$, 
which in turn is less than $1/(2 h_j \sigma)$.
For any given $r \in \R$,      
since the interval $(r-1,\, r+1)$ has length 2, one has 
by (\ref{eq7.64}) that
$P(r-1 <  h_j \sigma Z <  r+1) \leq 
1/(h_j \sigma) \leq \gamma/2$.
Hence by (\ref{eq7.67}), 
\begin{equation}
\label{eq7.68}
{\rm for\ each}\ r \in \R, \quad 
P\biggl(r -1 <  m^{-1/2} \sum_{k=1}^m h_j X^{(j)}_k < r+1\biggl)\
\leq\ \gamma.
\end{equation}

   Just for convenient notations, define the random 
variables $V$ and $W$ by
\begin{equation}
\label{eq7.69}
V\ :=\ m^{-1/2} \sum_{k=1}^m h_j X^{(j)}_k\ 
\quad {\rm and} \quad 
W\ :=\ m^{-1/2} \sum_{u \in \N - \{j\}}\, 
\sum_{k=1}^m h_u X^{(u)}_k. 
\end{equation}
Recall again Sub-step 3A (independence) as well as
the first and last sentences of Sub-step 3B.
Eq.\ (\ref{eq7.68}) simply says that 
$P(r-1 < V < r+1) \leq \gamma$ for each $r \in \R$.
Also, $m^{-1/2} \sum_{k=1}^m X_k = V + W$ 
by (\ref{eq7.3C1}) and (\ref{eq7.69}).
Further, the random variables $V$ and $W$ are independent,
and they are each discrete.
\smallskip

Let $T$ denote the set of all real numbers $w$ 
such that $P(W=w) > 0$.
By the preceding four sentences,
one has that for each real number $r$,
\begin{align}
P&\biggl(r-1 <  m^{-1/2}\sum_{k=1}^m X_k  < r+1\biggl)\
=\ P\Bigl(r-1 < V+W < r+1\Bigl) \nonumber \\ 
&=\ \sum_{w \in T} 
P\biggl(\{r-1 < V+W < r+1\} \cap \{W = w\}\biggl)\
=\ \sum_{w \in T} 
P\biggl(\{r-1 < V+w < r+1\} \cap \{W = w\}\biggl)
\nonumber \\
&=\ \sum_{w \in T} 
P\biggl(\{r-w-1 < V < r-w+1\} \cap \{W = w\}\biggl)
\nonumber \\ 
&=\ \sum_{w \in T} 
\Bigl[P\bigl(r-w-1 < V < r-w+1\bigl) \cdot P(W=w)\Bigl]\
\leq \sum_{w \in T} [\gamma \cdot P(W=w)]\
=\ \gamma \cdot 1.  \nonumber    
\end{align}
Thus (\ref{eq7.66}) holds.
That completes the proof of (\ref {eq3.13}),
and of statement (iii) in Theorem \ref{thm3.3}.
 \medskip

    {\bf Step 7.  Proof of Statement (iv)} [in Theorem \ref{thm3.3}].\ \
The argument will be divided into four ``sub-steps''.
\vskip 0.1 in

    {\bf Sub-step 7A.}\ \ 
Refer to Notations \ref{nt4.1}(A), including its last sentence.
Refer to Sub-step 2D(a) and the first sentence of Sub-step 3A.   
For each $j \in \N$, applying Lemma \ref{lem5.4}, let
$W^{(j)} := (W^{(j)}_1, W^{(j)}_2, W^{(j)}_3, \dots)$
be a sequence of independent, identically distributed,
discrete random variables, each having the probability function 
${\bf g}[\theta^*_j \zeps_j, \theta_j]$, 
such that
\begin{equation}
\label{eq7.7A1}
P \biggl(\, \sum_{k=1}^{I(j)} X^{(j)}_k\ \neq\ 
\sum_{u=1}^{I(j)} W^{(j)}_u \biggl)\
\leq\ 3\zeps_j,
\end{equation}
where here and below, the positive integer $I_j$ [see (\ref{eq7.2B6})] 
is also written as $I(j)$ for typographical convenience.
\vskip 0.1 in

Recall from Sub-step 2D(b)(c)(e) that as $j \to \infty$,
one has that $\theta^*_j \zeps_j \to 0$ and $\theta_j \to 0$,
and that $I_j \cdot \theta^*_j \zeps_j \to 1$.
By Lemma \ref{lem4.2}, 
$\theta_j\sum_{u=1}^{I(j)} W^{(j)}_u \to \mu_{P1sL}$ in 
distribution as $j \to \infty$.            
As a consequence of (\ref{eq7.7A1}) and Sub-step 2D(b), 
$(\theta_j \sum_{k=1}^{I(j)} X^{(j)}_k) 
- (\theta_j \sum_{u=1}^{I(j)} W^{(j)}_u)\, \to 0$
in probability as $j \to\infty$.
Hence by Slutsky's Theorem,
$\theta_j \sum_{k=1}^{I(j)} X^{(j)}_k\ \to\ \mu_{P1sL}$
in distribution as $j \to \infty$.
That is,
\begin{equation}
\label{eq7.7A2}
\frac {\theta_j} {h_j} \sum_{k=1}^{I(j)} h_j X^{(j)}_k\ \to\ \mu_{P1sL}\ 
{\rm in\ distribution\ as}\ j \to \infty. 
\end{equation}

     {\bf Sub-step 7B.}\ \ 
 By Sub-step 3A and Lemma \ref{lem5.3}(4)(9),
for each $u \in \N$ and each $n \in \N$,
the random variable $\sum_{k=1}^n X^{(u)}_k$ has mean 0 and 
second moment (variance) $< n \cdot 2\zeps_u/\theta_u$.    
Hence for each $j \geq 2$, by Sub-step 3A (independence),
both inequalities in Sub-step 2D(d), and 
then Sub-step 2D(g), the random variable 
$(\theta_j / h_j) \sum_{u=1}^{j-1} \sum_{k=1}^{I(j)} h_u X^{(u)}_k$
has mean 0, and its second moment (variance) satisfies
\begin{align}
\nonumber
E& \biggl[ \biggl(\frac {\theta_j} {h_j} 
\sum_{u=1}^{j-1} \sum_{k=1}^{I(j)} h_u X^{(u)}_k
\biggl)^2 \, \biggl]\
=\ \frac {\theta_j^2} {h_j^2} 
\sum_{u=1}^{j-1} \biggl[\, h_u^2\, 
\var \biggl(\, \sum_{k=1}^{I(j)} X^{(u)}_k \biggl) \biggl]\
\leq\ \frac {\theta_j^2} {h_j^2} 
\sum_{u=1}^{j-1} \Bigl(h_u^2 \cdot 
I_j \cdot 2\zeps_u/\theta_u \Bigl) \nonumber \\
\nonumber
&=\ \frac {\theta_j^2} {h_j^2} I_j 
\sum_{u=1}^{j-1} \Bigl(h_u^2 \cdot 
2\zeps_u/\theta_u \Bigl)\  
\leq\ \frac {\theta_j^2} {h_j^2} 
\cdot \frac {1} {\theta^*_j \zeps_j} 
\sum_{u=1}^{j-1} \Bigl(h_u^2 \cdot 
2\zeps_u/\theta_u \Bigl)\
\leq\ \frac {\theta_j} {h_j^2 \zeps_j}  
\sum_{u=1}^{j-1} \Bigl(h_u^2 \cdot 
2\zeps_u/\theta_u \Bigl)\
\leq\ \frac {2} {j}\ .  
\end{align}
As a consequence,
\begin{equation}
\label{eq7.7B1}
\frac {\theta_j} {h_j} 
\sum_{u=1}^{j-1} \sum_{k=1}^{I(j)} h_u X^{(u)}_k\
\to\ 0\ \, {\rm in\ probability\ as}\ j \to \infty.
\end{equation}

   {\bf Sub-step 7C.}\ \ 
If $j$ and $u$ are any positive integers such that $u \geq j+1$,
then $I_{u-1} \geq I_j$ by Sub-step 2D(c),
and hence by Sub-step 2D(g) and the first sentence
in Sub-step 3B one has that for each $k \in \Z$,  
$P(X^{(u)}_k \neq 0) = \zeps_u  \leq 9^{-u}/I_{u-1} \leq 9^{-u}/I_j$.
Hence for any positive integer $j$, 
\begin{align}
\nonumber
P\biggl( \frac {\theta_j} {h_j}
\sum_{u=j+1}^\infty \sum_{k=1}^{I(j)} 
h_uX^{(u)}_k \neq 0 \biggl)\
&\leq\ \sum_{u=j+1}^\infty \sum_{k=1}^{I(j)}
P(X^{(u)}_k \neq 0)\ 
\leq\ \sum_{u=j+1}^\infty \sum_{k=1}^{I(j)} 
(9^{-u}/I_j)\  \\
\nonumber
&=\ \sum_{u=j+1}^\infty 9^{-u}\
=\ 9^{-(j+1)}/(1 - 9^{-1})\ <\ 9^{-j}.   
\end{align}
Hence 
$
P\bigl((\theta_j/h_j) \sum_{u=j+1}^\infty \sum_{k=1}^{I(j)} 
h_uX^{(u)}_k \neq 0\bigl) \to 0
$ 
as $j \to \infty$.
As a consequence,
\begin{equation}
\label{eq7.7C1}
\frac {\theta_j} {h_j} \sum_{u=j+1}^\infty \sum_{k=1}^{I(j)} h_uX^{(u)}_k 
\to\ 0\, \ {\rm in\ probability\ as}\  j \to \infty.
\end{equation}

   {\bf Sub-step 7D.}\ \
Recall again (say) the last sentence of Sub-step 3B.
For each $j \geq 2$,  by (\ref{eq7.3C1}), the random variable 
$(\theta_j/h_j )\sum_{k=1}^{I(j)} X_k$
is the sum of the random variables on the left sides of
(\ref{eq7.7A2}), (\ref{eq7.7B1}), and (\ref{eq7.7C1}).    
       
Hence by (\ref{eq7.7A2}), (\ref{eq7.7B1}), (\ref{eq7.7C1}), and 
Slutsky's Theorem,  
$(\theta_j/h_j) \sum_{k=1}^{I(j)} X_k$
converges in distribution to $\mu_{P1sL}$ as $j \to \infty$.
Also, from Sub-step 2D(e), $\theta_j/h_j \to 0$ as $j \to \infty$.
By Sub-step 2D(c) (and a trivial further argument),
statement (iv) in Theorem \ref{thm3.3} holds.
That completes the proof of Theorem \ref{thm3.3}.

\section{Appendix}
\label{sc8}

   In this Appendix, (i) we retain all notations and notational conventions
from Section \ref{sc2} and Remark \ref{rem3.1}, 
but (ii) the notations from the rest of Section \ref{sc3}
and from Sections \ref{sc4}, \ref{sc5}, \ref{sc6}, and \ref{sc7} are void here.
In particular, here in this Appendix, symbols such as $\zeps$, $a_n$,
and $\zeta_n$ will be employed with meanings different from
those (of the same symbols) in Sections \ref{sc3}-\ref{sc7}.
\medskip

     Here in Section \ref{sc8}, all Markov chains will be strictly stationary 
and countable-state.
Unless explicitly stated, the state space need not be a subset of $\R$.
The state space can always simply be coded as, say, the set $\N$ of positive integers.
In any case (with or without such a coding), for such a strictly stationary, countable-state 
Markov chain $\xi := (\xi_k, k \in \Z)$, one can simply take, for each $k \in \Z$, the
$\sigma$-field $\sigma(\xi_k)$ to be the smallest $\sigma$-field containing as
members (elements) the events $\{\xi_k = s\}$ for elements $s$ of the state space;
and then for any nonempty subset $T \subset \Z$, one takes the $\sigma$-field
$\sigma(\xi_k, k \in T)$ to be $\bigvee_{k \in T} \sigma(\xi_k)$.
Then the equations for dependence coefficients for strong mixing conditions, 
such as (\ref{eq2.31})-(\ref{eq2.33}) and (\ref{eq2.41})-(\ref{eq2.43}), as well as
(\ref{eq8.14})-(\ref{eq8.15}) below, simply apply to this context.
\medskip

     In Remark \ref{rem2.7} in Section \ref{sc2}, three counterexamples from the literature 
were cited, as a review of the fact that if a given strictly stationary, non-reversible 
Markov chain satisfies geometric ergodicity (equivalently, absolute regularity with 
mixing rate $\beta(n) \to 0$ at least exponentially fast), and one considers a 
(nontrivial) real function of that Markov chain for which the second moment is finite,
the central limit theorem may fail to hold.
The purpose of this Appendix is to illustrate the fact that in at least one of those
sequences (the one cited from 
\cite {ref-journal-Bradley1983} and \cite {ref-journal-Bradley2007}), 
the ``function'' can be taken to be one-to-one,
with the resulting counterexample to the CLT being a Markov chain itself.
\medskip

   In this Appendix, the word ``countable'' will mean ``countably infinite''.
Both terms will be used freely.
\medskip
 
   In order to set up the discussions/arguments in this Appendix,  
we start with some preliminary background information.

\begin{background}
\label{bckg8.1}
Recall the underlying probability space 
$(\Omega, \cF, P)$ for this paper.
For any two $\sigma$-fields $\cA$ and $\cB$ ($\subset \cF$),
define the following classic coefficient of information:

\begin{equation}
\label{eq8.11}
I(\cA, \cB)\ :=\ \sup \sum_{i=1}^M \sum_{j=1}^N  
\biggl[ P(A_i \cap B_j) \cdot \log 
\biggl( \frac {P(A_i \cap B_j)} {P(A_i)P(B_j)} \biggl) \biggl]
\end{equation} 
where the supremum is taken over all pairs of (finite)
partitions $(A_1, A_2, \dots, A_M)$ and 
$(B_1, B_2, \dots, B_N)$ of $\Omega$
such that $A_i \in \cA$ for each $i$ and $B_j \in \cB$ for each $j$. 
Here interpret $0/0 := 0$ and $0 \log 0 := 0$.
By an application of Jensen's inequality, the quantity $I(\cA, \cB)$
takes its values in the set $[0, \infty] :=  [0, \infty) \cup \{\infty\}$
(it can be $\infty$).
If the $\sigma$-fields $\cA$ and $\cB$ are independent, then
$I(\cA, \cB) = 0$; otherwise, $I(\cA, \cB) > 0$.
\smallskip 

   In fact, for any two $\sigma$-fields $\cA$ and $\cB$,   
\begin{equation}
\label{eq8.12}
\beta(\cA, \cB)\ \leq\ [I(\cA, \cB)]^{1/2}.  
\end{equation}
(Here, if necessary, interpret $\infty^{1/2} = \infty$.)\ \ 
See e.g.\ [\cite {ref-journal-Bradley2007}, Vol.\ 1, Theorem 5.3(III)].  
A version of this inequality was given by Volkonskii and Rozanov 
\cite {ref-journal-VolkRoz}, 
and was at least in spirit attributed there to Pinsker.  
\medskip

   For a given purely atomic $\sigma$-field $\cG$ (with either finitely many or countably
many atoms), the ``entropy of $\cG$'' is the number in $[0, \infty]$ (it can be $\infty$)
defined by 
\begin{equation}
\label{eq8.13}
H(\cG)\ :=\ I(\cG, \cG)\ =\ -\sum_i P(G_i) \log P(G_i) 
\end{equation}
where $G_1, G_2, G_3, \dots$ are the atoms of $\cG$.
(The latter equality in (\ref{eq8.13}) is easy to verify.)
\medskip 

     Now suppose $X := (X_k, k \in \Z)$ is a (not necessarily Markovian) 
strictly stationary sequence of random variables.
For each positive integer $n$, define the following dependence coefficient:
\begin{equation}
\label{eq8.14}
I(n)\ =\ I_X(n)\ :=\ I(\cF_{-\infty}^0, \cF_n^\infty).
\end{equation}
Here we use the notations from Notations \ref{nt2.1} and Notations \ref{nt2.3}.
Analogs of the first two sentences after eq.\ (2.8) hold for these
dependence coefficients.
The strictly stationary random sequence $X$ is said to satisfy 
``information regularity'' if $I(n) \to 0$ as $n \to \infty$.
The information regularity condition was first studied by 
Volkonskii and Rozanov \cite {ref-journal-VolkRoz}, 
and was attributed there to Pinsker.
\smallskip

     For the case of a strictly stationary Markov chain (reversible or not), 
eq.\ (\ref{eq8.14}) takes the following well known augmented form, analogous to 
(\ref{eq2.41})-(\ref{eq2.43}):
\begin{equation}
\label{eq8.15}
I(n)\ =\ I_X(n)\ :=\ I(\cF_{-\infty}^0, \cF_n^\infty)\ =\ I(\sigma(X_0), \sigma(X_n)).
\end{equation}
(See e.g.\ [\cite {ref-journal-Bradley2007}, Vol.\ 1, Theorem 7.3(h)].)\ \  
\end{background}

     In [\cite {ref-journal-Bradley2007}, Vol. 3, Corollary 31.5], 
the following statement was given:
\zhfb

     {\bf Theorem A.}  {\it Suppose $(g_1, g_2, g_3, \dots)$ and $(h_1, h_2, h_3, \dots)$
are each a nonincreasing sequence of positive numbers, such that
$g_n \to 0$ and $h_n \to 0$ as $n \to \infty$.

   Then there exists a strictly stationary Markov chain $\zeta := (\zeta_k, k \in \Z)$
with a countably infinite state space $\Gamma$, and a function 
$f : \Gamma \to \R$, such that, defining the (real) random sequence
$X := (X_k, k \in \Z)$ by $X_k := f(\zeta_k)$ for $k \in \Z$, one has 
that the sequence  $X$ is strictly stationary, and the following statements hold:

     (a) For each $n \in \N$, 
$I_X(n) \leq I_\zeta(n) \leq g_n$ (and hence
$\alpha_\zeta(n) \leq \beta_\zeta(n) \leq g_n^{1/2}$ and
$\alpha_X(n) \leq \beta_X(n) \leq g_n^{1/2}$).

     (b) $H(\sigma(X_0)) \leq H(\sigma(\zeta_0)) \leq g_1$.
That is, the entropy of the marginal distribution of $\zeta_0$, 
and hence also that of the marginal distribution of $X_0$, 
is bounded above by $g_1$.

     (c) $E(X_0^2) < \infty$ and $EX_0 = 0$.
      
     (d)  $\var(X_1 + X_2 + \dots + X_n) \geq n^2 h_n$ for all $n \in \N$.

     (e) Eq.\ (3.3) (in Section 3 of this paper here) holds for the sequence $X$.
     
     (f) For every infinitely divisible law $\mu$ (on $(\R, \cR)$), there exists a
strictly increasing sequence \zhfb 
$(n(1), n(2), n(3), \dots)$ of positive integers, and
sequences of real numbers $(a_1, a_2, a_3, \dots)$ and $(b_1, b_2, b_3, \dots)$,
with $b_j \to \infty$ as $j \to \infty$, such that  
$[(X_1 + X_2 + \dots + X_{n(j)}) - a_j]/b_j \to \mu$ in distribution as $j \to \infty$.}
\zhfb

          This theorem, without mention of properties (b) and (e), is due to
the author [\cite {ref-journal-Bradley1983}, Theorem 2 and p.\ 95, Remark 2.1].
Property (e) was pointed out by the author in 
[\cite{ref-journal-Bradley1997}, p.\ 545].
Property (b) was incorporated into this theorem by the 
author [\cite{ref-journal-Bradley2007}, Vol.\ 3, Corollary 31.5]. 
\medskip

     Under the particular conditions in this theorem, it is easy to see that 
with a little work, the function $f$ can in fact be made to be one-to-one
(without affecting the other statements),
so that the sequence $X$ is itself a Markov chain. 
There are different ways of showing that.
Below, we shall illustrate one particular way of showing that, with a 
method that is somewhat long and indirect but quite gentle and 
unsophisticated, not getting much into technical complications in the 
original construction.
\zhfb

   But first, here is accordingly such a revised version of Theorem A.
\zhfb

     {\bf Theorem B.}\ \  {\it Suppose $(g_1, g_2, g_3, \dots)$ and $(h_1, h_2, h_3, \dots)$
are each a nonincreasing sequence of positive numbers, such that
$g_n \to 0$ and $h_n \to 0$ as $n \to \infty$.

   Then there exists a real, strictly stationary Markov chain 
$Y := (Y_k, k \in \Z)$ with a countably infinite state space ($\subset \R$),
such that the following statements hold:

     (a) For each $n \in \N$, 
$I_Y(n) \leq g_n$ (and hence
$\alpha_Y(n) \leq \beta_Y(n) \leq g_n^{1/2})$.

     (b) $H(\sigma(Y_0)) \leq g_1$.
That is, the entropy of the marginal distribution of $Y_0$ 
is bounded above by $g_1$.

     (c) $E(Y_0^2) < \infty$ and $EY_0 = 0$.
      
     (d)  $\var(Y_1 + Y_2 + \dots + Y_n) \geq n^2 h_n$ for all $n \in \N$.
     
      (e) Eq.\ (3.3) (in Section 3 of this paper here) holds for the sequence $Y$ ---
that is, with $X_k$ replaced by $Y_k$ for each $k$.
      
      (f) For every infinitely divisible law $\mu$ (on $(\R, \cR)$), there exists a
strictly increasing sequence \zhfb
$(n(1), n(2), n(3), \dots)$ of positive integers, and
sequences of real numbers  
$(a_1, a_2, a_3, \dots)$ and $(b_1, b_2, b_3, \dots)$,
with $b_j \to \infty$ as $j \to \infty$, such that  
$[(Y_1 + Y_2 + \dots + Y_{n(j)}) - a_j]/b_j \to \mu$ in distribution
as $j \to \infty$.}
\zhfb

     Of course with no loss of generality, property (d) can be formulated as follows:
There exists $N \in \N$ such that for all $n \geq N$, 
$\var(Y_1 + Y_2 + \dots + Y_n) \geq n^2 h_n$.
For then by multiplying the $Y_k$'s (and also the normalizing constants 
$a_n$ and $b_n$ in (f)) by a sufficiently large positive scalar, one can
ensure property (d) as stated (``for all $n \in \N$'') without affecting
any of the other stated properties. 
\smallskip

    Also, one can omit the assumption that the sequences 
$(g_1, g_2, g_3, \dots)$ and $(h_1, h_2, h_3, \dots)$ are monotonic 
(while retaining the critical assumption that both sequences converge to 0).
To obtain monotonicity for each sequence, one can, for each $n \in \N$, 
simply replace $g_n$ by $\min\{g_1, g_2, \dots, g_n\}$, and 
replace $h_n$ by $\max\{h_n, h_{n+1}, h_{n+2}, \dots\}$.
\zhfb
   
   Below (as indicated above), we shall illustrate one particular way to derive
Theorem B directly as a corollary of Theorem A.
But first, we shall review some standard elementary material from abstract 
algebra, and then state two technical lemmas, all for use in the proof.

\begin{background}
\label{bckg8.2}
{\bf (Algebraic background).}
In all discussion of fields $F \subset \R$ below,
it will be tacitly understood that the two operations will be
the usual addition and multiplication on $\R$ (restricted to $F$).
   
     Let $\Q$ denote the field of rational numbers.
Of course this field $\Q$ is countable.
 
The discussion here will focus on (at first not necessarily countable) sets $S$ 
(not necessarily fields) such that $\Q \subset S \subset \R$. 
\medskip  
   
   (A1) For any set $S$ such that 
$\Q \subset S \subset \R$, let $\cO(S)$ denote the set of all
real numbers $r$ of the form $r = s+t$, $r = s-t$, $r=s \cdot t$, or 
(if $t \neq 0$) $r = s/t$, where $s$ and $t$ are elements of $S$.       
\smallskip 
 
   (A2) For any set $S$ such that $\Q \subset S \subset \R$,
 one has that $S \subset \cO(S)$ (with possible equality).

\noindent [That is trivial:  If $s \in S$, then (since $\Q \subset S$),
$s = s \cdot 1$ (or $s + 0$) $\in \cO(S)$.]
\smallskip

    (A3)  If $\Q \subset S \subset \R$, then define
the sets $\cO^n(S)$ for $n \in \{0,1,2,\dots\}$ recursively as follows:
$\cO^0(S) := S$, and for any integer $n \geq 1$,
$\cO^n(S) := \cO(\cO^{n-1}(S))$.
[Thus $\cO^1(S) := \cO(S)$.]
\smallskip 
  
     (A4)  If $\Q \subset S \subset \R$, and $m$ and $n$ 
are integers such that $0 \leq m \leq n$, then       
$\cO^m(S) \subset \cO^n(S)$.

\noindent [For fixed $m \geq 0$, this follows by (A2) and (A3) and induction
on $n \geq m$.]
\smallskip

   (A5)  For any set $S$ such that $\Q \subset S \subset \R$,
the set $\bigcup_{n=0}^\infty \cO^n(S)$ is a field.
\smallskip

\noindent [Suppose $s, t \in \bigcup_{n=0}^\infty \cO^n(S)$. 
Let $m$ and $n$ be integers such that $s \in \cO^m(S)$ and 
$t \in \cO^n(S)$. 
Let $\ell := \max\{m,n\}$.
Then by (A4), $s,t \in \cO^\ell(S)$, and hence by 
(A1) and (A3), 
$s+t$, $s-t$, $s \cdot t$ and (if $t\neq 0$) $s/t$ are each an
element of $\cO^{\ell+1}(S)$.
Hence the set $\bigcup_{n=0}^\infty \cO^n(S)$ is closed under addition,
subtraction, multiplication, and division (division by 0 excluded).]
\smallskip

    (A6)  For any set $S$ such that $\Q \subset S \subset \R$,
define the notation ${\bf F}(S) := \bigcup_{n=0}^\infty \cO^n(S)$.
Then ${\bf F}(S)$ is the {\it smallest\/} field that contains $S$ as a subset.
That field ${\bf F}(S)$ will be referred to as the ``field generated by the set $S$''.
\smallskip

[The point is that for {\it any\/} field $G$ such that $S \subset G$, one has by 
(A1) and (A3) and induction on $n$ that $\cO^n(S) \subset G$ for every 
$n \in \{0,1,2,\dots\}$; and hence one has that ${\bf F}(S) \subset G$.]   
\zhfb
 
     (B1)  Suppose $\Q \subset S \subset \R$.  
If the set $S$ is countable, then the set $\cO(S)$ is countable.

\noindent [That just follows from (A1) and the fact that the Cartesian 
product $S \times S$ is countable.]
\smallskip

     (B2)  Suppose $\Q \subset S \subset \R$.
If the set $S$ is countable, then by (A1) and (A3) and induction, one has
that for each $n \geq 0$ , the set $\cO^n(S)$ is countable.
Hence the field ${\bf F}(S)$ is countable by (A6) (and the fact that 
the union of countably many countable sets is countable).
\end{background}

\begin{lemma}
\label{lem8.3}
(I) Suppose $\cA_1$, $\cA_2$, $\cB_1$, and $\cB_2$ are $\sigma$-fields 
(on the probability space $(\Omega, \cF, P)$) such that the $\sigma$-fields 
$\cA_1 \vee \cB_1$ and $\cA_2 \vee \cB_2$ are independent.
Then
\begin{equation}
\label{eq8.31}
 I(\cA_1 \vee \cA_2,\ \cB_1 \vee \cB_2)\ =\ 
I(\cA_1, \cB_1)\ +\ I(\cA_2, \cB_2).
\end{equation}   

     (II) If $\cG$ and $\cH$ are {\it independent\/} $\sigma$-fields
(on $(\Omega, \cF, P)$) that are each purely atomic, then 
$H(\cG \vee \cH) = H(\cG) + H(\cH)$.
\end{lemma} 
 
   For (I), see e.g.\ [[2], Vol.\ 1, Lemma 6.4]. 
Of course (II) follows from (\ref{eq8.31}) applied to the 
$\sigma$-fields $\cA_1 = \cB_1 = \cG$ and $\cA_2 = \cB_2 = \cH$.
Both (I) and (II) are old classic equalities from information theory;
for an early reference, see e.g.\ Pinsker \cite {ref-journal-Pinsker}.

\begin{lemma}
\label{lem8.4}
   Suppose $A$ is a positive number, and $x$ and $y$ are real
numbers such that $|x + y|\ <\ 1$.
Then either $|x| < A+1$ or $|y| > A$.
\end{lemma}

   {\bf Proof.}\ \ This is elementary.
One has that 
$|x| = |(-y) + (x + y)| \leq |-y| + |x+y| = |y| + |x+y|$; and hence
under the hypothesis, $|x| - |y| \leq |x+y| < 1$.
If it were the case that $|x | \geq A +1$ and $|y| \leq A$, 
then one would have $|x| - |y| \geq 1$, contradicting the preceding sentence.
Hence Lemma 8.4 holds. 
\zhfb

   {\bf Proof of Theorem B  (Derivation of Theorem B from Theorem A).}\ \ 
The argument will be divided into fourteen ``steps'', including three ``claims''.
\zhfb

    {\bf Step 1.}\ \ 
Refer to the nonincreasing sequences of positive numbers 
$(g_1, g_2, g_3, \dots)$ and $(h_1, h_2, h_3, \dots)$
in the statement of Theorem B.
Also, refer to the first sentence after Theorem B.
For slight extra tidiness of thought,
reducing finitely many of these numbers $g_n$ and $h_n$ if necessary, 
we shall assume without loss of generality that
\begin{equation}
\label{eq8.P11}
{\rm for\ all}\ n \in \N, \quad 0 < g_n \leq 1 \quad {\rm and} \quad
0 < h_n \leq 1.
\end{equation}
     
   Apply Theorem A with $g_n$ replaced by $g_n/2$ for each $n \in \N$.
Then by Theorem A, there exist (henceforth fixed) a strictly stationary 
random sequence $X := (X_k, k \in \Z)$ and a  
strictly stationary Markov chain $\zeta := (\zeta_k, k \in \Z)$
with a countably infinite state space $\Gamma$, and a function 
$f : \Gamma \to \R$, such that 
\begin{equation}
\label{eq8.P12}
X_k := f(\zeta_k) \quad {\rm for\ all}\ k \in \Z, 
\end{equation}
such that (also) the following six conditions are satisfied:

     (a) For each $n \in \N$, $I_\zeta(n) \leq g_n/2$ and hence
$\alpha_\zeta(n) \leq \beta_\zeta(n) \leq (g_n/2)^{1/2} < g_n^{1/2}$.

     (b) $H(\sigma(\zeta_0)) \leq g_1/2$.  

     (c) $E(X_0^2) < \infty$ and $EX_0 = 0$.
      
     (d)  $\var(X_1 + X_2 + \dots + X_n) \geq n^2h_n$ for all $n \in \N$.

     (e) Eq.\ (3.3) (in Remark 3.1 of this paper here) holds for the sequence $X$.
     
     (f) For every infinitely divisible law $\mu$ (on $(\R, \cR)$), there exist 
an infinite set $T \subset \N$,
a sequence $(a_n,\, n \in T)$ of real numbers, and
a sequence $(b_n,\, n \in T)$ of positive numbers, 
with $b_n \to \infty$ as $n \to \infty,\, n \in T$, such that 
$[(X_1 + X_2 + \dots + X_n) - a_n]/b_n \to \mu$ in distribution
as $n \to \infty,\, n \in T$.
\zhfb

     The formulation of property (f) here is equivalent to that of property (f)
in the statement of Theorem A.
Here we simply employ, for later convenience, a trivial change in
notational conventions.
\zhfb

     {\bf Step 2.}\ \ 
Of course $\lim_{x \to 1} (x \log x) = 0$ and $\lim_{x \to 0+} (x \log x) = 0$.
Hence \hfil\break 
$\lim_{p \to 0+}( [p \log p] + [(1-p) \log (1-p)]) = 0$.
Referring to (\ref{eq8.P11}), let $p$ be a number such that
\begin{equation}
\label{eq8.P21}
p \in (0,1/2)\ {\rm and}\ p\ {\rm is\ rational},\ {\rm and}\ \quad 
[p \log p] + [(1-p) \log (1-p)] \leq g_1/2. 
\end{equation}
Note that by (\ref{eq8.P21}), $p < 1/2 < 1-p$ and thus the 
values $p$ and $1-p$ are distinct.
\medskip 

     Let $\eta := (\eta_k, k \in \Z)$ be a sequence of independent, identically
distributed random variables,
with this sequence $\eta$ being independent of the Markov chain $\zeta$
(from the sentence containing (\ref{eq8.P12})), 
such that the random variables $\eta_k,\ k \in \Z$ take only the two 
values $-p$ and $1-p$, with the probability function given by
\begin{equation}
\label{eq8.P22}
P(\eta_0 = -p)\ =\ 1-p \indent {\rm and} \indent P(\eta_0 = 1-p)\ =\ p.
\end{equation}
Then by a simple calculation
\begin{equation}
\label{eq8.P23}
E\eta_0\ =\ 0 \indent {\rm and} \indent E(\eta_0^2)\ =\ \var(\eta_0)\ =\ p(1-p)\ < 1/4.
\end{equation}    
where the last inequality holds by e.g.\ (\ref{eq8.P21}) (its ``first third'') and a simple
calculation.
\smallskip
 
    Referring again to (\ref{eq8.P22}) and the sentence after (\ref{eq8.P21}), 
and also to the last ``third'' of (\ref{eq8.P21}), 
note that the entropy of the random variable $\eta_0$ satisfies
\begin{equation}
\label{eq8.P24}
H(\sigma(\eta_0))\ =\ [(1-p) \log (1-p)]\ +\ [p \log p]\ \leq\ g_1/2
\end{equation}            
By Lemma 8.3(II), eq.\ (\ref{eq8.P24}), and condition (b) 
after (\ref{eq8.P12}), one has that
\begin{equation}
\label{eq8.P25}
H(\sigma(\zeta_0, \eta_0))\, =\, H(\sigma(\zeta_0) \vee \sigma(\eta_0))\,
=\, H(\sigma(\zeta_0))\, +\, H(\sigma(\eta_0))\  
\leq\ g_1/2\, +\, g_1/2\ \leq\ g_1. \quad
\end{equation}

    Refer again to the sentence containing (\ref{eq8.P12}) and the 
entire sentence of (\ref{eq8.P22}).
By an elementary argument using the information in those two sentences, one
has that the sequence $((\zeta_k, \eta_k),\ k \in \Z)$ of random vectors is a 
strictly stationary Markov chain with countable state space 
$\Gamma \times \{-p,\, 1-p\}$. 
\zhfb

     {\bf Step 3.}\ \
Referring to the set $\Gamma$ and the function $f$ in the sentence of
eq.\ (\ref{eq8.P12}), let $\Gamma^*$ denote the range of $f$.
That is, $\Gamma^*$ is the set of all real numbers of the form $f(x)$
where $x \in \Gamma$.
The set $\Gamma^*$ is either finite or countably infinite.
(The set $\Gamma$ is countably infinite, but the function $f$ 
is not necessarily one-to-one.)  
\smallskip  
 
   In the rest of this step, we shall define a sequence $F_0, F_1, F_2, \dots$
of countable fields $\subset \R$, and a sequence $(u_1, u_2, u_3, \dots)$
of real numbers in the open unit interval $(0,1)$.
The definition will be recursive, and runs as follows:
\smallskip

   Refer to the second paragraph above (the first paragraph of Step 3 here).  
The set $\Q \cup \Gamma^*$ (where $\Q$ is the set of  
rational numbers) is countable.
Define the field $F_0 := {\bf F}(\Q \cup \Gamma^*)$, the smallest field containing
(as a subset) the set $\Q \cup \Gamma^*$.
By (B2) in Background \ref{bckg8.2}, this field $F_0$ is countable.
\smallskip

   Now for the recursion step, suppose $n$ is a positive number and the countable
field $F_{n-1} \subset \R$ has already been defined.
Of course the open unit interval $(0,1)$ is {\it un\/}countable,
and hence the set $(0,1) - F_{n-1}$ is nonempty.  
Let $u_n$ be an element of that set $(0,1) - F_{n-1}$.
Define the field $F_n := {\bf F}(F_{n-1} \cup \{u_n\})$, the smallest field
containing (as a subset) the set $F_{n-1} \cup \{u_n\}$.
By (B2) in Background \ref{bckg8.2}, this field $F_n$ is countable. 
[And of course $u_n \in (0,1)$.]
\smallskip

   That completes the recursive definition of the fields $F_n$ for $n \geq 0$
and the numbers $u_n \in (0,1)$ for $n \geq 1$.
\medskip

    One has that
\begin{equation}
\label{eq8.P31}
\Q \cup \Gamma^*\ \subset\ F_0\ \subset\ F_1\ \subset\ F_2\ \subset\
F_3\ \subset\ \dots.  
\end{equation}
Also, one has that
\begin{equation}
\label{eq8.P32}
{\rm for\ each}\ n \in \N, \quad u_n \in (0,1) \cap (F_n - F_{n-1}).
\end{equation}
It follows that these numbers $u_1,\ u_2,\ u_3, \dots$ are distinct. \zhfb
[Why?  If $1 \leq m < n$ are integers, then
$u_m \in F_m$ by (\ref{eq8.P32}) and hence $u_m \in F_{n-1}$ by (\ref{eq8.P31}),
but $u_n \notin F_{n-1}$ by (\ref{eq8.P32}), and hence $u_n \neq u_m$.]
\medskip 

   {\bf Step 4.}\ \
Recall again from the sentence containing (\ref{eq8.P12}) that 
the set $\Gamma$ is countable.
Let that set be represented as
\begin{equation}
\label{eq8.P41} 
\Gamma\ =\ \{\gamma_1, \gamma_2, \gamma_3, \dots\}   
\end{equation}
where each element of $\Gamma$ is listed exactly once in the right hand side.
\smallskip 

     Referring to (\ref{eq8.P41}), define the function 
${\bf v}: \Gamma \to \{u_1, u_2, u_3, \dots\}$ as follows:
\begin{equation}
\label{eq8.P42}
{\rm For\ each}\ n \in \N,\ \ \ {\bf v}(\gamma_n)\ :=\ u_n.
\end{equation}
Refer to the sentence after (\ref{eq8.P32}). 
This function ${\bf v}$ gives a one-to-one correspondence between the 
sets $\Gamma$ and $\{u_1, u_2, u_3, \dots\}$.
By (\ref{eq8.P41}), (\ref{eq8.P42}), and (\ref{eq8.P32}), one has that
\begin{equation}
\label{eq8.P43}
{\rm for\ every}\ t \in \Gamma,\ \ \ {\bf v}{(t)} \in (0,1).
\end{equation}
    
Recall again from the sentence containing (\ref{eq8.P12}) the function 
$f: \Gamma \to \R$. 
Recall the paragraph after (\ref{eq8.P25}).
Referring to (\ref{eq8.P41})-(\ref{eq8.P42}), define the function 
$\tau: \Gamma \times \{-p, 1-p\} \to \R$ as follows:
\begin{equation}
\label{eq8.P44}
{\rm For\ every}\ (t,z) \in \Gamma \times \{-p, 1-p\}, \quad
\tau \Bigl( (t,z) \Bigl)\ :=\ f(t) + {\bf v}(t) \cdot z.
\end{equation} 
That is [see (\ref{eq8.P41}) and (\ref{eq8.P42})], for each $n \in \N$, one has that \zhfb
$\tau((\gamma_n, -p)) = f(\gamma_n) + {\bf v}(\gamma_n) \cdot (-p) 
= f(\gamma_n) + u_n \cdot (-p)$ and \zhfb
$\tau((\gamma_n, 1-p)) = f(\gamma_n) + {\bf v}(\gamma_n) \cdot (1-p)   
= f(\gamma_n) + u_n \cdot (1-p)$.   
\medskip

     {\bf Claim 5.}\ \ 
{\it The function $\tau$ in (\ref{eq8.P44}) is one-to-one.}
\medskip

    {\bf Proof.}\ \ 
Suppose $(t,z)$ and $(t', z')$ are each an element of $\Gamma \times \{-p, 1-p\}$,
and that $(t,z) \neq (t',z')$. 
Our task is to show that
\begin{equation}
\label{eq8.P51}
\tau \Bigl( (t,z) \Bigl)\ \neq\ \tau \Bigl( (t',z') \Bigl).
\end{equation}
The argument will be divided into two cases, according to whether
$t = t'$ (and $z \neq z'$) or $t \neq  t'$ (and either $z = z'$ or $z \neq z'$).
\medskip

   {\bf Case 1:}\ \ $t = t'$ (and $z \neq z'$).
With the equality $t = t'$, the two distinct elements of 
$\Gamma \times \{-p, 1-p\}$ here are $(t,z)$ and $(t,z')$.       
Since $z$ and $z'$ are distinct elements of the set $\{-p, 1-p\}$, 
which has just two elements, the elements $z$ and $z'$ are,
in either order, $-p$ and $1-p$.
Thus the two elements of $\Gamma \times \{-p, 1-p\}$ here are,
in either order, $(t, -p)$ and $(t, 1-p)$.
\smallskip

    By (\ref{eq8.P44}), 
$\tau((t, -p)) = f(t) - {\bf v}(t) \cdot p$ and $\tau((t, 1-p)) = f(t) + {\bf v}(t) \cdot (1-p)$;
and thus by (\ref{eq8.P43}), $\tau((t,-p)) - \tau((t, 1-p)) = -{\bf v}(t) \neq 0$.
Thus $\tau((t, -p)) \neq \tau((t, 1-p))$.
Thus (\ref{eq8.P51}) holds for Case 1.
\medskip

   {\bf Case 2.}\ \ $t \neq t'$ (and either $z= z'$ or $z \neq z'$). 
Referring to the representation (\ref{eq8.P41}) of the 
set $\Gamma$, let $m$ and $n$ be the (distinct) positive 
integers such that $t = \gamma_m$ and
$t' = \gamma_n$.
Without loss of generality, assume that 
\begin{equation}
\label{eq8.P52} 
m\ <\ n.   
\end{equation}

    Then by (\ref{eq8.P44}) and (\ref{eq8.P42}) 
(regardless of whether $z = p$ or $z = 1-p$), one has that
\begin{equation}
\label{eq8.P53} 
\tau \Bigl( (t,z) \Bigl)\ =\ \tau \Bigl( (\gamma_m, z) \Bigl)\
=\ f(\gamma_m) + {\bf v}(\gamma_m) \cdot z\  
=\ f(\gamma_m) + u_m z.
\end{equation}
Similarly (regardless of whether $z' = p$ or $z' = 1-p$)
\begin{equation}
\label{eq8.P54}
\tau \Bigl( (t',z') \Bigl)\ =\ \tau \Bigl( (\gamma_n, z') \Bigl)\
=\ f(\gamma_n) + {\bf v}(\gamma_n) \cdot z'\   
=\ f(\gamma_n) + u_n z'.
\end{equation} 

   Now recall that the (not necessarily distinct) numbers $f(\gamma_m)$ and 
$f(\gamma_n)$ are elements of the set $\Gamma^*$
(from the first sentence of Step 3) and 
hence [see (\ref{eq8.P31})] are elements of the field $F_0$.
Also, recall from (\ref{eq8.P21}) that the number $p$ is rational, and hence 
$-p$ and $1-p$, are rational and are therefor elements of the field $F_0$
[see (\ref{eq8.P31}) again].
Thus the numbers $z$ and $z'$ (distinct or not), each being either $-p$ or $1-p$,
are elements of $F_0$.
Compiling the information so far in this paragraph, and referring to 
(\ref{eq8.P31}) and (\ref{eq8.P52}), one has (whether $z = z'$ or $z \neq z'$) 
that 
\begin{equation}
\label{eq8.P55}
\{ f(\gamma_m),\ f(\gamma_n),\ z,\ z'\}\ \subset\ F_0\ \subset\ F_m\ \subset\ F_{n-1}. 
\end{equation}

    Recall also from (\ref{eq8.P32}) that $u_m \in F_m$.
One now has [see (\ref{eq8.P55})] that each of the 
numbers $f(\gamma_m)$, $z$, and $u_m$ is an element of the field $F_m$.
Hence by (\ref{eq8.P53}), the number $\tau ( (t,z) )$ is an element of the field $F_m$,
and hence [see (\ref{eq8.P55}) again] an element of the field $F_{n-1}$.
\smallskip

   From (\ref{eq8.P55}), we already have that the numbers $f(\gamma_n)$ and $z'$ 
   are each an element of the field $F_{n-1}$.
However, from (\ref{eq8.P32}), one has that the number $u_n$ {\it fails\/}
to be an element of the field $F_{n-1}$.
It follows from (\ref{eq8.P54}) that the number $\tau( (t', z') )$ {\it fails\/} to be 
an element of the field $F_{n-1}$.
\smallskip
        
    [Why?  If it were the case that $\tau( (t', z') ) \in F_{n-1}$, 
then [since $f(t') = f(\gamma_n) \in F_{n-1}$] one would have by 
(\ref{eq8.P54}) that $u_n z' \in F_{n-1}$ and hence [since $z' \in F_{n-1}$ and 
$z'$ is either $-p$ or $1-p$ and $p \in (0,1/2)$ by (\ref{eq8.P21}) and 
therefore in either case $z' \neq 0$] one would have that $u_n \in F_{n-1}$, 
contradicting the second sentence of the preceding paragraph.]
\smallskip
 
   To repeat the last sentence of each of the second and third paragraphs above
this one, $\tau( (t,z) ) \in F_{n-1}$ and $\tau ( (t',z') ) \notin F_{n-1}$.
Hence (\ref{eq8.P51}) holds for Case 2.
That completes the proof of Claim 5.
\medskip

   {\bf Step 6.}\ \
By the paragraph of eq.\ (\ref{eq8.P12}), and the second paragraph of Step 2,
and eq.\ (\ref{eq8.P43}), for each $k \in \Z$, the random variables
$X_k$, $\eta_k$, and ${\bf v}(\zeta_k)$ are each real-valued. 
Refer to the paragraph after (\ref{eq8.P25}). 
On our probability space $(\Omega, \cF, P)$, 
define the sequence $Y := (Y_k, k \in \Z)$ of real-valued random variables
as follows:  For each $k \in \Z$ and each $\omega \in \Omega$
[see (\ref{eq8.P44}) and again (\ref{eq8.P12})], 
\begin{align}
\label{eq8.P61}
Y_k(\omega)\ :=\ \tau \Bigl( \Bigl(\zeta_k(\omega), \eta_k(\omega) \Bigl) \Bigl)\
&=\ f\Bigl( \zeta_k(\omega) \Bigl)\  
+\ {\bf v}\Bigl( \zeta_k(\omega) \Bigl)\, \cdot\, \eta_k(\omega) \nonumber\\
&=\ X_k(\omega)\ +\ {\bf v}\Bigl( \zeta_k(\omega) \Bigl)\, \cdot\, \eta_k(\omega). 
\end{align}

By the paragraph after (\ref{eq8.P25}), this random sequence $Y$ is strictly stationary.
By the same paragraph and Claim 5, this sequence $Y$ is also a Markov chain,
with a countably infinite state space.
That is, $Y$ is a real, strictly stationary, countable-state Markov chain.
\smallskip

     To complete the derivation of Theorem B from Theorem A, what remains
is to verify the six properties (a)-(f) in the statement of Theorem B.
But first, in connection with (\ref{eq8.P61}), we shall lay out in Claim 7 below 
some preliminary information involving products of certain random variables:
\zhfb

   {\bf Claim 7.}\ \ {\it Refer to item (c) after (\ref{eq8.P12}).  
The following statements hold:

     (a) For each $\omega \in \Omega$ and each $k \in \Z$, one has that
${\bf v}(\zeta_k(\omega)) \in (0,1)$ and $|\eta_k(\omega)| < 1$ and 
(hence) $|{\bf v}(\zeta_k(\omega)) \cdot \eta_k(\omega)| < 1$.

     (b) For each $k \in \Z$, one has that $E[{\bf v}(\zeta_k) \cdot \eta_k] = 0$.
     
     (c) For any two distinct integers $j$ and $k$, one has that
$E( [{\bf v}(\zeta_j) \cdot \eta_j] \cdot [{\bf v}(\zeta_k) \cdot \eta_k]) = 0.$

     (d) If $j$ and $k$ are any integers (not necessarily distinct), then
$E(X_j \cdot [{\bf v}(\zeta_k) \cdot \eta_k]) = 0$.

      (e) For any positive integer $n$, one has that
$E[(\sum_{k=1}^n [{\bf v}(\zeta_k) \cdot \eta_k])^2] \leq n$. 

      (f) For any positive integer $n$,
$E[(\sum_{j=1}^n X_j) (\sum_{k=1}^n [{\bf v}(\zeta_k) \cdot \eta_k])] = 0$.} 
\medskip

   {\bf Proof of Claim 7.}\ \ 
Statement (a) holds by (\ref{eq8.P43}) and (\ref{eq8.P21}) and 
the entire sentence of (\ref{eq8.P22}).

     For statements (b), (c), and (d), recall the entire sentence of (\ref{eq8.P22}),
including the independence conditions there and the fact [see (\ref{eq8.P23})]
that $E\eta_k = 0$ for each integer $k$. 

   Statement (b) holds because
$E[{\bf v}(\zeta_k) \cdot \eta_k] =  E[{\bf v}(\zeta_k)] \cdot E\eta_k
=  E[{\bf v}(\zeta_k)] \cdot 0 = 0$.

   Statement (c) holds because in its asserted equality, the left side equals
$E([{\bf v}(\zeta_j)] \cdot [{\bf v}(\zeta_k)]) \cdot E\eta_j \cdot E\eta_k
= E([{\bf v}(\zeta_j)] \cdot [{\bf v}(\zeta_k)]) \cdot 0 \cdot 0$. 
           
   For statement (d), recall from (\ref{eq8.P12}) that $X_j = f(\zeta_j)$,
and refer to statement (a) again [as well as, again, item (c) after (\ref{eq8.P12})].  
One has that
$$
E(X_j \cdot [{\bf v}(\zeta_k) \cdot \eta_k])\, 
=\, E[f(\zeta_j) \cdot {\bf v}(\zeta_k) \cdot \eta_k]\,  
=\, E[f(\zeta_j) \cdot {\bf v}(\zeta_k)] \cdot E\eta_k\, 
=\, E[f(\zeta_j) \cdot {\bf v}(\zeta_k)] \cdot 0\,
=\, 0.
$$

   By statement (a), $E([{\bf v}(\zeta_k) \cdot \eta_k]^2) \leq 1$ for each $k \in \Z$.
Hence statement (e) follows from statement (c).
Similarly statement (f) follows from statement (d).
That completes the proof of Claim 7.
\zhfb

   {\bf Step 8.  Proof of statement (a) in Theorem B.}\ \ 
Suppose $n \in \N$.
Then $I_Y(n) = I(\sigma(Y_0), \sigma (Y_n))$ by
(\ref{eq8.15}) together with (say) the third sentence 
after eq.\ (\ref{eq8.P61}).
Also,
 $\sigma(Y_k) = \sigma(\zeta_k, \eta_k) 
 = \sigma(\zeta_k) \vee \sigma(\eta_k)$
for each $k \in \Z$ by (\ref{eq8.P61}) and Claim 5.
Hence by Lemma \ref{lem8.3}(I) and the sentence of eq.\ (\ref{eq8.P22}), 
then eq.\ (\ref{eq8.15}) [together with the sentence of (\ref{eq8.P12})] 
and then item (a) after (\ref{eq8.P12}),
\begin{align}
\nonumber
I_Y(n)\ &=\ I\bigl(\sigma(Y_0),\, \sigma (Y_n) \bigl)\
=\ I\bigl( \sigma(\zeta_0) \vee \sigma (\eta_0),\ \sigma(\zeta_n) \vee \sigma(\eta_n) \bigl) \\
&=\ I\bigl(\sigma(\zeta_0), \sigma(\zeta_n) \bigl)\ 
+\ I(\bigl(\sigma(\eta_0), \sigma(\eta_n)\bigl)\
=\ I_\zeta(n) + 0\ <\ g_n.                    
\end{align}
Thus $I_Y(n) \leq g_n$.
Hence also $\alpha_Y(n) \leq \beta_Y(n) 
= \beta(\sigma(Y_0), \sigma(Y_n)) \leq g_n^{1/2}$ 
by (\ref{eq2.34}), (\ref{eq2.42}) [together with the third sentence 
after (\ref{eq8.P61})], and (8.2). 
Since $n \in \N$ was arbitrary, statement (a) in Theorem B holds.
\medskip

     {\bf Step 9.  Proof of statement (b) in Theorem B.}\ \   
By (\ref{eq8.P61}) and Claim 5 (again),
$\sigma(Y_0) = \sigma(\zeta_0, \eta_0)$.
Hence $H(\sigma(Y_0)) = H(\sigma(\zeta_0, \eta_0)) \leq g_1$ 
by (\ref{eq8.P25}).
Thus statement (b) in Theorem B holds.
\medskip

   {\bf Step 10.   Proof of statement (c) in Theorem B.}\ \     
From (\ref{eq8.P61}), one has that $Y_0 = X_0 + {\bf v}(\zeta_0) \cdot \eta_0$.
This has two consequences: 
First, $Y_0$ is the sum of two square-integrable random variables,
by Claim 7(a) together with item (c) after (\ref{eq8.P12}); and hence $E(Y_0^2) < \infty$.
Second, $EY_0 = 0$ by Claim 7(b) and (again) item (c) after (\ref{eq8.P12}).
Thus statement (c) in Theorem B holds. 
\medskip

     {\bf Step 11.   Proof of statement (d) in Theorem B.}\ \ 
Suppose $n$ is any positive integer.

   By strict stationarity [recall the paragraph after (\ref{eq8.P61})] 
and the equality $EY_0 = 0$ proved in Step 10 above,
one has that 
$E(Y_1 + Y_2 + \dots + Y_n) = 0$  and hence
$\var(Y_1 + Y_2 + \dots + Y_n) = E[(Y_1 + Y_2 + \dots + Y_n)^2]$. 
    
   Refer again to (\ref{eq8.P61}) for the representation     
$Y_k = X_k + {\bf v}(\zeta_k) \cdot \eta_k$ for $k \in \Z$.
One has that
$(Y_1 + Y_2 + \dots+ Y_n)  = 
(\sum_{j=1}^n X_j) + \sum_{k=1}^n ({\bf v}(\zeta_k) \cdot \eta_k)$ 
and hence
\begin{equation}
\nonumber
(Y_1 + Y_2 + \dots + Y_n)^2\ 
=\ \biggl(\, \sum_{j=1}^n X_j \biggl)^2\ 
+\ \biggl(\, \sum_{k=1}^n ({\bf v}(\zeta_k) \cdot \eta_k) \biggl)^2  
+\ 2 \cdot \biggl(\, \sum_{j=1}^n X_j \biggl) \cdot 
 \biggl(\, \sum_{k=1}^n ({\bf v}(\zeta_k) \cdot \eta_k) \biggl).
\end{equation}
As displayed here, the right hand side is the sum of three random variables,
of which (i) the first has expected value    
$E[(\sum_{j=1}^n X_j)^2] \geq n^2h_n$ by item (d) after (\ref{eq8.P12}),
(ii) the second (trivially being nonnegative) has expected value $\geq 0$, and 
(iii) the third has expected value 0 by Claim 7(f). 
Hence $E[(Y_1 + Y_2 + \dots + Y_n)^2] \geq n^2 h_n$.

     Since $n \in \N$ was arbitrary, statement (d) in Theorem B holds.    
\zhfb

    {\bf Step 12.  Proof of statement (e) in Theorem B.}\ \       
Recall the first paragraph after that of eq.\ (\ref{eq3.13}).
It suffices to carry out the argument for $b=1$ 
(with $X_k$ replaced by $Y_k$) in (\ref{eq3.13}).
That is, it suffices to show that the number
$\sup_{r \in \R} P(r-1 <  n^{-1/2} \sum_{k=1}^n Y_k < r+1)$
converges to $0$ as $n \to \infty$.
\smallskip

     Accordingly, let $\zeps > 0$ be arbitrary but fixed.
To complete the proof of statement (e) in Theorem B, it suffices to
prove that there exists $N \in \N$ such that 
\begin{equation}
\label{eq8.P121}
{\rm for\ all}\ n \geq N,\ \ \ 
\sup_{r \in \R} P\biggl(r-1\, <\, n^{-1/2}\sum_{k=1}^n Y_k\, <\, r+1\biggl)\ \leq\ \zeps.         
\end{equation}

   The event(s) in the left side of (\ref{eq8.P121}) can be rewritten as
$\{ |-r + n^{-1/2} \sum_{k=1}^n Y_k| < 1\}$.
With that in mind, for each $n \in \N$ and each $r \in \R$, define the
random variables
\begin{equation}
\label{eq8.P122}
T_{n,r}\ :=\ -r\ +\ n^{-1/2}\sum_{k=1}^n Y_k \indent {\rm and} \indent
U_{n,r} :=\ -r\ +\ n^{-1/2}\sum_{k=1}^n X_k;
\end{equation}
and for each $n \in \N$, define the random variable
\begin{equation}
\label{eq8.P123}
V_n\ :=\ n^{-1/2} \sum_{k=1}^n [{\bf v}(\zeta_k) \cdot \eta_k].
\end{equation}
Then by (\ref{eq8.P61}), for each $n \in \N$ and each $r \in \R$, 
\begin{equation}
\label{eq8.P124}
T_{n, r}\ =\ -r\ +\ n^{-1/2} \sum_{k=1}^n [ X_k + {\bf v}(\zeta_k) \cdot \eta_k]\
=\ U_{n,r} + V_n.
\end{equation} 

   Referring to the positive number $\zeps$ in the ``goal equation'' (\ref{eq8.P121}), 
let $A$ be a positive number sufficiently large that 
\begin{equation}
\label{eq8.P125}
1/A^2\ \leq\ \zeps/2.
\end{equation} 
By Chevyshev's inequality, eq.\ (\ref{eq8.P123}), Claim 7(e), 
and then eq.\ (\ref{eq8.P125}), one has that
\begin{equation}
\label{eq8.P126}
{\rm for\ each}\ n \in \N,\ \ \ 
P(|V_n| \geq A)\ \leq\ \frac {E(V_n^2)} {A^2}\ 
\leq\ \frac {1} {A^2}\ \leq\ \zeps/2.
\end{equation}

   Referring to item (e) right after (\ref{eq8.P12}), one has by (3.3) 
that the number \zhfb
$\sup_{r \in \R} P(r - [A+1] < n^{-1/2} \sum_{k=1}^n X_k < r + [A+1])$      
converges to 0 as $n \to \infty$.
That is, by (\ref{eq8.P122}),
$\lim_{n \to \infty} [\sup_{r \in \R} P(|U_{n,r}| < A+1)] = 0$.   
Accordingly, let $N$ be a positive integer such that
\begin{equation}
\label{eq8.P127}
{\rm for\ all}\ n \geq N,\ \ \ \sup_{r \in \R} P(|U_{n,r}| < A+1)\ \leq\ \zeps/2.
\end{equation} 
Our goal now is to prove that (\ref{eq8.P121}) holds [with the same 
integer $N$ there as in (\ref{eq8.P127}) here].
\smallskip 
 
    By (\ref{eq8.P124}) and Lemma 8.4, for each $n \in N$ and each $r \in \R$,
\begin{equation}
\nonumber
\{ |T_{n,r}|\ < 1 \}\ \subset\ \{ |U_{n,r}| < A+1 \} \cup \{ |V_n| > A \}.
\end{equation}     
Hence by (\ref{eq8.P126}) and (\ref{eq8.P127}), for each $n \geq N$ 
and each $r \in \R$, one has that
\begin{equation}
\nonumber  
P(|T_{n,r}| < 1)\ \leq\ P( |U_{n,r}| < A+1)\ +\ P(|V_n| > A)\ \leq\ 
(\zeps/2) + (\zeps/2)\ =\ \zeps. 
\end{equation}  
That is, by (\ref{eq8.P122}), eq.\ (\ref{eq8.P121}) holds.
That completes the proof of statement (e) in Theorem B.     
\zhfb

   {\bf Claim 13.}\ \ 
{\it Refer to item (f) after (\ref{eq8.P12}).
Suppose $\mu$ is a non-degenerate infinitely divisible law (on $(\R, \cR)$),
$G$ is an infinite subset of $\N$, 
and $a_n, n \in G$ are real numbers, 
and $b_n, n \in G$ are positive numbers,
such that
\begin{equation}
\label{eq8.P131}
\frac {(X_1 + X_2 + \dots + X_n) - a_n} {b_n}\ \to\ \mu\ \  
{\rm in\ distribution\ as}\ n \to \infty,\ n \in G.
\end{equation}
Then
\begin{equation}
\label{eq8.P132}
\lim_{n \to \infty,\, n \in G} b_n/n^{1/2}\ =\ \infty.
\end{equation}
}

     {\bf Proof of Claim 13.}\ \ 
Suppose the entire sentence of (\ref{eq8.P131}) holds.
Suppose (\ref{eq8.P132}) is false.
We shall seek a contradiction.
\smallskip

   Under the assumption that (\ref{eq8.P132}) is false, there exists an infinite subset
$H \subset G$, and a positive (finite) number $D$ such that
\begin{equation}
\label{eq8.P133}
{\rm for\ all}\ n \in H,\ \ \ b_n/n^{1/2}\ <\ D.
\end{equation}

   Referring to the entire sentence of (\ref{eq8.P131}), let $F_\mu$ denote the 
(cumulative) distribution function of the 
non-degenerate, infinitely divisible law $\mu$.
That is, $F_\mu (x) = \mu((-\infty, x])$ for $x \in \R$.
Let $x_0$ be a real number such that $0 < F_\mu(x_0) < 1$.    
(Such an $x_0$ exists because $\mu$ is non-degenerate.)\ \
Of course the set of continuity points of $F_\mu$ is dense in $\R$
(with the set of {\it dis\/}continuity points of $F_\mu$ being at most
countable).
Let $s$ and $t$ be {\it continuity\/} points of $F_\mu$ such that
\begin{equation}
\label{eq8.P134}
s < x_0 < t \indent {\rm and}\ \indent F_\mu(s) < F_\mu(x_0) < F_\mu(t).
\end{equation} 
(Of course this is redundant; the second condition implies the first.)\ \ 
Then
\begin{equation}
\label{eq8.P135}
\mu\bigl( (s,t) \bigl)\ =\ F_\mu(t) - F_\mu(s)\ >\ 0. 
\end{equation}

   Recall the entire sentence of (\ref{eq8.P133}).
By (\ref{eq8.P131}), the number 
$P(s < [(X_1 + X_2 + \dots + X_n) - a_n]/b_n < t)$ 
converges to $\mu((s,t))$ as $n \to \infty,\ n \in H$
(in fact as $n \to \infty,\ n \in G$).  
That is,
$P(sb_n + a_n < X_1 + X_2 + \dots + X_n < tb_n + a_n)$
converges to $\mu((s,t))$ as $n \to \infty,\ n \in H$.
That is [now also including (\ref{eq8.P135})], 
\begin{equation}
\label{eq8.P136}
\lim_{n \to \infty,\, n \in H} P\biggl(
\frac {sb_n + a_n} {n^{1/2}}
< \frac {X_1 + X_2 + \dots+ X_n} {n^{1/2}} 
< \frac {tb_n + a_n} {n^{1/2}} \biggl)\
=\ \mu \bigl( (s,t) \bigl)\ >\ 0. \quad
\end{equation}

   Referring to (\ref{eq8.P133}) and (\ref{eq8.P134}), define the positive number 
\begin{equation}
\label{eq8.P137}
c\ :=\ \frac {D \cdot (t-s)} {2}.
\end{equation}
Then by (\ref{eq8.P133}), one has that
\begin{equation}
\label{eq8.P138}
{\rm for\ every}\ n \in H, \ \ \ \frac {b_n} {n^{1/2}} \cdot \frac {t-s} {2}\ <\ c.
\end{equation}
For each $n \in H$, define the real number 
\begin{equation}
\label{eq8.P139}
r_n\ :=\ \frac {a_n} {n^{1/2}}\ +\ \frac {b_n} {n^{1/2}} \cdot \frac {s + t} {2}.
\end{equation}
Then for each $n \in H$, by (\ref{eq8.P139}) and (\ref{eq8.P138}),
\begin{equation}
\label{eq8.P1310}
r_n\ - c\ <\ 
\frac {a_n} {n^{1/2}}\ +\ \frac {b_n} {n^{1/2}} \cdot \frac {s + t} {2}\ 
-\ \frac {b_n} {n^{1/2}} \cdot \frac {t-s} {2}\
=\  \frac {a_n} {n^{1/2}}\ +\ \frac {b_n} {n^{1/2}} \cdot s  
\end{equation}
and 
\begin{equation}
\label{eq8.P1311}
r_n\ + c\ >\ 
\frac {a_n} {n^{1/2}}\ +\ \frac {b_n} {n^{1/2}} \cdot \frac {s + t} {2}\ 
+\ \frac {b_n} {n^{1/2}} \cdot \frac {t-s} {2}\
=\  \frac {a_n} {n^{1/2}}\ +\ \frac {b_n} {n^{1/2}} \cdot t. 
\end{equation}
The last terms in (\ref{eq8.P1310}) and (\ref{eq8.P1311}) are respectively 
the left and right ``endpoints'' in the term $P(\dots)$ in (\ref{eq8.P136}).
Hence by (\ref{eq8.P1310}), (\ref{eq8.P1311}), and (\ref{eq8.P136}) itself,
\begin{equation}
\label{eq8.P1312}
\liminf_{n \to \infty,\, n \in H} P\biggl(
r_n - c\
< \frac {X_1 + X_2 + \dots+ X_n} {n^{1/2}} 
<\ r_n + c \biggl)\
\geq\ \mu \bigl( (s,t) \bigl)\ >\ 0. 
\end{equation}

   However, by item (e) after (\ref{eq8.P12}) --- see eq.\ (\ref{eq3.13}) ---
\begin{equation}
\nonumber
\lim_{n \to \infty} P\biggl(
r_n - c\
< \frac {X_1 + X_2 + \dots+ X_n} {n^{1/2}} 
<\ r_n + c \biggl)\
= \ 0. 
\end{equation}   
That contradicts (\ref{eq8.P1312}).
Hence (\ref{eq8.P132}) holds after all.
That completes the proof of Claim 13.
\zhfb

     {\bf Step 14.  Proof of statement (f) in Theorem B.}\ \
Suppose $\mu$ is an infinitely divisible law (on $(\R, \cR)$). 
The proof of (f) will be divided into two cases according to whether
$\mu$ is degenerate (the frivolous case) or not.
\medskip

     {\bf Case 1:  $\mu$ is degenerate.}\ \   
That is, for some real number $r$ (now fixed),
$\mu$ is the point mass at $r$, say $\mu = \nu_r$
where $\nu_r(\{r\}) = 1$ .
\smallskip 

     For each positive integer $n$, let $t_n \geq 1$ be a positive number 
sufficiently large that $P(|Y_1 + Y_2 + \dots + Y_n|  \geq t_n) < 2^{-n}$.
For each positive integer $n$, define the positive number 
$b_n := 2^n t_n \geq 2^n$.
For convenience, for each positive integer $n$, define the random variable
$W_n := (Y_1 + Y_2 + \dots + Y_n)/b_n$.
Then for each positive integer $n$, 
$P(|W_n| \geq 2^{-n}) = P(|Y_1 + Y_2 + \dots + Y_n| \geq t_n) < 2^{-n}$.
\smallskip

     Hence $W_n \to 0$ in probability as $n \to \infty$.
Hence $W_n + r \to r$ in probability as $n \to \infty$.
Hence $W_n + r$ converges in distribution to $\nu_r$ as $n \to \infty$.
That is, $[(Y_1 + Y_2 + \dots. + Y_n) - (-b_nr)]/b_n$ converges 
in distribution to $\nu_r$ ($=\mu$) as $n \to \infty$.
Also (by the definition of $b_n$ above) $b_n \to \infty$ as $n \to \infty$.
Thus statement (f) in Theorem B holds (in a trivial way) for Case 1. 
\medskip

     {\bf Case 2.  The (infinitely divisible) law $\mu$ is non-degenerate.}\ \
 In the presentation of this argument, we shall employ an obvious analog
 of the paragraph at the end of Step 1.
 \smallskip    
     
   Apply item (f) after eq.\ (\ref{eq8.P12}).
Let $T$ be an infinite subset of $\N$,
and let $a_n,\, n \in T$ be real numbers, and let
$b_n,\, n \in T$ be positive numbers 
satisfying $b_n \to \infty$ as $n \to \infty,\, n \in T$, such that
\begin{equation}
\label{eq8.P141}
\frac {(X_1 + X_2 + \dots + X_n) - a_n} {b_n}\  \to \mu\ \  
{\rm in\ distribution\ as}\ n \to \infty,\, n \in T.
\end{equation}
To complete the proof of statement (f) in Theorem B for Case 2, it suffices to
prove that
\begin{equation}
\label{eq8.P142}
\frac {(Y_1 + Y_2 + \dots + Y_n) - a_n} {b_n}\  \to \mu\ \   
{\rm in\ distribution\ as}\ n \to \infty,\ n \in T.
\end{equation} 
 
     For (say) each $n \in T$, by Claim 7(e),
$E [ (\, n^{-1/2} \sum_{k=1}^n [{\bf v}(\zeta_k) \cdot \eta_k]\, )^2 ] \leq 1$.
Hence for any $n \in T$, 
multiplying both sides by $(n^{1/2}/b_n)^2$, one obtains the inequality 
$E [ (\, b_n^{-1} \sum_{k=1}^n [{\bf v}(\zeta_k) \cdot \eta_k]\, )^2 ] \leq (n^{1/2}/b_n)^2$. 
By Claim 13 and the entire sentence of (\ref{eq8.P141}), 
$n^{1/2}/b_n \to 0$ as $n \to \infty,\ n \in T$.
Hence  
$E [ (\, b_n^{-1} \sum_{k=1}^n [{\bf v}(\zeta_k) \cdot \eta_k]\, )^2 ] \to 0$ 
as $n \to \infty,\ n \in T$.  It follows that
\begin{equation}
\label{eq8.P143}
\frac {\sum_{k=1}^n [{\bf v}(\zeta_k) \cdot \eta_k]} {b_n}\ 
\to\ 0\ \ {\rm in\ probability\ as}\ n \to \infty,\, n \in T.
\end{equation}

   By (\ref{eq8.P61}), for each $n \in T$, the sum of the 
left side of (\ref{eq8.P141}) and the left side of (\ref{eq8.P143}) is
the left side of (\ref{eq8.P142}).
Hence by (\ref{eq8.P141}), (\ref{eq8.P143}), and Slutsky's Theorem,
one has that (\ref{eq8.P142}) holds.
That completes the proof of statement (f) in Theorem B for Case 2.
That completes the proof of statement (f) in Theorem B.
That completes the proof of Theorem B.

\end{document}